\input amstex
\documentstyle{amsppt}
\loadbold
\magnification=\magstep1
\pageheight{9.0truein}
\pagewidth{6.5truein}
\NoBlackBoxes
\TagsAsMath

\def\ZZ{{\Bbb Z}}

\def\CC{{\Bbb C}}
\def\AA{{\Bbb A}}
\def\B{{\Cal B}}
\def\C{{\Cal C}}

\def\L{{\Cal L}}

\def\O{{\Cal O}}
\def\P{{\Cal P}}
\def\R{{\Cal R}}
\def\U{{\Cal U}}
\def\SS{{\Cal S}}
\def\afrak{{\frak a}}
\def\bfrak{{\frak b}}
\def\dfrak{{\frak d}}
\def\gfrak{{\frak g}}
\def\hfrak{{\frak h}}
\def\pfrak{{\frak p}}
\def\lfrak{{\frak l}}
\def\nfrak{{\frak n}}

\def\qfrak{{\frak q}}

\def\wbar{\overline{w}}
\def\xbar{\overline{x}}
\def\Xbar{\overline{X}}

\def\til{\widetilde}

\def\ltil{\tilde{l}}

\def\util{\tilde{u}}
\def\vtil{\tilde{v}}

\def\Qtil{\til{Q}}
\def\Stil{\til{S}}
\def\Util{\til{U}}
\def\Vtil{\til{V}}

\def\Thetatil{\til{\Theta}}
\def\bfc{{\boldsymbol c}}
\def\bfe{{\boldsymbol e}}
\def\bfr{{\boldsymbol r}}
\def\bu{\bullet}

\def\leftmatrix{\left[ \smallmatrix \format \c&&\quad\c \\}
\def\rightmatrix{\endsmallmatrix\right]}

\def\leftmat{\left[ \smallmatrix}
\def\rightmat{\endsmallmatrix\right]}

\def\ep{\epsilon}
\def\Mn{M_n}
\def\Mmn{M_{m,n}}
\def\GLm{GL_m}
\def\GLn{GL_n}
\def\GLN{GL_N}
\def\GLt{GL_t}
\def\SLn{SL_n}

\def\CCx{\CC^{\times}}
\def\gl{\frak{gl}}
\def\glm{\frak{gl}_m}
\def\gln{\frak{gl}_n}
\def\glt{\frak{gl}_t}

\def\Omnt{\O^{m,n}_t}

\def\bfc{{\boldsymbol c}}
\def\bfr{{\boldsymbol r}}

\def\onen{\{1,\dots,n\}}

\def\ol{\overline}
\def\wn{w^n_\circ}
\def\wm{w^m_\circ}
\def\wt{w^t_\circ}
\def\wmn{w^{m+n}_\circ}
\def\wmwn{w^{m,n}_\circ}
\def\wnwmpair{(\wn,\wm)}
\def\wN{w^N_\circ}
\def\ci{\circ}
\def\Sigmnt{\Sigma^{m,n}_t}
\def\Smnt{S_{m+n}^{[-n,m]}[t]}
\def\Bww{\B^{w_1,w_2}}
\def\Stilmn{\Stil_{m,n}}

\def\Lie{\operatorname{Lie}}
\def\Cl{\operatorname{Cl}}
\def\rank{\operatorname{rank}}

\def\stab{\operatorname{Stab}}
\def\sign{\operatorname{sign}}
\def\dom{\operatorname{dom}}
\def\rng{\operatorname{rng}}
\def\tr{^{\operatorname{tr}}}

\def\diag{\operatorname{diag}}
\def\Ad{\operatorname{Ad}}
\def\Gr{\operatorname{Gr}}
\def\Inc{\operatorname{Inc}}
\redefine\Im{\operatorname{Im}}

\def\BFZ{{\bf 1}}
\def\Bor{{\bf 2}}
\def\Car{{\bf 3}}
\def\CP{{\bf 4}}
\def\Deo{{\bf 5}}
\def\DriICM{{\bf 6}}
\def\Dri{{\bf 7}}
\def\FZone{{\bf 8}}
\def\FZtwo{{\bf 9}}
\def\Ful{{\bf 10}}
\def\Fulbook{{\bf 11}}
\def\GSV{{\bf 12}}
\def\GLijm{{\bf 13}}
\def\Har{{\bf 14}}
\def\HLthree{{\bf 15}}
\def\HLn{{\bf 16}}
\def\HY{{\bf 17}}
\def\Hum{{\bf 18}}
\def\Jos{{\bf 19}}
\def\Karone{{\bf 20}}
\def\KS{{\bf 21}}
\def\Lau{{\bf 22}}
\def\Lu{{\bf 23}}
\def\LuY{{\bf 24}}
\def\Nor{{\bf 25}}
\def\Ren{{\bf 26}}
\def\Soi{{\bf 27}}
\def\V{{\bf 28}}
\def\W{{\bf 29}}

\topmatter

\title Poisson structures on affine spaces and flag varieties. I. Matrix
affine Poisson space
\endtitle

\rightheadtext{The matrix affine Poisson space}

\author K. A. Brown, K. R. Goodearl, and  
M. Yakimov 
\endauthor

\address Department of Mathematics, University of Glasgow,
Glasgow G12 8QW, Scotland\endaddress
\email kab\@maths.gla.ac.uk\endemail

\address Department of Mathematics, University of California,
Santa Barbara, CA 93106, USA\endaddress
\email goodearl\@math.ucsb.edu\endemail

\address Department of Mathematics, University of California,
Santa Barbara, CA 93106, USA\endaddress
\email yakimov\@math.ucsb.edu\endemail


\thanks The research of the first two authors was partially
supported by Leverhulme Research Interchange Grant F/00158/X,
and was begun during their participation in the Noncommutative Geometry
program at the Mittag-Leffler Institute in Fall 2003. They thank the
Institute for its fine hospitality. The research of the second author
was also partially supported by  National Science Foundation
grant DMS-9970159. The research of the third author was
partially supported by  NSF grant DMS-0406057 and a UCSB
junior faculty research  incentive grant. He
thanks the University  of Hong Kong for the warm hospitality
during a part of the preparation of this paper.
\endthanks

\abstract The standard Poisson structure on the rectangular
matrix variety $\Mmn(\CC)$ is investigated, via the orbits of
symplectic leaves under the action of the maximal torus
$T\subset GL_{m+n}(\CC)$. These orbits, finite in number, are
shown to be smooth irreducible locally closed subvarieties of
$\Mmn(\CC)$, isomorphic to intersections of dual Schubert
cells in the full flag variety of $GL_{m+n}(\CC)$. Three different
presentations of the $T$-orbits of symplectic leaves in
$\Mmn(\CC)$ are obtained -- (a) as pullbacks of Bruhat cells
in $GL_{m+n}(\CC)$ under a particular map; (b) in terms of
rank conditions on rectangular submatrices; and (c) as matrix
products of sets similar to double Bruhat cells in
$\GLm(\CC)$ and $\GLn(\CC)$. In presentation (a), the orbits
of leaves are parametrized by a subset of the Weyl group
$S_{m+n}$, such that inclusions of Zariski closures
correspond to the Bruhat order. Presentation (b) allows
explicit calculations of orbits. From presentation (c) it
follows that, up to Zariski closure, each orbit of leaves is
a matrix product of one orbit with a fixed column-echelon
form and one with a fixed row-echelon form. Finally,
decompositions of generalized double Bruhat cells in
$\Mmn(\CC)$ (with respect to pairs of partial permutation
matrices) into unions of $T$-orbits of symplectic leaves are
obtained.
\endabstract

\endtopmatter
\document

\head Introduction \endhead

\definition{0.1}
We investigate the geometry of the affine variety $\Mmn=\Mmn(\CC)$ of
complex $m\times n$ matrices in relation to its standard Poisson
structure (see \S1.5) and to the action of the torus of ``row and
column automorphisms''. Specifically, let $T$ denote the torus of
diagonal matrices in
$GL_{m+n}$, identified with $T_m\times T_n$ where $T_\ell$ denotes the
corresponding torus in
$GL_\ell$. There is a natural action of $T$ on $\Mmn$ which arises as the
restriction of the natural left action of $\GLm\times\GLn$ on $\Mmn$:
namely, $(a,b).x= axb^{-1}$ for $(a,b)\in T$ and $x\in \Mmn$. This action
of $T$ on
$\Mmn$ is by Poisson isomorphisms; in particular, the action of each
element of $T$ maps symplectic leaves of $\Mmn$ to symplectic leaves.
Thus, it is natural to look at $T$-orbits of symplectic leaves of 
$\Mmn$, which are regular Poisson submanifolds
of $\Mmn$, rather than at individual symplectic leaves. (Here and
throughout, we view the $T$-orbit of a symplectic leaf $\L$
as the set-theoretic union
$\bigcup_{t\in T} t.\L$, rather than as the family $(t.\L)_{t\in T}$
of symplectic leaves.) As advantages to this approach, we mention that
$T$-orbits of symplectic leaves are easier to identify than single
symplectic leaves, and these orbits exhibit direct relations with known
geometric and Lie-theoretic structures. For example, 
we prove that the $T$-orbits of
symplectic leaves in $\Mmn$ are isomorphic (as varieties) to
intersections of dual Schubert cells in the full flag variety of $GL_{m+n}$,
and each generalized double Bruhat cell in $\Mmn$ (corresponding to a
pair of partial permutation matrices) is a disjoint union of $T$-orbits
of symplectic leaves, containing one such orbit as an open dense
subset. One thus sees that the Poisson structure of $\Mmn$ is in some
ways similar to, but also more intricate than, that of the group
$\GLn$ -- for instance, as follows from the analysis of Hodges and
Levasseur \cite{\HLthree}, the orbits of symplectic leaves in $\GLn$
under left translation by the standard maximal torus are precisely the
double Bruhat cells.

In a sequel to this paper, we will investigate the relation between
the standard Poisson structures on different (partial) flag varieties of
a semisimple algebraic group. We will also relate the restriction
of the Poisson structure to various Poisson subvarieties with known
quadratic Poisson structures on affine spaces. A detailed study of
Poisson structures of the latter type associated to arbitrary Schubert
cells in flag varieties of semisimple groups will be presented as well.
\enddefinition

\definition{0.2}
Recall that a Poisson group structure on an algebraic group $G$ is thought 
of as the ``semiclassical limit'' of a quantization of $G$, a viewpoint
promulgated in particular by Drinfeld and his school (cf\.
\cite{\DriICM}). Relationships between the symplectic foliation of such a 
Poisson structure and the primitive spectrum of a quantized 
coordinate ring
$\O_q(G)$ are viewed under the heading of a generalized version of the
Kirillov--Kostant orbit method. In the work of Soibelman (e.g., 
\cite{\Soi}) on
compact groups $G$, this led to bijections between the symplectic leaves of
$G$ and the primitive ideals of $\O_q(G)$. Hodges and Levasseur
\cite{\HLthree, \HLn} then established analogous results for $G=\SLn$, which
were extended to all semisimple groups by Joseph \cite{\Jos}. We
take the corresponding viewpoint that the Poisson structure on
$\Mmn$ is the semiclassical limit of the structure of $\O_q(\Mmn)$,
and argue that the results of the present paper should correspond to the
framework of the primitive ideals in $\O_q(\Mmn)$. Specifically, we
conjecture that the sets of minors which define the $T$-orbits of symplectic
leaves in $\Mmn$ (obtainable from Theorem 4.2) should match the sets of
quantum minors which generate the prime ideals of $\O_q(\Mmn)$ invariant
under winding automorphisms (cf\. \cite{\GLijm}). Some relations between
$\Mmn$ and $\O_q(\Mmn)$ are already known. In particular, the set of
$T$-orbits of symplectic leaves in $\Mmn$, partially ordered by inclusions
of closures, is anti-isomorphic to the poset $T\text{-Spec\,} \O_q(\Mmn)$ of
winding-invariant prime ideals in $\O_q(\Mmn)$ -- our work shows that the
former poset is anti-isomorphic to the set
$$S_{m+n}^{\le\wmwn}= \left\{ y\in S_{m+n} \bigm| y\le \left( \smallmatrix
1 &2 &\cdots &n &n+1 &n+2 &\cdots &n+m \\ m+1 &m+2 &\cdots &m+n &1 &2
&\cdots &m \endsmallmatrix \right) \right\}$$
under the Bruhat order, while Launois \cite{\Lau, Theorem 5.6} has proved
that $S_{m+n}^{\le\wmwn}$ is isomorphic to $T\text{-Spec\,}
\O_q(\Mmn)$.
\enddefinition

\definition{0.3}
Before summarizing our main results, we indicate some notation,
beginning with $N=m+n$. By $\Gr(n,N)$ we denote the Grassmannian
of $n$-dimensional subspaces of an $N$-dimensional space.
We write $B^+_\ell$ and $B^-_\ell$ for the
standard Borel subgroups of any $GL_\ell$ (consisting of upper,
respectively lower, triangular matrices), and identify the Weyl group
of $GL_\ell$ with both the symmetric group $S_\ell$ and the group of
permutation matrices in $GL_\ell$. The symbol $w_\ci^\ell$ denotes the
longest element of $S_\ell$. For $0\le t\le \ell$, let $S^1_t$ and
$S^2_{\ell-t}$ denote the natural copies of $S_t$ and $S_{\ell-t}$
inside $S_\ell$, acting on the numbers $1,\dots,t$ and
$t+1,\dots,\ell$, respectively. Finally, for a Weyl group $W$
and a subgroup $W_1$ generated by simple reflections, $W^{W_1}$ denotes 
the set of minimal length representatives for left cosets in $W/W_1$.
Recall that each such coset has a unique representative in 
$W^{W_1}$ (cf. \cite{\Car, Proposition 2.3.3(i)}).
 
The following theorem summarizes Theorems 3.9 and 3.13.
\enddefinition

\proclaim{Theorem A} {\rm (a)} There are only finitely many $T$-orbits
of symplectic leaves in $\Mmn$, and they are smooth irreducible locally
closed subvarieties.

{\rm (b)} The $T$-orbits of symplectic leaves in $\Mmn$ can be
described as the sets
$$
\P_w = \biggl\{
x\in\Mmn \biggm| \leftmatrix \wn &0\\ x &\wm \rightmatrix \in 
B_N^+wB_N^+ \biggr\},
$$
where $w\in S_N$ and $w\ge \leftmat \wn &0\\ 0 &\wm \rightmat$ in the
Bruhat order.

{\rm (c)} The closure of $\P_w$ equals the disjoint union of the $\P_z$
for $z\le w$.

{\rm (d)} As an algebraic variety, $\P_w$ is isomorphic to the
intersection of dual Schubert cells
$$
B_N^-. \leftmatrix \wn &0\\ 0 &\wm \rightmatrix B_N^+ \cap B_N^+. wB_N^+
$$
in the full flag variety $\GLN/B_N^+$. \qed\endproclaim

\definition{0.4}
Fulton \cite{\Ful} has given computational descriptions of Bruhat cells
$B_N^+wB_N^+$ in terms of ranks of rectangular submatrices. We apply
his results to the sets $\P_w$, to characterize exactly which matrices
$x$ lie in each $\P_w$, in terms of ranks of rectangular submatrices of
$x$. See Theorem 4.2 for the precise statement.
\enddefinition

\definition{0.5}
The results of Theorem A are obtained by embedding $\Mmn$ in 
the Grassmannian $\Gr(n,N)$ which, equipped with an 
appropriate Poisson structure, becomes a Poisson
homogeneous space for the standard Poisson algebraic group 
$\GLN$. (For details on Poisson homogeneous spaces 
for Poisson algebraic groups see \cite{\Dri} or Section 1.) 
This approach provides, in addition, a natural Poisson compactification
of $\Mmn$ which, in particular, suggests an approach to the problem of
studying the spectrum of $\O_q(M_{m,n})$ via noncommutative projective
geometry.

A completely different viewpoint is 
obtained by focussing, as we do in
Sections 5 and 6, on the sets $\Omnt$ of matrices in $\Mmn$ with a
fixed rank $t$. Each $\Omnt$ is a Poisson homogeneous space for the
natural action of $\GLm\times\GLn$ (equipped with an appropriate 
Poisson group structure). The latter group is the Levi factor of the
maximal parabolic subgroup of $\GLN$ defining $\Gr(n,N)$. The key results of
this approach, taken from Theorems 5.11, 6.1, 6.4 and Corollary 6.5, are as
follows.
\enddefinition

\proclaim{Theorem B} Fix a nonnegative integer $t\le \min\{m,n\}$.

{\rm (a)} The $T$-orbits of symplectic leaves within $\Omnt$ can be
described as the sets
$$
\P^t_{(y,v,z,u)} = \bigcup \Sb \tau\in S^1_t\\
z\tau\le y,\ v\tau^{-1}\le u \endSb
\bigl( B_m^+yB_m^+\cap B_m^-z\tau \bigr). \leftmatrix I_t &0\\ 0&0
\rightmatrix .\bigl( \tau^{-1}B_n^-u^{-1}B_n^-\cap v^{-1}B_n^+ \bigr),
$$
where $(y,v,z,u)\in S_m^{S^2_{m-t}} \times S_n^{S^1_t}
\times S_m^{S^1_t} \times S_n^{S^2_{n-t}}$ and $z\le y$, $v\le u$.

{\rm (b)} For $(y,v,z,u)$ as in {\rm (a)}, the set
$$
\C_{y,z}.\R_{u,v}= \bigl( B_m^+yB_m^+\cap B_m^-z \bigr). \leftmatrix I_t
&0\\ 0&0 \rightmatrix .\bigl( B_n^-u^{-1}B_n^-\cap v^{-1}B_n^+ \bigr)
$$
is dense in $\P^t_{(y,v,z,u)}$.

{\rm (c)} The sets $\C_{y,z}= \bigl( B_m^+yB_m^+\cap B_m^-z \bigr).
\leftmatrix I_t\\ 0 \rightmatrix$ are $(T_m\times T_t)$-orbits of
symplectic leaves of $M_{m,t}$, and each of the sets consisting of all
matrices in
$M_{m,t}$ with rank $t$ and a given column-echelon form is a disjoint
union of certain $\C_{y,z}$.

{\rm (d)} The sets $\R_{u,v}= \leftmatrix I_t
&0 \rightmatrix .\bigl( B_n^-u^{-1}B_n^-\cap v^{-1}B_n^+ \bigr)$ are
$(T_t\times T_n)$-orbits of symplectic leaves of $M_{t,n}$, and each of
the sets consisting of all matrices in $M_{t,n}$ with rank $t$ and a
given row-echelon form is a disjoint union of certain $\R_{u,v}$.
\qed\endproclaim

The descriptions of torus orbits of symplectic leaves in $\Mmn$
given in part (b) of Theorem A and part (a) of Theorem B are matched
in Theorem 5.11 and Proposition 5.9.

\definition{0.6}
Finally, we study the decomposition of $\Mmn$ into generalized double
Bruhat cells
$$
\Bww= B_m^+w_1B_n^+ \cap B_m^-w_2B_n^-,
$$
for partial permutation matrices $w_1$, $w_2$. If $w_1$ and $w_2$ have
the same rank $t$ (which is necessary for $\Bww$ to be nonempty), there
are unique decompositions
$$\xalignat2
w_1 &= y \leftmatrix I_t &0\\ 0&0 \rightmatrix v^{-1}  &w_2 &= z
\leftmatrix I_t &0\\ 0&0 \rightmatrix u^{-1}  \endxalignat
$$
where $y\in S_m^{S^2_{m-t}}$, $v\in S_n^{S^1_t S^2_{n-t}}$,
$z\in S_m^{S^1_t S^2_{m-t}}$, and $u\in S_n^{S^2_{n-t}}$ (see Lemma
7.3). The following results are given in Theorem 7.4.
\enddefinition

\proclaim{Theorem C} Let $w_1,w_2\in \Mmn$ be partial permutation
matrices with rank $t$, decomposed as above.

{\rm (a)} $\Bww$ is nonempty if and only if $z\le
y$ and $v\le u$, in which case it is a $T$-stable Poisson
subvariety of $\Mmn$, and a smooth irreducible locally closed
subvariety.

{\rm (b)} The partition of $\Bww$ into $T$-orbits of symplectic
leaves is given by
$$
\Bww= \bigsqcup\, \biggl\{ \P^t_{(y,v\tau_2,z\tau_1,u)} \biggm|
\matrix \tau_1\in S^2_{m-t}\subseteq S_m,\ z\tau_1\le y \\
\tau_2\in S^2_{n-t}\subseteq S_n,\ v\tau_2\le u \endmatrix \biggr\}.
$$

{\rm (c)} $\P^t_{(y,v,z,u)}$ is Zariski open and dense in $\Bww$.
\qed\endproclaim

\definition{0.7}
Let us also note that the standard Poisson algebraic group $GL_m$
is a $T$-stable Poisson subvariety of $M_{m,m}$. Thus the  
$T$-orbits of symplectic leaves of $GL_m$ (which are the same as 
the $T_m$-orbits of leaves) comprise a subset of the $T$-orbits
of symplectic leaves of $M_{m,m}$. The former are the double
Bruhat cells $B_m^+ w_1 B_m^+ \cap B_m^- w_2 B_m^-$ of $GL_m$, for
$w_1, w_2 \in S_m$. They 
were studied in detail by Fomin and Zelevinsky in \cite{\FZone},
who in a joint work with Berenstein also proved 
\cite{\BFZ} that their rings of regular functions 
provide important examples of upper {\it{cluster algebras}}
\cite{\FZtwo}.
Our results in particular show that the double Bruhat cells
in $GL_m$ are special cases of intersections of dual 
Schubert cells on the full flag variety of $GL_{2m}$.
It would be very interesting to understand whether 
any intersection of dual Schubert cells on the 
full flag variety of an arbitrary reductive algebraic group 
gives rise to a cluster algebra in the sense of Fomin 
and Zelevinsky \cite{\FZtwo}. If this is true, it will
imply that any $T$-orbit of symplectic leaves
of $\Mmn$ is the spectrum of a cluster algebra.
\enddefinition

\definition{0.8}
We conclude the introduction with some remarks on our
notation and conventions. All manifolds and algebraic varieties considered 
in this paper are over the field
of complex numbers.

Given an algebraic group $G$ with tangent Lie algebra
$\gfrak$, we denote by $L(\gamma)$ and
$R(\gamma)$ the left and right invariant
multi-vector fields on $G$ corresponding to $\gamma \in \wedge
\gfrak$. If $G$ acts on a smooth  quasiprojective
variety $M$, we will denote by 
$$
\chi  : \wedge \gfrak \rightarrow
\Gamma(M, \wedge T M)
\tag 0.1
$$
the extension of the infinitesimal 
action of $\gfrak$ on $M$ to $\wedge \gfrak$.
In the special case of the 
left and right multiplication actions
of $G$ on itself ($g.a=ga$ and $g.a = a g^{-1}$),
the above infinitesimal actions will be denoted
by
$$
\chi^R, \chi^L : \wedge \gfrak \rightarrow
\Gamma(G, \wedge T G).
\tag 0.2
$$
Note that for $\gamma \in \wedge \gfrak$,
$$\xalignat2
\chi^L(\gamma) &=R(\gamma)
 &\chi^R(\gamma) &= (-1)^{\epsilon(\gamma)} L(\gamma),
\tag 0.3\endxalignat
$$
where $\epsilon(\gamma)$ is the parity of $\gamma$.

If $Y$ is a locally closed subvariety of an algebraic
variety $X$ and $Z\subseteq Y$, we will denote the closure
of $Z$ in $Y$ by 
$\Cl_Y(Z)$. By a {\it stratification\/} of an
algebraic  variety $X$ we mean a partition of $X$ into smooth, irreducible,
locally closed subvarieties,
$X = \bigsqcup_{\alpha \in A} X_\alpha$, such that
for each $\alpha \in A$, we have
$\Xbar_\alpha = \bigsqcup_{\beta \in A(\alpha)} X_\beta$ for some
index set
$A(\alpha) \subseteq A$. 

We will use the following convention to distinguish 
double cosets from orbits of cosets. {\it For any subgroups $C$ 
and $D$ of a group $G$\/}:
\roster
\item"(i)" {\it The $(C,D)$ double coset of $g \in G$ will 
be denoted by $C g D$\/};
\item"(ii)" {\it The $C$-orbit of $g D$ in $G/D$
will be denoted by $C . g D$.\/}
\endroster 
The adjoint action of $g\in G$ on $h\in G$ will be written as
$\Ad_g(h)= ghg^{-1}$. 
\enddefinition

\head 1. Poisson algebraic groups 
and Poisson homogeneous spaces \endhead

We begin with background and notation for Poisson algebraic groups 
and Poisson homogeneous spaces, and then characterize the symplectic
leaves and their orbits in certain Poisson homogeneous spaces.

\definition{1.1\. Poisson varieties} 
Recall that a {\it Poisson manifold\/} is a pair $(X,\pi)$
consisting of a smooth manifold $X$ together with a Poisson bivector
field
$\pi \in \wedge^2 T X$, that is, $[\pi,\pi]=0$ where $[.,.]$
denotes the Schouten bracket. A (not necessarily closed)
submanifold $Y$ of
$X$  is called a {\it Poisson submanifold\/} if 
$\pi_y \in \wedge^2 T_y Y$ for all $y \in Y.$ 
In this case $(Y, \pi|_Y)$ is a Poisson manifold as well.
A (not necessarily closed) submanifold $Y$ of
$X$  is called a {\it complete Poisson submanifold\/} if
it is stable under all Hamiltonian flows. 
Any complete Poisson submanifold is 
a Poisson submanifold. The converse is not
necessarily true but, if $(X, \pi)$ is a Poisson manifold which is
partitioned into a disjoint union of Poisson submanifolds 
$X = \bigsqcup_{\alpha \in A} Y_\alpha,$
then all $Y_\alpha$ are complete Poisson submanifolds, 
see \cite{\HY, Lemma 3.2}.

The Poisson manifold $(X, \pi)$ is  
{\it regular\/}, respectively {\it symplectic\/},
if $\rank \pi$ is constant, respectively $\rank \pi = \dim X$. 
A {\it symplectic leaf\/} of $(X, \pi)$ is a maximal connected 
(not necessarily closed)
symplectic submanifold. It is well known that
any Poisson manifold $(X, \pi)$ can be decomposed into 
a disjoint union of its symplectic leaves, see e.g. 
\cite{\W, \V}. 
Note that a (not necessarily closed) submanifold $Y$ of
$X$ is a complete Poisson submanifold if and only if
it is a union of symplectic leaves of $(X, \pi)$.  

Let us also recall that a map
$\phi : (X, \pi) \rightarrow (Z, \pi')$ 
between two Poisson manifolds is called 
a {\it Poisson map\/} if
$\phi_*(\pi_x) = \pi'_{\phi(x)}$ for all $x \in X$.
For instance if $Y$ is a Poisson submanifold 
of $(X, \pi)$, the natural inclusion
$i : (Y, \pi|_Y) \hookrightarrow (X, \pi)$ 
is Poisson.

{\it{All Poisson manifolds considered 
in this paper will be {\rm(}complex{\rm)} smooth quasiprojective 
Poisson varieties. The symplectic leaves of a smooth 
quasiprojective Poisson variety are not necessarily
algebraic, i.e., smooth irreducible locally closed 
subvarieties. We will see below that this is the 
case for many Poisson varieties admitting 
appropriate transitive algebraic group actions.}}
\enddefinition

\definition{1.2\. Poisson algebraic groups
and Manin triples} A {\it Poisson algebraic group\/}
is an algebraic group $G$ equipped with a Poisson 
bivectorfield $\pi \in \wedge^2 TG$ such that the 
map 
$$(G, \pi) \times (G, \pi) \rightarrow (G, \pi)$$
is Poisson. The tangent Lie algebra 
$\gfrak = \Lie(G)$ of a Poisson algebraic group
$(G, \pi)$ has a canonical Lie bialgebra structure; see \cite{\CP,
\S1.3} and \cite{\KS, \S3.3} for details. 

Recall that a {\it Manin triple of Lie algebras\/} is a triple 
$(\dfrak, \afrak, \bfrak)$ with the following properties:
\roster
\item $\dfrak$ is a Lie algebra, $\afrak$ and $\bfrak$ 
are Lie subalgebras of $\dfrak$, and $\dfrak$ is the 
vector space direct sum of $\afrak$ and $\bfrak$.
\item $\dfrak$ is equipped with a nondegenerate invariant 
bilinear form with respect to which both $\afrak$ 
and $\bfrak$ are Lagrangian (i.e., maximal isotropic) subspaces.
\endroster

To any Lie bialgebra $\gfrak$ one associates 
the Manin triple $(D(\gfrak), \gfrak, \gfrak^*)$.
Here $D(\gfrak)$ and $\gfrak^*$ are the underlying Lie algebras 
of the double and the dual Lie bialgebras of $\gfrak$. 
The bilinear form on $D(\gfrak)$ is given by
$\langle x + \alpha, y + \beta \rangle= \beta(x)+\alpha(y)$
for $x, y \in \gfrak$, $\alpha, \beta \in \gfrak^*$.
\enddefinition

\proclaim{1.3\. Definition} A \underbar{Manin triple of algebraic
groups} is a triple $(D, A, B)$ of algebraic groups  such that 
$A$ and $B$ are algebraic subgroups of $D$ and 
$(\Lie(D), \Lie(A), \Lie(B))$ is a Manin triple of
Lie algebras. 
\endproclaim

Fix a Manin triple of algebraic groups $(D, A, B)$.  
Then $D$ has a canonical Poisson algebraic group structure 
with a Poisson bivector 
field given by
$$
\pi^D = L (r) - R(r)=\chi^R(r) - \chi^L(r) \quad \text{where}
\quad r= \sum_i x_i \wedge x^i \in \wedge^2 \Lie D  
$$
in the notation (0.2)--(0.3) for left and right 
invariant multi-vector fields $L(.)$, $R(.)$
on $D$ and infinitesimal actions $\chi^L(.)$, 
$\chi^R(.)$ of $\Lie D$ on $D$. 
Here $\{x_i\}$ and $\{x^i\}$ are dual bases of
$\Lie(A)$ and $\Lie(B)$, respectively, with respect
to the nondegenerate bilinear form
on $\Lie(D)$.

The groups $A$ and $B$ are Poisson subvarieties of $D$.
The Poisson algebraic group $(D, \pi^D)$ is a {\it double\/} 
of $(A, \pi^D|_A)$, and $(B, -\pi^D|_B)$ is a 
{\it dual\/} Poisson algebraic group of 
$(A, \pi^D|_A)$; cf\. \cite{\CP, \S1.4} and \cite{\KS,
\S3.3}.

{\it{We will say that a Poisson algebraic group $(G, \pi)$ 
is a part of a Manin triple of algebraic groups 
$(D, G, F)$ if the Poisson structure $\pi$ coincides 
with the Poisson structure $\pi^D|_G$ induced from $D$.}}  

\definition{1.4\. Standard Poisson structures on reductive
algebraic groups} 
Let $G$ be a complex reductive algebraic group. 
The standard Poisson structure on $G$, turning it into a 
Poisson algebraic group, is defined as follows. 
Fix two opposite Borel subalgebras $\bfrak^\pm$ of 
$\gfrak = \Lie G$ and set $\hfrak = \bfrak^+ \cap \bfrak^-$
for the corresponding Cartan subalgebra of $\gfrak$. 
Fix a nondegenerate bilinear invariant form 
$\langle.,.\rangle$ on $\gfrak$ 
for which the square of the length of a long root is equal
to 2. Choose sets of root vectors $\{e_\alpha\}$
and $\{f_\alpha\}$, spanning respectively 
the nilradicals $\nfrak^+$ and $\nfrak^-$ 
of $\bfrak^+$ and $\bfrak^-$, normalized by $\langle e_\alpha, f_\alpha
\rangle =1$.

The standard $r$-matrix of $\gfrak$ is given by
$$
r = \sum_\alpha e_\alpha \wedge f_\alpha 
\tag 1.1
$$
and the corresponding standard Poisson structure on 
$G$ is defined by
$$
\pi = L(r) - R(r)= \chi^R(r) - \chi^L(r),
\tag 1.2
$$ 
in the notation (0.2)--(0.3).

The standard $r$-matrix on $G = GL_N$ is
$$
r^N = \sum_{1 \leq i < j \leq n} E_{ij} \wedge E_{ji}
\in \wedge^2 \gl_N
\tag 1.3
$$  
where the $E_{ij}$ are the standard elementary matrices.

{\it{By abuse of notation, $\GLN$ will denote
the algebraic group $\GLN$ equipped with the 
standard Poisson structure $\pi^N$ from {\rm(1.2)}, associated 
to the $r$-matrix $r^N$ {\rm(1.3)}. By $\GLN^\bu$ 
we will denote the Poisson algebraic 
group $(\GLN, - \pi^N)$.}}

Any standard (complex) reductive Poisson algebraic group 
$(G, \pi)$ is a part of the Manin triple 
$(G \times G, \Delta(G), F)$ where $\Delta(G)$
is the diagonal of $G \times G$ and 
$$
F = \{ (h u^+, h^{-1} u^-) \mid h \in T,\, u^\pm \in U^\pm \}
\subseteq B^+ \times B^-,
\tag 1.4
$$
where $B^\pm$ are the Borel subgroups of $G$ corresponding 
to $\bfrak^\pm$, $U^\pm$ are their unipotent radicals, and 
$T= B^+ \cap B^-$ is the corresponding maximal torus of
$G$. For the standard Poisson structure on $G$,
$$
\gfrak^* = \Lie F = 
\{ (h + n^+, -h + n^-) \mid h \in \hfrak,\, 
n^\pm \in \nfrak^\pm \} \subseteq \bfrak^+ \oplus \bfrak^-.
\tag 1.5
$$
The nondegenerate invariant bilinear
form on 
$\Lie(G \times G) \cong \gfrak \oplus \gfrak$,
used in the Manin triple of Lie algebras
$(\gfrak \oplus \gfrak, \Delta(\gfrak), \gfrak^*)$,
is
$$ \langle (x_1, x_2), (y_1, y_2) \rangle =
\langle x_1, y_1 \rangle -
\langle x_2, y_2 \rangle,
\tag 1.6
$$
where in the right hand side 
$\langle., .\rangle$
denotes the bilinear form on $\gfrak$, 
fixed above.
\enddefinition

\definition{1.5\. Matrix affine Poisson spaces} 
The {\it{$m\times n$ matrix affine Poisson space}} 
is the affine space
$\AA^{mn}$, identified with the space $\Mmn$ 
of all $m\times n$ complex matrices. The standard Poisson
structure on $\Mmn$ is given by 
$$
\pi^{m,n} = \sum_{i, k=1}^m \sum_{j, l=1}^n 
\bigl( \sign(k-i) + \sign(l-j) \bigr) x_{il} x_{kj}
\frac{\partial}{\partial x_{ij}} \wedge
\frac{\partial}{\partial x_{kl}} 
\tag 1.7
$$
in terms of the standard coordinate functions $x_{ij}$ 
on $\Mmn$. {\it{By abuse of notation, $\Mmn$ will 
denote the matrix affine Poisson space, 
thus dropping the symbol for the Poisson
structure {\rm(1.7)} on $\Mmn$.}}

Note that $\GLm$ acts on $\Mmn$ by left multiplication
($g.x = gx$ for $g \in \GLm$, $x \in \Mmn$), and $\GLn$ acts 
on $\Mmn$ by (inverted) right multiplication
($g.x = x g^{-1}$ for $g \in \GLm$, $x \in \Mmn$). 
The extensions of the corresponding infinitesimal
actions of $\gl_m$ and $\gl_n$ on $\Mn$ to
$\wedge \gl_m$ and $\wedge \gl_n$ 
will be denoted by
$$
\chi^L : \wedge \gl_m \rightarrow 
\Gamma( \Mmn, T \Mmn) \quad 
\text{and} \quad
\chi^R : \wedge \gl_n \rightarrow
\Gamma( \Mmn, T \Mmn).
$$
Note that in the case $m=n$ these extend the infinitesimal
actions $\chi^L$ and $\chi^R$ 
of $\gl_m$ on $\GLm \subseteq M_m$, defined in (0.2).  

By direct computation one shows that the Poisson structure
(1.7) on $\Mmn$ is also given by the formula
$$
\pi^{m,n}= \chi^R(r^n) - \chi^L(r^m)
\tag 1.8
$$
in terms of the standard $r$-matrix $r^N$ for $GL_N$, 
see (1.3).

Note that $\GLn$ is a Poisson subvariety of $M_{n,n}$.
\enddefinition

\definition{1.6\. Poisson homogeneous spaces} 
Fix a Poisson algebraic group $(G, \pi)$
and set $\gfrak = \Lie(G)$. 
A {\it Poisson $(G, \pi)$-space\/} is a smooth quasiprojective
Poisson variety $(M, \pi_M)$ equipped with a morphic $G$-action for
which
$$(G, \pi) \times (M, \pi_M) \rightarrow (M, \pi_M)$$
is a Poisson morphism.

A {\it Poisson homogeneous space for $(G, \pi)$\/} is a Poisson 
$(G, \pi)$-space $(M, \pi^M)$ for which $M$ is a 
homogeneous $G$-space. (Recall that any homogeneous
space of an algebraic group is a smooth quasiprojective 
variety \cite{\Bor, Theorem 6.8}.) To each $m \in M$, one associates
the {\it Drinfeld subalgebra\/} \cite{\Dri}
$$\lfrak_m = \{ x + \alpha \in D(\gfrak) \mid 
x \in \gfrak,\, \alpha \in \gfrak^*,\, 
\alpha|_{\gfrak_m} = 0,\, 
\alpha \rfloor \pi_m^M = x + \gfrak_m \}$$
of the double $D(\gfrak)$,
where $\gfrak_m$ denotes the Lie algebra of the stabilizer
$G_m= \stab_G(m)$, the tangent space $T_m M$ is identified with 
$\gfrak/\gfrak_m$, and the Poisson bivectorfield $\pi_m^M$ is thought
of as an element of $\wedge^2 (\gfrak/\gfrak_m)$. Note that
$$ 
\gfrak_m = \gfrak \cap \lfrak_m.  \tag 1.9
$$

The Drinfeld subalgebras $\lfrak_m$ 
are moreover Lagrangian subalgebras
of the double $D(\gfrak)$, equipped with the canonical nondegenerate 
invariant bilinear form \cite{\Dri; \KS, Proposition 6.2.15}. The map,
associating to
$m \in M$ its  Drinfeld subalgebra $\lfrak_m \subseteq D(\gfrak)$, is
$G$-equivariant:
$$\lfrak_{g m} = \Ad_{g} (\lfrak_m)$$ 
where $\Ad_{g}$ refers to the adjoint action of $G$ on $D(\gfrak)$.
\enddefinition

\proclaim{1.7\. Definition} A Poisson homogeneous $(G, \pi)$-space
$(M, \pi^M)$ will be called \underbar{algebraic} if the  Drinfeld
subalgebra of some $m \in M$ is the tangent Lie algebra of an
algebraic subgroup
$L_m \subseteq D$.  
\endproclaim

Because of the $G$-equivariance of the map $m \mapsto \lfrak_m$,
if the condition in the definition is satisfied for one 
point $m \in M$, then it holds for any $m \in M$. 

An important type of Poisson homogeneous space $(M, \pi^M)$ is 
the class of those for which $\pi^M$ vanishes at some 
point of $M$. In the rest of this subsection 
we describe those.

An algebraic subgroup $Q$ of a Poisson algebraic group
$(G, \pi)$ will be called an
{\it{almost Poisson algebraic subgroup}\/}
if
$$ \pi_q \in T_q Q \wedge T_q G$$
for all $q \in Q$.
(Recall that if $\pi_q \in \wedge^2 T_q Q$ for all $q \in Q$,
then $Q$ is called a Poisson algebraic subgroup of $(G, \pi)$.)
Fix an almost Poisson algebraic subgroup $Q$ of $(G, \pi)$,
and consider the projection
$$p : G \rightarrow G/Q, \quad p(g) = g Q.$$
Then 
$$\pi_{gq} - R_q (\pi_g) \in 
L_g (T_{q} Q) \wedge T_{gq} G  
$$
for all $g \in G$, $q \in Q$,
and the rule
$$
\pi^{G/Q}_{gQ} = p_*( \pi_g), \quad g \in G
\tag 1.10
$$
gives a well-defined Poisson structure $\pi^{G/Q}$ on $G/Q$.
The pair $(G/Q, \pi^{G/Q})$ is a Poisson homogeneous 
space of $(G, \pi)$ and $\pi^{G/Q}$ vanishes at the base 
point $e Q$ of $G/Q$.

\proclaim{1.8\. Theorem} Fix a Poisson algebraic group $(G, \pi)$.

{\rm(a)} Any Poisson homogeneous $(G, \pi)$-space $(M, \pi^M)$ with 
the property that the Poisson bivectorfield $\pi^M$ vanishes 
at some point $m \in M$ is isomorphic to $(G/Q, \pi^{G/Q})$ for
$Q=\stab_G(m)$ which is an almost Poisson algebraic
subgroup of $(G, \pi)$. 

{\rm(b)} For an almost Poisson algebraic subgroup $Q$ of $G$,
the Drinfeld Lagrangian subalgebra of the base point $e Q$
of the Poisson homogeneous space $(G/Q, \pi^{G/Q})$ is
$$
\lfrak = \qfrak + \qfrak^\perp
\tag 1.11
$$
where $\qfrak = \Lie Q$ and $\qfrak^\perp$ 
refers to the orthogonal subspace to $\qfrak \subseteq \gfrak$
in $\gfrak^*$.

{\rm(c)} A connected algebraic subgroup $Q$ of $(G, \pi)$ is 
an almost Poisson algebraic subgroup if and only if 
the orthogonal complement 
$\qfrak^\perp \subseteq \gfrak^*$ 
is a subalgebra of the dual Lie bialgebra $\gfrak^*$ 
of $\gfrak$ {\rm(}as in part {\rm(b)}, $\qfrak = \Lie Q${\rm)}. 

{\rm(d)} A connected algebraic subgroup $Q$ of $(G, \pi)$ is a
Poisson algebraic subgroup if and only if $\qfrak^\perp$ is an 
ideal in $\gfrak^*$. \qed
\endproclaim

Parts (a) and (d) of this theorem can be found, e.g., 
in \cite{\KS, page 52 and Proposition 6.2.3}; parts (b) and (c) are
well known.

Below we gather some results on symplectic leaves 
of algebraic Poisson homogeneous spaces. Fix a Poisson algebraic group
$(G, \pi)$ which is a part of a Manin triple of algebraic groups
$(D, G, F)$, as defined in \S1.3. Fix also an algebraic
Poisson homogeneous 
$(G, \pi)$-space with connected stabilizer subgroups $G_m$ (see
\S1.6).  Such a homogeneous space has the form
$G/N$ where $N$ is a connected subgroup of $G$ and the Drinfeld
Lagrangian subalgebra of $\Lie (D)$ corresponding to the base
point $e N \in G/N$  integrates to an algebraic subgroup $L
\subseteq D$. Note that 
$$
N = (G \cap L)^\circ,
$$
the identity component of $G\cap L$,
because of (1.9) and the connectedness of $N$.
Consider the composition of maps
$$
\Pi : G/N @>{\mu}>> G/(G \cap L) @>{\,\cong\,}>> G. L 
\subseteq D/L,
\tag 1.12
$$ 
where $\mu$ is the map $gN \mapsto g(G\cap L)$. 

\proclaim{1.9\. Theorem} 
Assume that $(G, \pi)$ is a Poisson algebraic group
which is a part of a Manin triple of algebraic groups $(D, G, F)$.
Let $(G/N, \pi')$ be an algebraic 
Poisson homogeneous $(G, \pi)$-space with 
connected stabilizer subgroups for which the Drinfeld Lagrangian
subalgebra  of the base point $e N$ is $\Lie L$ for an algebraic
subgroup 
$L$ of $G$.

Then the symplectic leaves of $G/N$ are the connected components
{\rm(}i.e., irreducible components{\rm)} of the inverse images under
$\Pi$ of the 
$F$-orbits on $D/L$, and all of them are smooth irreducible 
locally closed subvarieties of $G/N$.
\endproclaim

Note that some $F$-orbits on $D/L$ might not intersect 
the image of $\Pi$, but when an $F$-orbit on $D/L$ intersects 
the image of $\Pi$, the intersection is transversal
since the Lie algebras of $G$ and $F$ span $\Lie D$. 
Below we will consider only those $F$-orbits on $D/L$ 
that intersect the image of $\Pi$.

\demo{Proof} Any $F$-orbit on $D/L$ is a smooth locally closed
subvariety, see e.g. \cite{\Bor, Proposition 1.8}. 
Thus, its inverse image under $\Pi$
(if it is nontrivial) is a locally closed subvariety of $G/N$.
Each intersection of an $F$-orbit on $D/L$ with $\Im \Pi = G.L$ 
is a transversal intersection of group orbits and therefore is
smooth. As a consequence its inverse image under the \'etale map 
$\Pi : G/N \rightarrow G. L$ is smooth as well.

Finally, the connected components of the (nontrivial)
inverse images of $F$-orbits are known to be 
symplectic leaves of $(G/N, \pi)$
due to results of Lu \cite{\Lu} and Karolinsky \cite{\Karone}
in the differential category. Since \cite{\Lu, \Karone}
assume that $D = F G$, we sketch another approach.
Consider the bivector field 
$\chi(r) \in \Gamma(D/L, \wedge^2 T D/L)$ where 
$r \in \wedge^2 \Lie D$ is the $r$-matrix for the Poisson structure 
on $D$, see Definition 1.3, and $\chi(.)$ refers to the 
natural infinitesimal action of $\Lie D$ on $D/L$.
It was proved in \cite{\LuY}
that $\chi(r)$ is a Poisson bivectorfield and that 
the connected components of the intersections of any $F$ and $G$ 
orbits on $D/L$ are symplectic leaves of $\chi(r)$.
It is straightforward to show that the map 
$\Pi : (G/N, \pi') \rightarrow (D/L, \chi(r))$ is Poisson.
The statement now follows from the fact that
$\Pi : G/N \rightarrow G.L$ is \'etale.
\qed\enddemo

In the remainder of this section, we gather some results
on orbits of symplectic leaves in Poisson homogeneous spaces. 
In the setting of Theorem 1.9, assume 
that $H$ is a subgroup of $G$ that normalizes $F \subseteq D$.
Then the Poisson structure $\pi$ on $G$ vanishes on $H$, 
see \cite{\LuY}, and as a consequence $H$ acts by Poisson 
isomorphisms on any Poisson homogeneous 
$(G, \pi)$-space $(M, \pi^M)$. This in particular means 
that each element $h \in H$ maps symplectic leaves
of $(M, \pi^M)$ to symplectic leaves. The $H$-orbits of symplectic
leaves are characterized in the following 
theorem which is adapted from
\cite{\LuY}. Let us first note that since $H$ normalizes 
$F \subseteq D$, the product $H F$  is an algebraic subgroup of $D$.   

\proclaim{1.10\. Theorem} In the setting of Theorem {\rm1.9}, the 
$H$-orbits of symplectic leaves of the Poisson homogeneous 
space $G/N$ are the irreducible components of the inverse   
images under $\Pi$ of the $H F$-orbits on $D/L$ {\rm(}see
{\rm(1.12))}, and all of them are smooth irreducible 
locally closed subvarieties of $G/N$.
\endproclaim

\demo{Proof} Fix $y \in G. L = \Im(\Pi) \subseteq D/L$. 
The intersection of $\Im \Pi =G.L$ with $F y$
is transversal because the Lie algebras of $HG$ and $F$ 
span $\Lie(D)$. Therefore $\Im \Pi \cap F y$ is a 
smooth and locally closed subset of $D/L$. The second statement
follows from the fact that both $G.L$ and $F y$ are 
locally closed subsets of $D/L$ (as orbits of 
algebraic groups). Let $\P$ be an irreducible component
of $\Pi^{-1}(HFy)$. It is a smooth, irreducible, 
locally closed subset of $G/N$ because $\Pi: G/N \rightarrow G.L$ is an
\'etale morphism, recall (1.12).

We need to show that $\P = H \SS$ for some irreducible component
$\SS$ of $\Pi^{-1}(F y)$. First, note that for two distinct 
irreducible components $\SS_1$ and $\SS_2$ of $\Pi^{-1}(Fy)$,
either $H \SS_1= H \SS_2$ or $H \SS_1$ and $H \SS_2$
are disjoint. Since the map $\Pi$ is $H$-equivariant,
$$
\Pi^{-1}(H F y) =  H \Pi^{-1}(F y). \tag 1.13  
$$
As a consequence,
$$
\P = H \SS_1 \sqcup \ldots \sqcup H \SS_m
$$
for some irreducible components $\SS_i$ of 
$\Pi^{-1}(F y)$, lying inside $\P$. All that we 
need to show now is that $m=1$. Since 
$\P$ is irreducible it is sufficient to show that 

{\it{For each irreducible component $\SS$ of $\Pi^{-1}(F y)$,
the set $H \SS$ is an open subset of $\P$.}}

We show this in the rest of the proof. Let 
$x' \in G. L =\Im \Pi$. Since $H$ normalizes $F$,
$$ 
T_{x'} (HF x') =
T_{x'} (H x') +
T_{x'} (F x').
$$
The intersections 
$HF x' \cap G x'$ 
and
$F x' \cap G x'$ 
in $D/L$ 
are transversal
because the Lie algebras
of $F$ and $G$ span $\Lie(D)$. Taking into account
this and the facts that $H$ is a subgroup of $G$ 
and $\Im \Pi = G . x' \supset H . x'$  gives
$$
T_{x'} (HF x' \cap \Im \Pi) =
T_{x'} (H x' ) +
T_{x'} (F x' \cap \Im \Pi).
$$
Since $\Pi$ is an \'etale map,
recall (1.12), we obtain
$$
T_x \P = T_x (H x) + T_x \SS \qquad\text{for all\ } x \in \SS.
$$
If $f : H \times \SS \rightarrow \P$ denotes the map
$(h, x) \mapsto hx$,
then the above equality implies that $d f$ is surjective
at any point of $H \times \SS$. As a consequence of this,
the morphism $f$ is smooth and thus flat, because 
$H \times \SS$ and $\P$ are nonsingular, see 
\cite{\Har, \S III, Proposition 10.4}. 
The latter implies that $f$ is open, 
see \cite{\Har, \S III, Problem 9.1}. Therefore 
the image of $f$ (which is nothing but $H \SS$) is
an open subset of $\P$. 

In fact, since we work over $\CC$ the last 
statement is almost immediate: the fact that 
the differential of 
$f : H \times \SS \rightarrow \P$ 
is surjective everywhere implies that 
the image of $f$ is open in the classical topology.
But $\Im f$ is also a constructible subset of 
$\P,$ thus it is a Zariski open subset.
\qed\enddemo

\head 2. Intersections of Bruhat and Schubert cells 
\endhead

Our main results rely on certain combinatorial and geometric
information about intersections of Bruhat and Schubert cells,
which we develop in this section.

\definition{2.1\. Bruhat and Schubert cells}
Let $G$ be a complex reductive algebraic group.
As in \S1.4, fix two opposite Borel
subgroups $B^\pm$ of $G$ and set $T = B^+ \cap B^-$
for the corresponding maximal 
torus  of $G$. Denote the projection 
to the flag variety by
$$\eta : G \rightarrow G/B^+. \tag 2.1$$
Recall that the $(B^\pm, B^\pm)$-double 
cosets of $G$ are called 
{\it{Bruhat cells}\/} of $G$ and the $B^\pm$-orbits
on $G/B^+$ are called {\it{Schubert cells}\/}
of $G/B^+$.

Let $U^\pm$ be the unipotent
radical of $B^\pm$. Denote by $W$ the Weyl group 
of $(G, T)$, by $\leq$ the Bruhat 
order on $W$, and by $l(.)$ 
the length function on $W$.
For each $w \in W$, fix a
representative $\dot{w}$ in the normalizer 
of $T$. When the result of a 
formula involving some $\dot{w}$ does not
depend on the particular representative $\dot{w}$
of $w$, the notation for such a representative 
will be omitted.
As a consequence of the Bruhat lemma,
all Bruhat cells of $G$ are uniquely
represented in the form $B^\pm w B^\pm$ 
for some $w \in W$ and all Schubert cells 
of $G/B^+$ are uniquely represented in the 
form $B^\pm . w B^+$ for some $w \in W$.

For each $w \in W$, define the following subgroups 
of $U^\pm$:
$$ 
U_w^- = U^- \cap \Ad_w(U^-) \quad 
\text{and} \quad 
U_w^0 = \Ad_w^{-1}(U^-) \cap U^+. \tag 2.2
$$
Recall that $U^-$, $U_w^-$, and $U_w^0$ are 
affine spaces (and closed subvarieties of $G$),
and as such,
$$
U_w^- \times \Ad_w(U_w^0) \cong U^-, \tag 2.3
$$
with the isomorphism given by group multiplication (e.g., 
see \cite{\Bor, \S14.12, p\. 193}).

In Theorem 2.3, for all $y, z \in W$
we describe the structure of the 
locally closed subvarieties 
$$
B^- z \cap B^+ y B^+, \quad
U^- \dot{z} \cap B^+ y B^+, \quad \text{and}
\quad
U_z^- \dot{z} \cap B^+ y B^+ \tag 2.4 
$$
of the intersection of Bruhat cells 
$B^- z B^+ \cap B^+ y B^+$ in terms of the 
intersection of the dual Schubert cells
$$
\B_{z,y}=
B^-. z B^+ \cap B^+. y B^+ \subseteq G/B^+.
\tag{2.5}
$$
The first two varieties in (2.4) are smooth 
due to the transversality of the 
intersections 
$(\Lie U^- + \Lie B^+ = \Lie G)$.
It will be shown in Theorem 2.3 that
the third variety in (2.4) is also smooth.  
In Theorem 2.5, we describe the Zariski 
closures in $G$ of the sets in (2.4).

First recall the following result of Deodhar, 
\cite{\Deo, Corollary 1.2}:
\enddefinition

\proclaim{2.2\. Proposition} {\rm [Deodhar]} For $y, z \in W$, the  
intersection $\B_{z,y}= B^-. z B^+ \cap B^+. y B^+$ of dual Schubert
cells   is nonempty if and only if $y \geq z$ in 
the Bruhat order of $W$. In that case, the intersection is 
a smooth irreducible locally closed subvariety
of $G/B^+$ of dimension $l(y) - l(z)$. \qed
\endproclaim

The smoothness in the second part of the Proposition 
is a direct consequence of the transversality
of the intersection. The harder result in the second 
part is the irreducibility. It follows from 
a stratification of the intersection by smooth 
irreducible locally closed subvarieties isomorphic 
to $\CC^n \times \left( \CC^\times \right)^m$,
$n, m \in \ZZ$, obtained by Deodhar \cite{\Deo, Theorem 1.1},  
in which only one set has dimension equal 
to $l(y) - l(z)$. A direct consequence of the first
part of the 
Proposition is that the intersection of Bruhat cells
$B^- z B^+ \cap B^+ y B^+$ is nonempty if and only 
$y \geq z$.

\proclaim{2.3\. Theorem} Let
$y, z \in W$ with $y \geq z$.

{\rm(a)} The projection $\eta : G \rightarrow G/B^+$
restricts to a biregular isomorphism of 
affine spaces
$$
\eta : U^-_z\dot{z}  @>{\,\cong\,}>> B^-. z B^+.
\tag 2.6 
$$
The set $U^-_z \dot{z} \cap B^+ y B^+$
is a smooth irreducible locally closed subset of $G$, and 
$\eta$ further restricts to a biregular
isomorphism 
of quasiprojective varieties
$$
\eta : U^-_z \dot{z} \cap B^+ y B^+ @>{\,\cong\,}>>
B^-. z B^+ \cap B^+. y B^+.
\tag 2.7
$$

{\rm(b)} The group multiplication in $G$ restricts to 
biregular isomorphisms of quasiprojective varieties
$$
\left( U^-_z \dot{z} \cap B^+ y B^+ \right) \times U^0_z
@>{\,\cong\,}>> U^- \dot{z} \cap B^+ y B^+ 
\tag 2.8
$$
and
$$
\left( U^-_z \dot{z} \cap B^+ y B^+ \right) \times 
U^0_z \times T
@>{\,\cong\,}>> B^- z \cap B^+ y B^+.
\tag 2.9
$$
\endproclaim

\demo{Proof} 
(a)  The first statement (2.6) is well known. (E.g., see \cite{\Bor,
Theorem 14.12(b)} for the analogous isomorphism $U^+\cap \Ad_w(U^-)
\rightarrow B^+.wB^+$.) Because
$U^-_z$ is a closed subvariety of
$G$, to complete the proof of part (a), 
all that we need to show is
$$
\eta \left(
U^-_z \dot{z} \cap B^+ y B^+ \right) =
B^-. z B^+ \cap B^+. y B^+. \tag 2.10
$$
It is obvious that
$$
\eta \left(
U^-_z \dot{z} \cap B^+ y B^+ \right)
\subseteq 
B^-. z B^+ \cap B^+. y B^+.
$$
But
$$
\eta \left( B^- z B^+ \cap B^+ y B^+ \right) = 
B^-. z B^+ \cap B^+. y B^+,
$$
and $B^-zB^+ \subseteq U_z^-\dot{z}B^+$ because of (2.3), so that
$B^- z B^+ \cap B^+ y B^+ \subseteq 
(U_z^-\dot{z} \cap B^+ y B^+) B^+$.
The surjectivity in (2.10) now follows from the isomorphism (2.6).

(b) First note that the right action of 
$U_z^0\subseteq B^+ \cap \Ad_z^{-1} U^-$ on $G$ preserves the 
intersection on the right hand side of (2.8), that is, 
$$
U^- \dot{z} \cap B^+ y B^+ \supset 
\left( U^-_z \dot{z} \cap B^+ y B^+ \right) U^0_z.
$$
To show the opposite inclusion, let
$$
g \in U^- \dot{z} \cap B^+ y B^+.
$$ 
Multiplying (2.3) on the right by $\dot{z}$, we get that
$$
g = g_1 u \quad \text{for some} \quad
g_1 \in U^-_z \dot{z} \; \text{and} \;
u \in U^0_z.
$$
Since $B^+ y B^+ U_z^0 =B^+ y B^+$, we obtain that 
$g_1 = g u^{-1} \in 
U^-_z \dot{z} \cap B^+ y B^+$ and thus
$$
g = g_1 u \in 
\left( U^-_z \dot{z} \cap B^+ y B^+ \right) U^0_z.
$$
Therefore
$$
U^- \dot{z} \cap B^+ y B^+ =
\left( U^-_z \dot{z} \cap B^+ y B^+ \right) U^0_z,
$$
which together with (2.3) implies (2.8).

In a similar way one proves (2.9), using
(2.8) and 
$$
B^- z \cap B^+ y B^+ =
\left( U^- \dot{z} \cap B^+ y B^+ \right) T. \qquad\square
$$ 
\enddemo

The following Theorem combines and summarizes 
Proposition 2.2 and Theorem 2.3.

\proclaim{2.4\. Theorem} For any $y, z \in W$ with $y \leq z$, the
sets
$U_z^- \dot{z} \cap B^+ y B^+$, $U^- \dot{z} \cap B^+ y B^+$
and $B^- z \cap B^+ y B^+$ are smooth irreducible locally closed
subvarieties of the intersection of Bruhat cells 
$B^- z B^+ \cap B^+ y B^+ \subseteq G$. They
are related to the intersection 
of dual Schubert cells 
$\B_{z,y}= B^- . z B^+ \cap B^+ . y B^+ \subseteq G/B^+$
by the following biregular isomorphisms, obtained as compositions 
of the isomorphisms {\rm(2.7)--(2.9)}:
$$
\align
U^-_z \dot{z} \cap B^+ y B^+ &\cong
\B_{z,y}
\\
U^- \dot{z} \cap B^+ y B^+ &\cong 
\B_{z,y} \times U^0_z
\\
B^- z \cap B^+ y B^+ &\cong 
\B_{z,y} \times U^0_z \times T. \qquad\square
\endalign
$$
\endproclaim

The first of the intersections above will play an important role in
the following section. We label it as follows
$$
\U_{z,y}= U^-_z \dot{z} \cap B^+ y B^+ \tag 2.11
$$
for $y,z\in W$.

\proclaim{2.5\. Theorem} For any $y, z \in W$ with $y \leq z$,
the Zariski closures of the three locally closed subsets of $G$
considered in Theorem {\rm 2.4} are given by
$$\align
\overline{U^-_z \dot{z} \cap B^+ y B^+} &=
U^-_z \dot{z} \cap \overline{B^+ y B^+} =
\bigsqcup \Sb w \in W \\ z \leq w \leq y \endSb
U^-_z \dot{z} \cap B^+ w B^+ \tag\roman{a} \\ 
\overline{U^- \dot{z} \cap B^+ y B^+} &=
U^- \dot{z} \cap \overline{B^+ y B^+} =
\bigsqcup \Sb w \in W \\ z \leq w \leq y \endSb
U^- \dot{z} \cap B^+ w B^+ \tag\roman{b} \\
\overline{B^- z \cap B^+ y B^+} &=
B^- z \cap \overline{B^+ y B^+} =
\bigsqcup \Sb w \in W \\ z \leq w \leq y \endSb
B^- z \cap B^+ w B^+. \tag\roman{c}
\endalign$$
\endproclaim

In the proof of Theorem 2.5, we will need
the following algebrogeometric fact.

\proclaim{2.6\. Lemma} Let $\bigsqcup_{\alpha \in A} X_\alpha$ 
be a  stratification {\rm(}cf\. {\rm\S0.8{\rm)} of a smooth algebraic
variety
$X$, and
$Y$  a smooth, irreducible, locally closed 
subvariety of $X$ that 
intersects all the strata $X_\alpha$ transversely. 
Then
$$
\Cl_Y(Y \cap X_\alpha) = Y \cap \Xbar_\alpha
$$
for all $\alpha \in A$.
\endproclaim

\demo{Proof} Fix $\alpha \in A$. Then 
$\Xbar_\alpha = 
\bigsqcup_{\beta \in A(\alpha)} X_\beta$ for some subset
$A(\alpha) \subseteq A$, and 
$\dim X_\beta< \dim X_\alpha$ for all $\beta \in A(\alpha) \setminus
\{\alpha\}$. 

Because $\Xbar_\alpha$ is a closed subvariety of $X$ 
that contains $X_\alpha$, the set $\Cl_Y(Y \cap
X_\alpha)$  equals the union of those irreducible components
of 
$$
Y \cap \Xbar_\alpha = \bigsqcup_{\beta \in A(\alpha)} 
Y \cap X_\beta
$$ 
that meet 
$Y\cap X_\alpha$. On one hand, the dimension
of any irreducible component of $Y\cap \Xbar_\alpha$
is greater than or equal to $\dim Y + \dim X_\alpha - \dim X$;
see \cite{\Har, Chapter I, Proposition 7.1 and Theorem 7.2}. 
On the other hand, for all
$\beta \in A(\alpha) \setminus
\{\alpha\}$, 
$$
\dim(Y \cap X_\beta) = \dim Y + \dim X_\beta - \dim X
< \dim Y + \dim X_\alpha - \dim X
$$
because of the transversality of the intersection of 
$Y$ with $X_\beta$. Therefore each irreducible component
of $Y \cap \Xbar_\alpha$ meets 
$Y \cap X_\alpha$, which completes the proof of the lemma. \qed
\enddemo

\demo{Proof of Theorem {\rm2.5}} The second equalities in 
(a)--(c) follow from Proposition 2.2, Theorem 2.3, and the well 
known fact for the closures of Bruhat cells,
$$
\overline{B^+ y B^+} =   
\bigsqcup \Sb w \in W \\ w \leq y \endSb
B^+ w B^+.
$$

The first equalities in (b) and 
(c) are obtained by applying 
Lemma 2.6 to the Bruhat decomposition 
$G = \bigsqcup_{w \in W} B^+ w B^+$ of the group $G$
and taking $Y = U^- \dot{z}$ and 
$Y= B^- z$, respectively. In both cases, 
the intersection of $Y$ with any Bruhat
cell $B^+ w B^+$ is transversal since 
$\Lie U^-$ and $\Lie B^+$ span $\Lie G$. Moreover,
in both cases $Y$ is a closed subvariety 
of $G$ and $\Cl_Y(Z)$ coincides with
$\overline{Z}$ for any subset $Z$ of $Y$.

The first equality in (a) cannot be proved in exactly 
the same way because 
$\Lie U^-_z$ and $\Lie B^+$ do not span 
$G$. We apply Lemma 2.6 to the stratification 
of the flag variety $G/B^+$ by Schubert cells
$B^+ . w B^+$, and take $Y=B^- . z B^+$. 
This gives us
$$
\Cl_{B^- . z B^+}
(B^- .z B^+ \cap B^+ .y B^+) = 
B^- .z B^+ \cap 
\overline{B^+ .y B^+}.
$$   
Applying the biregular isomorphisms (2.6) and (2.7),
one obtains
$$
\Cl_{U^-_z \dot{z}}(U^-_z \dot{z} \cap B^+ yB^+) =
U^-_z \dot{z} \cap \overline{B^+ yB^+}.
$$
Since $U^-_z .z$ is a closed subvariety 
of $G$ we can replace the left hand side 
with $\overline{U^-_z \dot{z} \cap B^+ yB^+}$.
This completes the proof of (a). \qed
\enddemo

\head 3. A first approach to $\Mmn$ through a Poisson 
structure on $\Gr(n,m+n)$ 
\endhead

Throughout this section, fix positive integers $m$ and $n$, with
$$N=m+n.$$
We derive a description of the orbits of symplectic leaves in
$\Mmn$ under a natural action of the maximal torus of $\GLN$, by
embedding $\Mmn$ in a Grassmannian Poisson homogeneous space,
$\Gr(n,N)$.

\definition{3.1\. Generalities on $\GLN$}
The Borel subgroups
of $\GLN$ consisting of upper and lower triangular matrices
will be respectively denoted by $B^+$ and $B^-$. Let
$U^\pm$ be their unipotent radicals. The 
maximal torus of $\GLN$ consisting of diagonal
matrices will be denoted by $T$. In situations where it is helpful
to indicate that we are working with subgroups of the $N\times N$
general linear group, we will label the above Borel and Cartan
subgroups of $\GLN$ as $B^\pm_N$ and $T_N$. However, we
reserve subscripts on $U^\pm$ for a different purpose -- see
(3.3) below.

We will need to describe a number of sets of matrices given in
block form, and it will be convenient to use a block form of set
notation for the purpose. For example, if $A$, $B$, $C$, $D$ are
subsets of
$M_n$, $M_{n,m}$, $M_{m,n}$, $M_{m}$ respectively, we set
$$\leftmatrix A&B\\ C&D \rightmatrix= \biggl\{ \leftmatrix a&b\\
c&d \rightmatrix \in M_N \biggm| a\in A,\, b\in B,\, c\in C,\, d\in
D \biggr\}.$$
In case one of the sets $A$, $B$, $C$, $D$ is a singleton, we may
omit the corresponding braces from the notation. Thus, for
instance, the notation $\leftmatrix I_n &M_{n,m}\\ 0&I_{m}
\rightmatrix$ indicates the unipotent subgroup $\biggl\{
\leftmatrix I_n &b\\ 0 &I_{m} \rightmatrix \biggm| b\in M_{n,m}
\biggr\}$ of $\GLN$, where
$I_n$ and $I_{m}$ are the identity matrices of sizes $n$ and
$m$.

Define the following maximal parabolic
subgroup  of $\GLN$:
$$
P_n= \leftmatrix GL_n &M_{n, m}\\ 0 &GL_{m} \rightmatrix.
\tag 3.1
$$
Let $L_n$ be the Levi factor of $P_n$ containing $T$,
and $U_n^+$ the unipotent radical of $P_n$.
Denote by $U_n^-$ the unipotent radical 
of the parabolic subgroup 
of $\GLN$ opposite to $P_n$. Explicitly, $L_n =L_n^1
L_{m}^2$ where
$$\xalignat2
L_n^1 &=
\leftmatrix GL_n&0 \\ 0&I_{m}\rightmatrix \cong GL_n 
&L_{m}^2 &=
\leftmatrix I_n&0 \\ 0&GL_{m} \rightmatrix \cong GL_{m}
\tag 3.2 \endxalignat 
$$
and
$$\xalignat2
U_n^+ &= \leftmatrix I_n &M_{n, m} \\ 0&I_{m}\rightmatrix
&U_n^- &= \leftmatrix I_n & 0 \\ M_{m, n} &I_{m}\rightmatrix.
\tag 3.3\endxalignat
$$
Let $\bfrak^\pm$, $\hfrak$, 
$\pfrak_n$, 
$\lfrak_n$, $\lfrak_n^1$, $\lfrak_{m}^2$, and 
$\nfrak_n^\pm$ denote the Lie algebras of 
$B^\pm$, $T$, $P_n$, $L_n$, $L_n^1$,
$L_{m}^2$, and $U_n^\pm$. 
The Lie algebras $\nfrak_n^+$ and 
$\nfrak^-_n$ are naturally identified as vector 
spaces with $M_{n, m}$ and $M_{m, n}$. 
The exponential maps 
$\exp : \nfrak_n^\pm \rightarrow U_n^\pm$ are
bijective and are explicitly given by
$$\xalignat2
\nfrak_n^- &\cong M_{m, n} \ni x \mapsto
\leftmatrix I_n & 0 \\ x &I_{m}\rightmatrix
&\nfrak_n^+ \cong M_{n, m} \ni y \mapsto
\leftmatrix I_n & y \\ 0 &I_{m}\rightmatrix.
\tag 3.4\endxalignat
$$

The Weyl group of $\GLN$ is isomorphic to 
the symmetric group $S_N$. The maximal length 
element of $S_N$ will be denoted 
by $\wN$. Explicitly, we have 
$\wN= \left(\smallmatrix 1&2&\cdots&N\\ N&N-1&\cdots&1
\endsmallmatrix\right)$.

{\it{For $k = 1, \ldots, N$,
we will denote by $S^1_k$ and $S^2_k$ the subgroups  
of $S_N$ that are isomorphic to $S_k$ and permute
respectively the first and the last $k$ indices.}}
In other words:
$$\aligned
S_k^1 &= 
\{w\in S_N \mid w(i)=i \text{\ for all\ } i>k\} \\
S_k^2 &= 
\{w\in S_N \mid w(i)=i \text{\ for all\ } i\le N-k\}.
\endaligned\tag 3.5
$$
In this notation, the Weyl groups of $L_n^1$ and 
$L_{m}^2$ are identified respectively 
with the subgroups $S_n^1$ and $S_{m}^2$
of the Weyl group $S_N$ of $\GLN$. 
The Weyl group of the Levi factor $L_n$ is 
identified with the subgroup 
$S_n^1  S_{m}^2$ of $S_N$.
 
Denote by $\wnwmpair \in S_N$ the product of the maximal 
length elements of $S_n^1$ and $S_{m}^2$. In other words,
this is the maximal length element of the Weyl group of the 
Levi factor $L_n$. Set
$$
\wmwn = w^N_\ci \wnwmpair.
\tag 3.6
$$
It is the maximal length representative in $S_N$ of the coset
$\wN (S_n^1  S_{m}^2)$. 

For a given $w \in S_N$, define the following
subsets of $S_N$:
$$\gather
\aligned S_N^{\leq w} &= \{ y \in S_N \mid y \leq w \} \\
S_N^{\ge w} &= \{ y \in S_N \mid y \ge w \} \endaligned \tag 3.7 \\
S_N^{[-n,m]} = \{ y \in S_N \mid -n\le s(i)-i\le m \text{\ for all\
} i=1, 2,
\ldots, N \}. \tag 3.8
\endgather
$$  
In Lemma 3.12, we will show that the subsets $S_N^{[-n,m]}$ and
$S_N^{\le \wmwn}$ of $S_N$ coincide. This set will  enter as a
parametrizing set for the set of $T$-orbits of symplectic  leaves
of the matrix affine Poisson space $M_{m, n}$.

Finally, consider the embedding 
$$\xalignat2
S_N &\hookrightarrow N(T),
 &S_N &\ni w \mapsto (a_{ij})= (\delta_{i w(j)}) \in N(T),
\tag 3.9\endxalignat
$$
which is a section for the projection 
$N(T) \rightarrow N(T)/T \cong S_N$.
{\it{By abuse of notation we will identify 
$S_N$ with its image in $N(T)$, and thus use the same letter $w$ to
denote the permutation matrix in $N(T)$ corresponding to a
permutation $w \in S_N$.}} 
Under this identification, the maximal 
length element $\wN \in S_N$ 
corresponds to the (unit) anti-diagonal 
matrix. Moreover, we have
$$
\wmwn = \leftmatrix 0 &\wm\\ \wn &0
\rightmatrix \leftmatrix \wn &0\\ 0 &\wm \rightmatrix =
\leftmatrix 0 &I_{m}\\ I_n &0 \rightmatrix .
$$
\enddefinition

\definition{3.2\. $GL_{N}/P_n$ and $\Gr(n,N)$} Recall the
natural isomorphism
$\Gr(n,N) \cong GL_{N}/P_n$.
\enddefinition

\proclaim{Proposition} {\rm(a)} The orthogonal complement
of $\pfrak_n$ in the dual Lie bialgebra $\gl_{N}^*$
for the standard Lie bialgebra structure on $\gl_{N}$ 
{\rm(}recall {\rm(1.5))} is 
$\pfrak_n^\perp = \nfrak_n^+ \oplus \{0\}$.

{\rm(b)} The parabolic subgroup $P_n$ of $GL_{N}$
is a Poisson algebraic subgroup for the standard 
Poisson structure on $GL_{N}$.

{\rm(c)} The pair 
$\bigl( \Gr(n,N) \cong GL_{N}/P_n,\, - \chi(r^{N}) \bigr)$ is
a  Poisson homogeneous space for the standard Poisson 
algebraic group $GL_{N}$. Here $r^{N}$ is the standard
$r$-matrix {\rm(1.3)} for $\gl_{N}$ and $\chi$ denotes the 
infinitesimal action for the left multiplication of 
$GL_{N}$ on $\Gr(n,N)$.

{\rm(d)} The Drinfeld Lagrangian subalgebra of the base point 
$e P_n$ of the Poisson homogeneous space $\bigl( GL_{N}/P_n,\,
-\chi(r^{N})
\bigr)$ is 
$$\aligned
\ol{\lfrak}_n &=
\biggl\{ 
\bigg( \leftmatrix a & b_1 \\ 0&c \rightmatrix, 
\leftmatrix a & b_2 \\ 0&c \rightmatrix \bigg) 
\biggm| a \in \gl_n,\, c \in \gl_m,\,
b_i \in M_{n, m} \biggr\} \\
 &\subseteq 
\gl_{N} \oplus \gl_{N}
\cong D(\gl_{N}).
\endaligned\tag 3.10
$$
It is the tangent Lie algebra of the algebraic subgroup
$$
\ol{L}_n =
\biggl\{ 
\bigg( \leftmatrix a & b_1 \\ 0&c \rightmatrix, 
\leftmatrix a & b_2 \\ 0&c \rightmatrix \bigg) 
\biggm| a \in GL_n,\, c \in GL_m,\,
b_i \in M_{n, m} \biggr\} 
\tag 3.11
$$
of $GL_{N} \times GL_{N}$; in particular,
$\bigl( \Gr(n, N),\, - \chi(r^{N}) \bigr)$ is an algebraic 
Poisson homogeneous space for the standard 
Poisson algebraic group $GL_{N}$. Moreover,
$$
\Delta(GL_{N}) \cap \ol{L}_n = 
\Delta(P_n).
\tag 3.12
$$ 

{\rm(e)} Each intersection of a $B^+$- and a $B^-$-orbit 
on $\Gr(n,N)$ is a locally closed Poisson 
subvariety of $\bigl( \Gr(n, N),\, - \chi(r^{N}) \bigr)$.
\endproclaim

\demo{Proof} (a) It is straightforward to check that 
$\nfrak^+_n \oplus \{0\} \subseteq \gl^*_{N} \subseteq \gl_{N}
\oplus
\gl_{N}$ is orthogonal to $\Delta(\pfrak_n)$ with respect to the 
bilinear form (1.6), recall \S1.4. 
The statement now follows from the fact that the sum of the dimensions 
of $\pfrak_n$ and $\nfrak_n^+$ is 
equal to
$\dim \gl_{N}$.

Part (b) follows from Theorem 1.8(d) and the first part. 

(c) Consider the 
projection 
$p : GL_{N} \rightarrow \GLN/P_n$ and the Poisson
structure (1.10) for $\GLN/P_n$. Since the standard matrices 
$E_{ij}$ belong to $\pfrak_n$
for $i<j$, we have 
$$
p_*(\chi^L(r^{N})) =0.
$$
Thus in the present situation 
the Poisson structure (1.10) is exactly 
$- \chi(r^{N})$. Now part (c) follows from the 
discussion before Theorem 1.8.

(d) Since the Poisson structure $-\chi(r^{N})$ 
vanishes at the base point $e P_n$ of 
$GL_{N}/P_n$, according to Theorem 1.8(b) 
the  Drinfeld Lagrangian subalgebra
of the double 
$D(\gl_n) \cong \gl_n \oplus \gl_n$ is 
$\Delta(\pfrak_n) + \pfrak_n^\perp$. A simple 
computation leads to (3.10). 
The rest of part (d) is straightforward and 
will be omitted.

(e) Observe that the subgroup $T\subseteq GL_{N}\subseteq
D(GL_{N})$ normalizes the subgroup $F\subseteq D(GL_{N})$ (recall
(1.4)), and that $TF= B^+\times B^-$.
Theorem 1.10 implies that the $T$-orbits of 
symplectic leaves of
$\bigl( GL_{n+ m}/P_n,\, - \chi(r^{N}) \bigr)$
are the irreducible components of the inverse images of the 
$(B^+ \times B^-)$-orbits on 
$D(\GLN)/\ol{L}_n$ under the  
map 
$$
GL_{N}/P_n @>{\Delta}>> 
D(\GLN)/\ol{L}_n
$$
(cf.~(1.12)), which is an embedding because 
of (3.12). It is obvious (because $\ol{L}_n \subseteq P_n\times P_n$)
that each such  inverse image falls within a single 
intersection of a $B^+$- and $B^-$-orbit
on $\Gr(n,N)$. Thus, the latter are finite unions of $T$-orbits
of symplectic leaves and hence Poisson subvarieties of 
$\bigl( \Gr(n, N),\, - \chi(r^{N}) \bigr)$.
\qed\enddemo

Throughout the remainder of the section, we shall always assume
that $\Gr(n,N)\cong GL_{N}/P_n$ has been equipped with the
Poisson structure $-\chi(r^{N})$.

\definition{3.3\. The open $B^-$-orbit on $\Gr(n,N)$}
The $B^-$-orbit through the base point of $GL_{N}/P_n$
is a Zariski open subvariety. According to Proposition 3.2(e),
it is a Poisson subvariety of 
$GL_{N}/P_n$. Moreover,
the open orbit $B^-. P_n \subseteq GL_{N}/P_n$
is an affine space which 
is isomorphic to 
$U_n^-$ by 
$$
U^-_n \ni u \mapsto u P_n;
$$
in particular, $B^-.P_n= U^-_n.P_n$.
Composing this map with the exponential
map 
$$
\exp : \Mmn \cong \nfrak^-_n 
@>{\cong}>> U^-_n
$$
induces an isomorphism of affine spaces
$$\xalignat2
\Mmn &@>{\cong}>>  U^-_n.P_n \subseteq GL_{N}/P_n, &x &\mapsto 
\leftmatrix I_n & 0 \\ x & I_m \rightmatrix P_n.
\tag 3.13\endxalignat
$$
We consider a twisted version of this 
isomorphism:
$$\xalignat2
\Psi : &\Mmn @>{\cong}>> U^-_n . P_n,
&x &\mapsto \exp(x \wn) P_n = 
\leftmatrix I_n & 0 \\ x \wn & I_m \rightmatrix P_n.
\tag 3.14\endxalignat
$$
Recall that $\wn$ denotes the maximal length element of 
$S_n$ and its representative in the normalizer of 
the diagonal subgroup of $GL_n$, as fixed in \S3.1. 
\enddefinition 

The restriction of the Poisson structure $-\chi(r^{N})$
to $U^-_n. P_n$ was computed by
Gekhtman, Shapiro, and Vainshtein
in \cite{\GSV}. The following result can be deduced from 
their computations, but
we offer a more geometric proof.

\proclaim{3.4\. Proposition} The map 
$\Psi : \Mmn \rightarrow U^-_n . P_n$ is a Poisson
isomorphism between  the matrix affine Poisson space $\Mmn$
and the Poisson subvariety $U^-_n . P_n$ of $GL_{N}/P_n$. 
\endproclaim

We break the standard $r$-matrix for $\gl_{N}$ into 
three terms as follows
$$
r^{N} = 
\sum_{1 \leq i < j \leq n} E_{ij} \wedge E_{ji} +
\sum_{n< i < j \leq N} E_{ij} \wedge E_{ji} +
\sum_{i\leq n <j} E_{ij} \wedge E_{ji}
\in \wedge^2 \gl_{N}
\tag 3.15
$$
and denote them by $r^{N}_1$, $r^{N}_2$, and $r^{N}_3$,
respectively. First we establish an auxiliary result.

\proclaim{3.5\. Lemma} In the above notation,
$$
\chi(r^{N}_3)|_{U^-_n . P_n}=0.
$$
\endproclaim

\demo{Proof} We shall use the label (3.13)$^{-1}$ to refer to the
inverse isomorphism $U_n^-.P_n \rightarrow \Mmn$ of the isomorphism
(3.13).
Since $U^-_n$ is abelian and 
$E_{j+n, i} \in \nfrak^-_n$ for 
$i\le n$, $j \leq m$, under the isomorphism (3.13)$^{-1}$ we have
$$
\chi(E_{j+n, i})|_{U_n^-.P_n} \mapsto \frac{\partial}{\partial
x_{ji}} 
\quad \text{for} \quad 
i\le n,\, j \leq m.
\tag 3.16
$$
By a direct computation, one checks that 
for $x \in \Mmn$, $y \in M_{n,m}$, and a small 
$\ep \in \CC$,
$$
\leftmatrix I_n & \ep y \\ 0 &I_m\rightmatrix .
\leftmatrix I_n & 0 \\ x &I_m\rightmatrix P_n 
= \leftmatrix I_n & 0 \\ x( I_n + \ep y x)^{-1} &I_m\rightmatrix P_n
= \leftmatrix I_n & 0 \\ x - \ep x y x + O(\ep^2) &I_m\rightmatrix P_n.
$$
This implies that under the isomorphism (3.13)$^{-1}$,
$$
\chi(E_{i, j+n})|_{U_n^-.P_n} \mapsto - \sum_{k=1}^m \sum_{l=1}^n
x_{ki} x_{jl}
\frac{\partial}{\partial x_{kl}} 
\quad \text{for} \quad 
i\le n,\, j \leq m.
\tag 3.17
$$
Combining (3.16) and (3.17), we see that under the isomorphism
(3.13)$^{-1}$,
$$
\chi(r^{N}_3)|_{U_n^-.P_n} \mapsto 
- \sum_{k,j=1}^m \sum_{i, l=1}^n
x_{ki} x_{jl}
\frac{\partial}{\partial x_{kl}} 
\wedge
\frac{\partial}{\partial x_{ji}} 
=0. \qquad\square
$$
\enddemo

\definition{3.6\. Actions of $L_n$ on $U^-_n. P_n$ and
$\Mmn$, and a proof of Proposition 3.4} 
Since the Levi factor $L_n$ normalizes $U^-_n$, it preserves the
open $B^-$-orbit $U^-_n. P_n$ on $\Gr(n,N)$ (recall \S3.1 
for notation). Via the isomorphism 
(3.13), this induces an action of 
$GL_m \times GL_n \cong L^2_n \times L^1_n=L_n$ on the affine space
$\Mmn$. It is given by $(a, b). x = ax b^{-1}$ for
$a \in GL_m$, $b \in GL_n$, $x \in \Mmn$,
which is checked 
by a direct computation:  
$$
\leftmatrix b & 0 \\ 0 & a\rightmatrix .
\leftmatrix I_n & 0 \\ x & I_m\rightmatrix P_n =
\leftmatrix I_n & 0 \\ a x b^{-1} & I_m\rightmatrix P_n.
$$
This action of $GL_m \times GL_n$ on $\Mmn$ 
breaks into the actions of $GL_m$ and $GL_n$ on $\Mmn$ from 
\S1.5, used to define the standard Poisson 
structure $\pi_{m,n}$ on $\Mmn$. In \S5.2, we will 
consider this from a Poisson point of view.
\enddefinition

\demo{In the rest of {\rm \S3.6} we prove Proposition {\rm 3.4}}
The terms $r^{N}_1$ and $r^{N}_2$ 
of the standard $r$-matrix on $\gl_{N}$, see (3.15), are 
respectively equal to the pushforwards of $r^n$ and $r^m$ 
under $\gl_n \cong \lfrak^1_n \hookrightarrow \gl_{N}$
and $\gl_m \cong \lfrak^2_m \hookrightarrow \gl_{N}$. From the above 
discussion it follows that
under the isomorphism (3.13)$^{-1}$,
$$
-\chi(r^{N}_1 + r^{N}_2)|_{U^-_n. P_n} \mapsto
- \chi^L(r^m) - \chi^R(r^n).
$$  
(Recall from \S1.5 that
$\chi^L(.)$ and $\chi^R(.)$ denote 
the infinitesimal actions 
of $\gl_m$ and $\gl_n$ on $\Mmn$.)
Since the maximal length element $\wn \in S_n$ 
satisfies $\Ad_{\wn}(E_{ij}) = E_{n+1-i,n+1-j}$, we have 
$\Ad_{\wn}(r^n) = -r^n$, and thus
$$
\Psi_*(\pi^{m,n}) = 
-\chi(r^{N}_1 + r^{N}_2)|_{U^-_n. P_n}
$$
(see (1.8) and (3.14) for the definitions of $\pi^{m,n}$ and
$\Psi$). Now Proposition 3.4 follows from Lemma 3.5.
\qed
\enddemo

\definition{3.7\. A Poisson homogeneous space of $B^-$}
One can use Proposition 3.4 to identify $\Mmn$ with a (full) Poisson
homogeneous space of $B^-$. First, recall the well known fact that $B^-$
is a Poisson algebraic subgroup of $\GLN$. Since $P_n$ is also a Poisson
algebraic subgroup of $\GLN$ (cf. Proposition 3.2 (b)), we get that
$$
B^- \cap L_n = B^- \cap P_n
$$
is a Poisson algebraic subgroup of $\GLN$ (and thus
of $(B^-, \pi^N|_{B^-})$ as well).
According to Theorem 1.8, one obtains a natural structure
of a Poisson homogeneous space for $(B^-, \pi^N|_{B^-})$ on
$B^-/(B^-\cap L_n)$ by equipping it with the Poisson bivectorfield
$\nu_*(\pi^N|_{B^-})$ where $\nu$ is the projection
$\nu : B^- \rightarrow B^-/(B^- \cap L_n)$.
\enddefinition

\proclaim{Corollary} The map
$$
\Mmn \cong \nfrak^-_n \ni x \mapsto \exp(x \wn) (B^- \cap L_n)
$$
is a Poisson isomorphism between the matrix affine Poisson
space $\Mmn$ and the Poisson homogeneous space
$(B^-/(B^-\cap L_n), \nu_*(\pi^N|_{B^-}) )$ of $(B^-, \pi^N|_{B^-})$.
\endproclaim

One can use this corollary instead of Proposition 3.4 in obtaining
the results in \S 3.8, but Proposition 3.4 is conceptually more
important since it provides a natural compactification of the matrix
affine Poisson space.

\demo{Proof} The map $\Psi$ provides a Poisson isomorphism of $\Mmn$
with the complete Poisson subvariety $U_n^-. P_n$ of
$(\GLN/P_n, -\chi(r^N)),$ cf\. Proposition 3.4. The latter is
a $B^-$-orbit with stabilizer $B^- \cap L_n$ of the base point $P_n$
and thus can be identified
with the homogeneous space $B^-/(B^-\cap L_n)$. Under this
identification, the Poisson structure $-\chi(r^N)|_{U_n^-.P_n}$ is
matched with the Poisson structure $\nu_*(\pi^N|_{B^-})$ because both
are pushforwards of the standard Poisson structure $\pi^N$ on $\GLN$.
The corollary now follows from the fact that
$x \mapsto \exp(x \wn) (B^- \cap L_n)$
is nothing but the map $\Psi$ when we identify $U_n^-.P_n$ and
$B^-/(B^-\cap L_n)$.
\qed\enddemo

\definition{3.8} Recall the notation $W^V$ for the set of minimal
length representatives for left cosets of a subgroup $V$ of a Weyl
group $W$.
\enddefinition

\proclaim{Lemma} The set $S_{N}\times S_{N}^{S^1_n S^2_m}$ is a
complete, irredundant set of representatives for the $(B^+\times
B^-, \ol{L}_n)$ double cosets in $GL_{N}\times GL_{N}$.
\endproclaim

\demo{Proof} We apply Theorem A.1. For that purpose, let $G=
GL_{N}\times GL_{N}$, choose $B^+\times B^-$ and $B^-\times
B^+$ to be the positive and negative Borel subgroups of $G$,
respectively, and consider the parabolic subgroup $P= P_n\times
P_n$ of $G$, which contains $B^+\times B^+$. There is a Levi
decomposition $P= L_0N$ where $L_0= L_n\times L_n \supset T\times
T$ and $N= U^+_n\times U^+_n$, and we put $L_0= L_n^{\ell} L_n^r$
where $L_n^{\ell}= L_n\times \{I\}$ and $L_n^r= \{I\}\times L_n$.
There is an isomorphism $\Theta: L_n^{\ell} \rightarrow L_n^r$
given by $\Theta(a,I)= (I,a)$, and we observe that the simple factors
$F \times \{I\}$ of $L_n^{\ell}$ (where $F= L_n^1$, $L_m^2$) satisfy
$$
\Theta \bigl( (F\times \{I\}) \cap (B^-\times B^+) \bigr)= (\{I\}\times
F) \cap (B^+\times B^-).
$$
Let $\pi_j: P\rightarrow P/N \cong L_0
\rightarrow L_n^j$ (for $j=\ell,r$) denote the natural projections,
and observe that the subgroup
$$
R= \{ p\in P \mid \Theta\pi_{\ell}(p)= \pi_r(p) \}
$$
coincides with $\ol{L}_n$. Since the Weyl group of $L_n^r$,
considered as a subgroup of the Weyl group of $G$, is just
$\{1\}\times (S^1_n S^2_m) \subseteq S_{N}\times S_{N}$, Theorem
A.1 implies that the set 
$$
(S_{N}\times S_{N})^{\{1\}\times
(S^1_n S^2_m)}= S_{N}\times S_{N}^{S^1_n S^2_m}
$$
is a complete,
irredundant set of representatives for the $(B^+\times
B^-, \ol{L}_n)$ double cosets in $G$.
\qed\enddemo

\definition{3.9\. $T$-orbits of symplectic leaves in $\Mmn$} Since
the image of $GL_m\times GL_n\cong L_n \subseteq GL_{N}$ contains
the torus $T$, the action of $GL_m\times GL_n$ on $\Mmn$ given in
\S3.6 incorporates an action of $T$ on $\Mmn$. Specifically,
if $T_m$ and $T_n$ denote the maximal tori consisting of diagonal
matrices in $GL_m$ and $GL_n$ respectively, then $(a,b).x=
axb^{-1}$ for $a\in T_m$, $b\in T_n$, $x\in \Mmn$.
\enddefinition

\proclaim{Theorem} There are only finitely many 
$T$-orbits of symplectic leaves on the matrix affine Poisson 
space $\Mmn$. They are smooth irreducible locally closed subvarieties of
$\Mmn$, and they are parametrized by
$S_{N}^{\ge
\wnwmpair}$,  recall {\rm (3.7)}. The $T$-orbit of symplectic
leaves corresponding to $w \in S_{N}^{\ge \wnwmpair}$ is   
explicitly given by 
$$
\P_w = \biggl\{
x\in\Mmn \biggm| \leftmatrix \wn &0\\ x &\wm \rightmatrix \in 
B^+wB^+ \biggr\}.
\tag 3.18
$$
As an algebraic variety, $\P_w$ is biregularly isomorphic
to $\B_{\wnwmpair, w}$.
\endproclaim

\demo{Proof}
We will make use of the isomorphism $\Psi$ (see (3.14)) of
Proposition 3.4 between the matrix affine Poisson space $\Mmn$ and 
the $T$-stable Poisson subvariety
$U^-_n. P_n$ of $GL_{N}/P_n$.
Recall that $U_n^-.P_n= B^-.P_n$ is open in $GL_{N}/P_n$.
The isomorphism $\Psi$ is not
$T$-equivariant, but we have
$$\align
\Psi((a,b).x) &= \leftmatrix I_n &0\\ axb^{-1}\wn &I_m \rightmatrix
P_n \\
 &= \leftmatrix I_n &0\\ ax\wn(\wn b\wn)^{-1} &I_m \rightmatrix
P_n= \leftmatrix \wn b\wn &0\\ 0 &a \rightmatrix . \Psi(x)
\endalign
$$
for $a\in T_m$, $b\in T_n$, $x\in \Mmn$, whence $\Psi$ and
$\Psi^{-1}$ preserve $T$-orbits. Consequently, $\Psi$ maps
$T$-orbits of symplectic leaves in $\Mmn$ to $T$-orbits of
symplectic leaves in $U_n^-.P_n$.

Firstly, Theorem 1.10 (applied with $H=T$, as in the proof of
Proposition 3.2(e)) implies that the
$T$-orbits of symplectic leaves of $U^-_n. P_n$
are smooth locally closed subvarieties, and so the same is true for
$\Mmn$. The map (1.12) in the present situation  is 
$$\xalignat2
\Delta : GL_{N}/P_n &\hookrightarrow 
(GL_{N} \times GL_{N})/ \ol{L}_n,
 &\Delta(g P_n) &= (g,g) \ol{L}_n
\endxalignat 
$$
(see the proof of Proposition 3.2(e)). From Theorem 1.10, we also
know that the $T$-orbits of symplectic leaves of $U^-_n. P_n$ are
those irreducible components of inverse images of $(B^+ \times 
B^-)$-orbits on $(GL_{N} \times GL_{N})/
\ol{L}_n$  under $\Delta$ that lie inside $U^-_n. P_n$.

The set of $(B^+ \times B^-)$-orbits on 
$(GL_{N} \times GL_{N})/ \ol{L}_n$ is in one to one correspondence with 
the set of $(B^+ \times B^-, \ol{L}_n)$ double cosets in 
$GL_{N} \times GL_{N}$. According to Lemma 3.8, the latter set
is parametrized by 
$S_{N} \times S_{N}^{S^1_n S^2_m}$. Therefore, the 
$(B^+ \times B^-)$-orbits on $(GL_{N} \times GL_{N})/ \ol{L}_n$
are the sets
$$\xalignat2
(B^+ \times B^-) . (w_1,w_2)\ol{L}_n, & 
&w_1 &\in S_{N},\
w_2 \in S_{N}^{S^1_n S^2_m},
\tag 3.19\endxalignat
$$ 
and all such sets are distinct. Observe that 
$$
\Delta^{-1} \left( 
(B^+ \times B^-) . (w_1,w_2)\ol{L}_n \right) \subseteq
B^- . w_2 P_n.
$$
If $w_2 \in S_{N}^{S^1_n S^2_m}$ and $w_2 \neq 1$,
then $B^- w_2 P_n \cap B^- P_n = \emptyset$ because of the Bruhat 
lemma. Thus, only the $\Delta$-inverse images of the sets
(3.19) with $w_2 = 1$ might intersect $U^-_n . P_n$ 
nontrivially. 

The intersection with $U^-_n . P_n$ of the 
$\Delta$-inverse image of the set (3.19) 
with $w_2=1$ 
consists of $u P_n \in GL_{N}/P_n$
for those $u \in U^-_n$ for which 
$$ 
u = b^+ w_1 l u^+_1 = b^- l u^+_2
\tag 3.20
$$
for some $b^\pm \in B^\pm$, $l \in L_n$,
$u^+_i \in U_n^+$. From these equalities, 
one obtains $l \in L_n \cap B^-$ 
and $u_2^+=e$. Conversely, if $u = b^+ w_1 l u^+_1$ for some
$b^+\in B^+$, $l\in L_n\cap B^-$, $u_1^+\in U_n^+$, we can also
write $u= b^-l$ where $b^-= ul^{-1}\in B^-$. Thus,
$$
\Delta^{-1} \left( (B^+ \times B^-) . (w_1,1)\ol{L}_n \right)
\cap U^-_n . P_n =
\left( U^-_n \cap B^+ w_1 (L_n\cap B^-) U^+_n \right). P_n.
$$
Next, observe that $(L_n\cap B^-)U_n^+= \wnwmpair B^+\wnwmpair$
and that $U^-_n = U^-_{\wnwmpair}$ (recall (2.2)). Thus, setting
$w= w_1\wnwmpair$ and recalling the notation (2.11), we have
$$\aligned
\Delta^{-1} \left( (B^+ \times B^-) . (w_1,1)\ol{L}_n \right)
\cap U^-_n . P_n &=
\left( U^-_n \cap B^+ w_1 (L\cap B^-) U^+_n \right). P_n  \\
&= \left( U^-_{\wnwmpair} \cap B^+ w B^+ \wnwmpair \right). P_n \\
&= \U_{\wnwmpair, w}. P_n
\endaligned \tag 3.21
$$
(since $\wnwmpair\in P_n$).
According to Theorem 2.4, $\U_{\wnwmpair, w}$ is 
irreducible. Therefore, the set (3.21) is a single $T$-orbit 
of symplectic leaves of $U_n^-.P_n$. The fact that the $T$-orbits 
of symplectic leaves of the matrix affine
Poisson space are the sets (3.18) follows 
by applying the Poisson isomorphism 
$\Psi : \Mmn \rightarrow U^-_n. P_n$
to (3.21). Namely, since $U_n^-\cap P_n= \{I\}$, we compute that
$$\align
\Psi^{-1} \bigl( \U_{\wnwmpair, w}. P_n \bigr) &= \biggl\{ x\in\Mmn
\biggm| \leftmatrix I_n &0\\ x\wn &I_m \rightmatrix P_n \in \bigl(
U_n^-\cap B^+wB^+\wnwmpair \bigr) .P_n \biggr\} \\
 &= \biggl\{ x\in\Mmn
\biggm| \leftmatrix I_n &0\\ x\wn &I_m \rightmatrix \in U_n^-\cap
B^+wB^+\wnwmpair \biggr\} \tag 3.22 \\
 &= \biggl\{ x\in\Mmn
\biggm| \leftmatrix I_n &0\\ x\wn &I_m \rightmatrix \leftmatrix \wn
&0\\ 0 &\wm \rightmatrix \in B^+wB^+ \biggr\} = \P_w.
\endalign
$$
Moreover, $\P_w\cong \U_{\wnwmpair, w}\cong \B_{\wnwmpair, w}$ by
Theorem 2.4. Irreducibility thus follows from Proposition
2.2 of Deodhar. Finally,
$\P_w$ is nonempty if and only if
$\B_{\wnwmpair, w}$ is nonempty, which occurs precisely when $w\ge
\wnwmpair$, by Proposition 2.2.  
\qed\enddemo

\definition{3.10} Since $\leftmatrix \wn &0\\ \Mmn &\wm \rightmatrix
\subseteq B^-\wnwmpair B^-$, the set $\P_w$ described in (3.18) can
be written as the inverse image of $B^-\wnwmpair B^-\cap B^+wB^+$
under the map $\Omega: \Mmn\rightarrow GL_{N}$ given by $x\mapsto
\leftmatrix \wn &0\\ x &\wm \rightmatrix$. It is known that the
$T$-orbits of symplectic leaves in $GL_{N}$ coincide with the
double Bruhat cells $B^-yB^-\cap B^+wB^+$ (e.g., this follows from
the results of
\cite{\HLthree, Appendix A}). The following statement is thus an immediate
consequence: {\it The
$T$-orbits of symplectic leaves in
$\Mmn$ are precisely the nonempty $\Omega$-inverse images of the $T$-orbits
of symplectic leaves in $GL_{N}$\/}. 
The lifting $\Omega$ of $\Psi$ is neither $T$-equivariant nor
Poisson, and because of this one cannot approach Theorem 3.9
directly using
$\Omega$.
\enddefinition

\definition{3.11. Alternative descriptions of
$S_{N}^{\ge\wnwmpair}$} It is convenient to describe the Bruhat
order on $S_N$ in terms of relations between sets of integers, as
follows. First, if $I$ and $J$ are $t$-element subsets of
$\{1,\dots,N\}$, list their elements in ascending order, say
$$\xalignat2
I &= \{i_1< i_2< \cdots< i_t\} &J &= \{j_1< j_2< \cdots< j_t\},
\endxalignat
$$
and then define $I\le J$ if and only if $i_l\le j_l$ for
$l=1,\dots,t$. For $y,z\in S_N$, we have
$$
y\le z \quad\iff\quad \{y(1),\dots,y(p)\} \le \{z(1),\dots,z(p)\} 
\text{\ for\ } p=1,\dots,N \tag 3.23
$$
(e.g., \cite{\Fulbook, Exercise 8, p\. 175}). For $I$ and $J$ as above,
it is clear that $I\le J$ if and only if $\wN(I)\ge \wN(J)$. Hence,
$$
y\le z \quad\iff\quad \wN y\ge \wN z
$$
for any $y,z\in S_N$.

In particular, a permutation $w\in S_{N}$ satisfies $w\ge
\wnwmpair$ if and only if $w^{N}_\ci w\le w^{N}_\ci \wnwmpair =
\wmwn$ (recall (3.6)). Thus,
$$
S_{N}^{\ge\wnwmpair} = w^{N}_\ci S_{N}^{\le\wmwn}. \tag 3.24
$$

The following description of $S_{N}^{\le\wmwn}$ is known in the
case $m=n$; we thank Jon McCammond for bringing the
result to our attention. This result also appears in \cite{\Lau,
Proposition 1.3}; we provide a proof for the reader's convenience.
Recall (3.8) for the notation
$S_{N}^{[-n,m]}$.
\enddefinition

\proclaim{3.12\. Lemma} $S_{N}^{\le\wmwn} = S_{N}^{[-n,m]}$ and
$$
S_{N}^{\ge \wnwmpair} = \{y\in S_{N} \mid n \le y(i)+i-1 \le
m+2n \text{\ for all\ } i=1,\dots,N\}.
$$
\endproclaim

\demo{Proof} Since the second statement follows immediately from
the first via (3.24), we need only prove the first statement.

First, consider $s\in S_{N}^{\le\wmwn}$ and $j\in
\{1,\dots,N\}$. If
$j\le n$, then
$$
s(\{1,\dots,j\}) \le \wmwn(\{1,\dots,j\}) = \{m+1,\dots,m+j\}, \tag
3.25
$$
whence $s(j)\le m+j$.
On the other hand, if $j>n$, then 
$$
s(\{1,\dots,j-1\}) \le \wmwn(\{1,\dots,j-1\}) = 
\{1,\dots,j-1-n,m+1,\dots,N\},
$$
    from which we see that
$\{1,\dots,j-1-n\} \subseteq s(\{1,\dots,j-1\})$, and consequently
$s(j)\ge j-n$. We automatically have $s(j)\ge j-n$ when $j\le n$,
and $s(j)\le m+j$ when $j>n$. Thus, $s\in S_{N}^{[-n,m]}$.

Conversely, let $s\in S_{N}^{[-n,m]}$ and $j\in
\{1,\dots,N\}$. If $j\le n$, then $s(i)\le
i+m\le j+m$ for $i\le j$, whence $s(\{1,\dots,j\}) \subseteq
\{1,\dots,j+m\}$, and consequently (3.25) holds. On the other hand,
if $j>n$, then $s(i)\ge i-n> j-n$ for $i>j$, whence $\{1,\dots,j-n\}
\subseteq s(\{1,\dots,j\})$, and consequently
$$
s(\{1,\dots,j\}) \le \{1,\dots,j-n,m+1,\dots,N\} =
\wmwn(\{1,\dots,j\}).
$$
Therefore $s\le \wmwn$. 
\qed\enddemo

In the last result of this section, we describe the 
Zariski closures of the $T$-orbits of symplectic leaves
of the matrix affine Poisson space.
 
\proclaim{3.13\. Theorem} The Zariski closures of the 
$T$-orbits of symplectic leaves of the matrix affine Poisson space
$\Mmn$, see Theorem {\rm 3.9}, are given by
$$
\ol{\P}_w = \bigsqcup \Sb z\in S_{N}\\ \wnwmpair\le z\le w \endSb
\P_z = \biggl\{
x\in\Mmn \biggm| \leftmatrix \wn &0\\ x &\wm \rightmatrix \in 
\ol{B^+wB^+} \biggr\}. \tag 3.26
$$
Consequently, the inclusions between the Zariski closures
of the $T$-orbits of symplectic leaves {\rm (3.18)} on $\Mmn$
correspond to the Bruhat order on 
$S_{N}^{\ge\wnwmpair}$. 
\endproclaim

\demo{Proof} As noted in the proof of Theorem 3.9, $U_n^-=
U^-_{\wnwmpair}$. Since
$\wnwmpair
\in P_n$, the isomorphism between
$U_n^-$ and $U_n^-.P_n$ (recall \S3.3) yields a corresponding
isomorphism between $U_n^- \wnwmpair$ and $U_n^-.P_n$.

Now let $w\in S_{N}^{\ge\wnwmpair}$. According to (3.22), we have
$$
\P_w= \Psi^{-1}( \U_{\wnwmpair,w}.P_n ).
$$
Invoking the isomorphisms $\Psi: \Mmn \rightarrow U^-_n.P_n$ and $U_n^-
\wnwmpair
\rightarrow U_n^-.P_n$, we obtain
$$\aligned
\ol{\P}_w &= \Psi^{-1} \bigl( \Cl_{U_n^-.P_n}(
\U_{\wnwmpair,w}.P_n) \bigr) \\
 &= \Psi^{-1} \bigl( \bigl[ \Cl_{U_n^- \wnwmpair}
(\U_{\wnwmpair,w}) \bigr] .P_n \bigr) \\
 &= \Psi^{-1} \bigl( \ol{\U}_{\wnwmpair,w} .P_n \bigr).
\endaligned \tag 3.27
$$
By Theorem 2.5(a),
$$
\ol{\U}_{\wnwmpair,w}= \bigsqcup \Sb z\in S_{N}\\ \wnwmpair\le
z\le w \endSb \U_{\wnwmpair,z}. \tag 3.28
$$
The first equality of (3.26) follows from (3.27) and (3.28). Since
$\ol{B^+wB^+}$ is the disjoint union of the cells $B^+zB^+$ for
$z\le w$, we have
$$
\biggl\{
x\in\Mmn \biggm| \leftmatrix \wn &0\\ x &\wm \rightmatrix \in 
\ol{B^+wB^+} \biggr\} = \bigsqcup \Sb z\in S_{N}\\ z\le w \endSb
\P_z
\,,
$$
which yields the second equality of (3.26) because $\P_z$ is empty
when $z\not\ge \wnwmpair$ (recall the end of the proof of Theorem
3.9).
\qed\enddemo

\head 4. Computational description of $T$-orbits of symplectic
leaves \endhead

As in the previous section, we fix positive integers $m$, $n$, and
$N=m+n$. We derive a description of the $T$-orbits $\P_w$ of symplectic
leaves in $\Mmn$ in terms of ranks of rectangular submatrices.

\definition{4.1\. Descriptions of $B^+wB^+$ and $B^-wB^-$} In order
to give computational descriptions of the sets $\P_w$ in (3.18) and
$\ol{\P}_w$ in (3.26), we rely on the computational descriptions
of $B^-wB^+$ and its closure given by Fulton in \cite{\Ful}; these
descriptions are easily modified to deal with $B^+wB^+$. Since we
will also make use of the corresponding descriptions in $\Mmn$ and
$M_{n,m}$, we give a general version of these results.

 Let $1\le a\le b\le k$ and $1\le c\le d\le l$.
For $x\in M_{k,l}$, we write $x_{[a,\dots,b;c,\dots,d]}$ to denote
the submatrix of $x$ with rows $a,\dots,b$ and columns
$c,\dots,d$. Recall from \S3.1 that we use $B^\pm_k$ and
$B^\pm_l$ to denote the standard Borel subgroups of $GL_k$ and
$GL_l$. The closures in the proposition below denote Zariski
closures in the matrix variety $M_{k,l}$.
\enddefinition

\proclaim{Proposition} {\rm [Fulton]} Let $k$ and $l$ be positive
integers and
$x,w\in M_{k,l}$.

{\rm (a)} $x\in B_k^+wB_l^+$ if and only if $\rank(
x_{[p,\dots,k;1,\dots,q]})= \rank( w_{[p,\dots,k;1,\dots,q]})$ for
all $p=1,\dots,k$ and $q=1,\dots,l$.

{\rm (b)} $x\in \overline{B_k^+wB_l^+}$ if and only if $\rank(
x_{[p,\dots,k;1,\dots,q]}) \le \rank( w_{[p,\dots,k;1,\dots,q]})$
for all $p=1,\dots,k$ and $q=1,\dots,l$.

{\rm (c)} $x\in B_k^-wB_l^-$ if and only if $\rank(
x_{[1,\dots,p;q,\dots,l]})= \rank( w_{[1,\dots,p;q,\dots,l]})$ for
all $p=1,\dots,k$ and $q=1,\dots,l$.

{\rm (d)} $x\in \overline{B_k^-wB_l^-}$ if and only if $\rank(
x_{[1,\dots,p;q,\dots,l]}) \le \rank( w_{[1,\dots,p;q,\dots,l]})$
for all $p=1,\dots,k$ and $q=1,\dots,l$.
\endproclaim

\demo{Proof} (a) Observe that $x\in B_k^+wB_l^+$ if and only if
$w^k_\ci x\in B_k^-w^k_\ci wB_l^+$. The result of \cite{\Ful, p\.
390, second display} shows that $w^k_\ci
x\in B_k^-w^k_\ci wB_l^+$ if and only if 
$$\rank \bigl( (w^k_\ci
x)_{[1,\dots,p;1,\dots,q]} \bigr) = \rank \bigl( (w^k_\ci
w)_{[1,\dots,p;1,\dots,q]} \bigr)$$
for $p=1,\dots,k$ and $q=1,\dots,l$. Part (a)
follows.

(b) This follows from \cite{\Ful, Proposition 3.3(a)} in the same
manner as (a). 

(c) and (d) follow similarly.
\qed\enddemo

\definition{4.2\. Description of $\P_w$} 
Recall the notation $\P_w$ from (3.18) for
$T$-orbits of symplectic leaves in $\Mmn$.

It will be convenient to write some matrices $w\in M_N$ in the
following block form:
$$\xalignat2 
w &= \left[\matrix w_{11} &w_{12}\\ w_{21} &w_{22}
\endmatrix\right],
 &&\pmatrix w_{11} \in M_n &w_{12}\in M_{n,m} \\
w_{21}\in \Mmn &w_{22}\in M_m \endpmatrix.
\tag 4.1 \endxalignat
$$
\enddefinition

\proclaim{Theorem} Let $x\in\Mmn$ and $w\in
S_{N}^{\ge\wnwmpair}$, and write $w=
\leftmatrix w_{11} &w_{12}\\ w_{21} &w_{22} \rightmatrix$ as in
{\rm (4.1)}. Then $x\in\P_w$ if and only if the following four
conditions hold:

{\rm (a)} $\rank(x_{[p,\dots,m;1,\dots,q]})=
\rank\bigl( (w_{21})_{[p,\dots,m;1,\dots,q]} \bigr)$ for
$p=1,\dots,m$, $q=1,\dots,n$.

{\rm (b)} $\rank(x_{[1,\dots,p;q,\dots,n]})= \rank\bigl( (\wm
w_{12}\tr \wn)_{[1,\dots,p;q,\dots,n]} \bigr)$ for
$p=1,\dots,m$, $q=1,\dots,\allowmathbreak n$.

{\rm (c)} $\rank(x_{[1,\dots,m;p,\dots,q]})= q+1-p -
\rank\bigl( (\wn w_{11})_{[p,\dots,n;p,\dots,q]} \bigr)$ for
$2\le p\le q\le n$.

{\rm (d)} $\rank(x_{[p,\dots,q;1,\dots,n]})= q+1-p -
\rank\bigl( (w_{22}\wm)_{[p,\dots,q;1,\dots,q]} \bigr)$ for
$1\le p\le q\le m-1$.

\noindent Furthermore, $x\in \ol{\P}_w$ if and only if conditions
{\rm (a)--(d)} hold with each rank equality replaced by $\le$.
\endproclaim

\demo{Proof} We shall repeatedly use the following easy
observation: whenever a partial permutation matrix $u$ is
partitioned into blocks, the rank of $u$ equals the sum of the
ranks of the blocks.

Set $\xbar= \leftmatrix \wn &0\\ x &\wm
\rightmatrix$. In view of Theorem 3.9 and Proposition 4.1(a), we
have
$x\in\P_w$ if and only if
$$
\rank( \xbar_{[r,\dots,N;1,\dots,s]} )= \rank(
w_{[r,\dots,N;1,\dots,s]} ) \tag 4.2
$$
for all $r,s= 1,\dots,N$. Observe that (4.2) holds automatically
if $r=1$ (in which case both sides equal $s$), or if $s=N$ (in
which case both sides equal $N+1-r$). We shall consider (4.2) in a
number of separate cases.

{\bf Case 1}: $s\le n< r$. Set $p=r-n$ and $q=s$, and note
that
$$\xalignat2
\xbar_{[r,\dots,N;1,\dots,s]} &= x_{[p,\dots,m;1,\dots,q]}
&w_{[r,\dots,N;1,\dots,s]} &= (w_{21})_{[p,\dots,m;1,\dots,q]}.
\endxalignat
$$
Hence, (4.2) holds for $s\le n< r$ if and only if (a) holds.

{\bf Case 2}: $r,s\le n$ and $r+s\le n+1$. Since $r\le n+1-s$, we
have
$$
\xbar_{[r,\dots,N;1,\dots,s]} = \leftmatrix 0\\ w^s_\ci\\
x_{[1,\dots,m;1,\dots,s]} \rightmatrix
$$
(where the $0$ block is present only if $r< n+1-s$). It follows
that $\xbar_{[r,\dots,N;1,\dots,s]}$ has rank $s$ in this case.
Since $w\in
S_{N}^{\ge\wnwmpair}$, Lemma 3.12 says that $w(j)\ge n+1-j$ for
all $j$. For $j\le s$, we obtain $w(j)\ge n+1-s\ge r$, and so
$w_{[r,\dots,N;1,\dots,s]}$ has a $1$ in each of its $s$ columns.
Hence, $w_{[r,\dots,N;1,\dots,s]}$ has rank $s$, and thus (4.2)
always holds in the present case, independent of $x$.

{\bf Case 3}: $r,s\le n$ and $r+s>n+1$. Set $p=n+2-r$ and $q=s$,
so that $2\le p\le q$. We have
$$
\xbar_{[r,\dots,N;1,\dots,s]} = \leftmatrix
 w^{p-1}_\ci &0\\
x_{[1,\dots,m;1,\dots,p-1]} &x_{[1,\dots,m;p,\dots,q]}
\rightmatrix,
$$
and so $\rank( \xbar_{[r,\dots,N;1,\dots,s]} ) = p-1+ \rank(
x_{[1,\dots,m;p,\dots,q]} )$. For $j\le p-1$, we have $w(j)\ge
n+1-j\ge n+2-p= r$, which implies that
$w_{[r,\dots,N;1,\dots,p-1]}$ has rank $p-1$. Hence,
$$\align
\rank( w_{[r,\dots,N;1,\dots,s]} ) &= p-1+ \rank(
w_{[r,\dots,N;p,\dots,q]} ) \\
 &= p-1+q+1-p - \rank( w_{[1,\dots,r-1;p,\dots,q]} ) \\
 &= q - \rank\bigl( (w_{11})_{[1,\dots,n+1-p;p,\dots,q]} \bigr) \\
 &= q - \rank\bigl( (\wn w_{11})_{[p,\dots,n;p,\dots,q]} \bigr).
\endalign
$$
Therefore, (4.2) holds for $r,s\le n$ and $r+s>n+1$ if and only if
(c) holds.

{\bf Case 4}: $r,s>n$ and $r+s>m+2n$. Set $t=N+1-r$. Then $s\ge
n+t$, and so
$$
\xbar_{[r,\dots,N;1,\dots,s]} = \leftmatrix
x_{[r-n,\dots,m;1,\dots,n]} &w^t_\ci &0 \rightmatrix
$$
(where the $0$ block is present only if $s>n+t$). Hence,
$\xbar_{[r,\dots,N;1,\dots,s]}$ has rank $t$ in this case. Lemma
3.12 says that $w(j)\le m+2n+1-j$ for all $j$, and so for $j\ge
s+1$, we get $w(j)\le m+2n-s\le r-1$. Consequently, the nonzero
entries in rows $r,\dots,N$ of $w$ must occur in columns
$1,\dots,s$, from which we obtain $\rank(
w_{[r,\dots,N;1,\dots,s]} ) = N+1-r= t$. Therefore (4.2) always
holds in the present case.

{\bf Case 5}: $r,s>n$ and $r+s\le m+2n$. Set $p=r-n$ and $q=
N-s$, so that $1\le p\le q\le m-1$. Now
$$
\xbar_{[r,\dots,N;1,\dots,s]} = \leftmatrix
x_{[p,\dots,q;1,\dots,n]} &0 \\ x_{[q+1,\dots,m;1,\dots,n]}
&w^{s-n}_\ci \rightmatrix,
$$
and so $\rank( \xbar_{[r,\dots,N;1,\dots,s]} ) = s-n + \rank(
x_{[p,\dots,q;1,\dots,n]} )$. As in Case 4, for $j\ge s+1$, we have
$w(j)\le m+2n-s= n+q$, whence $w_{[n+q+1,\dots,N;1,\dots,s]}$ has
rank $(N+1)-(n+q+1)= m-q= s-n$. Hence,
$$\align
\rank( w_{[r,\dots,N;1,\dots,s]} ) &= s-n + \rank(
w_{[r,\dots,n+q;1,\dots,s]} ) \\
 &= s-n+ n+q+1-r - \rank( w_{[r,\dots,n+q;s+1,\dots,N]} ) \\
 &= s-n+ q+1-p - \rank\bigl( (w_{22})_{[p,\dots,q;m+1-q,\dots,m]}
\bigr) \\
 &= s-n+ q+1-p - \rank\bigl( (w_{22}\wm)_{[p,\dots,q;1,\dots,q]}
\bigr)
\endalign
$$
Therefore, (4.2) holds for $r,s>n$ and $r+s\le m+2n$ if and only
if (d) holds.

{\bf Case 6}: $2\le r\le n+1$ and $n\le s< N$. Set $p= N-s$ and
$q= n+2-r$, so that $1\le p\le m$ and $1\le q\le n$. We have
$$
\xbar_{[r,\dots,N;1,\dots,s]} = \leftmatrix w^{q-1}_\ci &0 &0 \\
x_{[1,\dots,p;1,\dots,q-1]} &x_{[1,\dots,p;q,\dots,n]} &0 \\
x_{[p+1,\dots,m;1,\dots,q-1]} &x_{[p+1,\dots,m;q,\dots,n]}
&w^{m-p}_\ci \rightmatrix
$$
(where the left column, respectively bottom row, is present only
if $q>1$, respectively $p<m$), and so
$\xbar_{[r,\dots,N;1,\dots,s]}$ has rank $q-1+m-p +
\rank( x_{[1,\dots,p;q,\dots,n]} )$. On the other hand,
$$\align
\rank( w_{[r,\dots,N;1,\dots,s]} ) &= s - \rank(
w_{[1,\dots,r-1;1,\dots,s]} ) \\
 &= s-(r-1) + \rank( w_{[1,\dots,r-1;s+1,\dots,N]} ) \\
 &= s+1-r + \rank\bigl( (w_{12})_{[1,\dots,n+1-q;m+1-p,\dots,m]} \\
 &= q-1+m-p + \rank\bigl( (\wm
w_{12}\tr \wn)_{[1,\dots,p;q,\dots,n]} \bigr).
\endalign
$$
Thus, (4.2) holds for $2\le r\le n+1$ and $n\le s< N$ if and
only if (b) holds.

Therefore (4.2) holds for $r,s=1,\dots,N$ if and
only if (a), (b), (c), (d) all hold, and we have established the
desired characterization of $\P_w$. The characterization of
$\ol{\P}_w$ follows from the information in Cases 1--6 together
with Theorem 3.13 and Proposition 4.1(b).
\qed\enddemo

\definition{4.3} Let $x\in \Mmn$ and $w= \leftmatrix w_{11}
&w_{12}\\ w_{21} &w_{22} \rightmatrix \in S_{N}^{\ge\wnwmpair}$
as in Theorem 4.2. According to Proposition 4.1(a)(c), the first
two conditions of the theorem are equivalent to the conditions
$x\in B_m^+ w_{21} B_n^+$ and $x\in B_m^-
\wm w_{12}\tr\wn B_n^-$, respectively. The corollary below follows
immediately. As we shall see in Example 4.5, the inclusion (4.3)
is typically proper.
\enddefinition

\proclaim{Corollary} Let $w\in
S_{N}^{\ge\wnwmpair}$, and write $w=
\leftmatrix w_{11} &w_{12}\\ w_{21} &w_{22} \rightmatrix$ as in
{\rm (4.1)}. Then 
$$
\P_w \subseteq B_m^+ w_{21} B_n^+ \cap B_m^-
\wm w_{12}\tr\wn B_n^-. \qquad\square \tag 4.3
$$
\endproclaim

\definition{4.4\. Example} Let $m=n$ and $u,v\in S_n$, and set $w=
\leftmatrix 0 &u\\ v &0 \rightmatrix \in S_N$. Via Lemma 3.12, it
is easily checked that $w\in S_{2n}^{\ge(\wn,\wn)}$. Let us use
Theorem 4.2 to compute $\P_w$ in this case.

Let $x\in M_n$. As discussed in \S4.3, conditions (a) and (b)
of the theorem require that $x\in B_n^+vB_n^+ \cap B_n^-
\wn u\tr\wn B_n^-$. In particular, $x$ must be invertible.
Conditions (c) and (d) of the theorem require
$$\xalignat2
\rank(x_{[1,\dots,n;p,\dots,q]}) &= q+1-p &&(2\le p\le q\le n) \\
\rank(x_{[p,\dots,q;1,\dots,n]} &= q+1-p &&(1\le p\le q< n).
\endxalignat
$$
These conditions hold automatically for $x\in\GLn$. Therefore, we
conclude that
$$
\P_{\leftmat 0 &u\\ v &0 \rightmat} = B_n^+vB_n^+ \cap B_n^-
\wn u\tr\wn B_n^-,
$$
a double Bruhat cell in $\GLn$. This recovers the previously known
description of $T$-orbits of symplectic leaves in $\GLn$ (cf\.
\cite{\HLthree, Appendix A} for the parallel case of $\SLn$).
\enddefinition

\definition{4.5\. Example} We give an example to show that
conditions (c) and (d) of Theorem 4.2 are typically not redundant,
i.e., (4.3) is typically a proper inclusion.

Take $m=n=3$, and consider the permutation matrix
$$
w = \leftmat 0&0&0&0&0&1\\ 0&1&0&0&0&0\\ 0&0&1&0&0&0\\
0&0&0&0&1&0\\ 0&0&0&1&0&0\\ 1&0&0&0&0&0 \rightmat \in M_6.
$$
Write $w= \leftmatrix w_{11} &w_{12}\\ w_{21} &w_{22}
\rightmatrix$ as in (4.1), and note that $w^3_\ci w_{12}\tr
w^3_\ci = w_{12}$. For $x\in M_3$, conditions (a) and (b) of
Theorem 4.2 require that
$$\align
\rank(x_{[p,\dots,3;1,\dots,q]}) &= \rank\bigl(
(w_{21})_{[p,\dots,3;1,\dots,q]} \bigr) = 1 \\
\rank(x_{[1,\dots,p;q,\dots,3]}) &= \rank\bigl(
(w_{12})_{[1,\dots,p;q,\dots,3]} \bigr) = 1
\endalign
$$
for $p,q=1,2,3$. These requirements boil down to $x_{31}, x_{13}
\ne 0$ and $\rank(x)=1$. It follows that $x_{11}, x_{33} \ne 0$.
Consequently,
$$
B_3^+ w_{21} B_3^+ \cap B_3^- w_{12} B_3^-= \biggl\{ x\in
\leftmatrix
\CCx &\CC &\CCx \\ \CC &\CC &\CC \\ \CCx &\CC &\CCx \rightmatrix
\biggm| \rank(x) = 1 \biggr\}.
$$

Next, observe that $w^3_\ci w_{11}= \leftmat 0&0&1\\ 0&1&0\\
0&0&0 \rightmat$ and $w_{22} w^3_\ci= \leftmat 0&1&0\\
0&0&1\\ 0&0&0 \rightmat$. Condition (c) of Theorem 4.2 requires
that
$$\align
\rank(x_{[1,2,3;2]}) &= 1 - \rank\bigl( (w^3_\ci w_{11})_{[2,3;2]}
\bigr) = 0 \\
\rank(x_{[1,2,3;2,3]}) &= 2 - \rank\bigl( (w^3_\ci
w_{11})_{[2,3;2,3]} \bigr) = 1 \\
\rank(x_{[1,2,3;3]}) &= 1 - \rank\bigl( (w^3_\ci w_{11})_{[3;3]}
\bigr) = 1.
\endalign
$$
The first equation means that the middle column of $x$ must be
zero; the other equations follow from the previous conditions.
Finally, condition (d) of Theorem 4.2 requires that
$$\align
\rank(x_{[1;1,2,3]}) &= 1 - \rank\bigl( (w_{22} w^3_\ci)_{[1;1]}
\bigr) = 1 \\
\rank(x_{[1,2;1,2,3]}) &= 2 - \rank\bigl( (w_{22}
w^3_\ci)_{[1,2;1,2]} \bigr) = 1 \\
\rank(x_{[2;1,2,3]}) &= 1 - \rank\bigl( (w_{22} w^3_\ci)_{[2;1,2]}
\bigr) = 1.
\endalign
$$
The last equation means that the middle row of $x$ must be
nonzero, while the other equations follow from the previous
conditions. 

We conclude that
$$
\P_w = \biggl\{ x\in \leftmatrix \CCx &0 &\CCx\\ \CCx &0 &\CCx\\
\CCx &0 &\CCx \rightmatrix \biggm| \rank(x) = 1 \biggl\},
$$
which is properly contained in $B_3^+ w_{21} B_3^+ \cap B_3^-
w_{12} B_3^-$. In fact, one can show that the latter intersection
is a disjoint union of four $T$-orbits of symplectic leaves,
corresponding to matrices of rank $1$ whose middle row or middle
column is zero or nonzero.
\enddefinition

\head 5. A second approach to $\Mmn$ 
by rank stratification \endhead

As above, fix positive integers $m$, $n$, and $N=m+n$. We investigate
the $T$-orbits of symplectic leaves of matrices with a given rank $t$,
which leads to a new description of orbits of leaves, 
quite different from Theorem 3.9.

\definition{5.1\. The set of rank $t$ matrices} Fix a nonnegative
integer $t\le \min\{m,n\}$, and set
$$
\Omnt = \{x\in\Mmn \mid \rank(x) =t\}.
\tag 5.1
$$
If $x\in \Omnt$, then $x\in \P_w$ for some $w\in
S_{N}^{\ge\wnwmpair}$ (Theorem 3.9), and Corollary 4.3 shows
that $\P_w \subseteq B_m^+w_{21}B_n^+$ for some partial permutation
matrix $w_{21}\in \Mmn$. Clearly $\rank(w_{21})= \rank(x)=t$,
whence $B_m^+w_{21}B_n^+ \subseteq \Omnt$, and so $x\in \P_w\subseteq
\Omnt$. Therefore, $\Omnt$ is a union of $T$-orbits of symplectic
leaves. Note that when $w\in S_{N}$ is written in the form
$\leftmat w_{11} &w_{12}\\ w_{21} &w_{22} \rightmat$ as in (4.1),
we have
$$
\rank(w_{21})= | w^{-1}(\{n+1,\dots,N\}) \cap \{1,\dots,n\} |.
$$
Hence, we define
$$
S_{N}^{\ge\wnwmpair}[t] = \bigl\{ w\in S_{N}^{\ge\wnwmpair}
\bigm| | w^{-1}(\{n+1,\dots,N\}) \cap \{1,\dots,n\} | =t \bigr\},
\tag 5.2
$$
so that we can state
$$
\Omnt= \bigsqcup \Sb w\in S_{N}^{\ge\wnwmpair}[t] \endSb \P_w \,.
\tag 5.3
$$
This statement invites us to view the matrix affine Poisson space
$\Mmn$ as stratified by matrix rank, and to analyze the $T$-orbits
$\P_w$ of symplectic leaves with special attention to their matrix
ranks. This analysis, carried out in the present section,
leads to new descriptions of the orbits $\P_w$.
\enddefinition

\definition{5.2\. $\Omnt$ as a Poisson homogeneous space} Under
the natural action of the group $G= \GLm\times\GLn$ on $\Mmn$,
given by $(a,b).x= axb^{-1}$, the set $\Omnt$ is the $G$-orbit of
the matrix
$$
I^{m,n}_t = \leftmatrix I_t &0_{t,n-t}\\ 0_{m-t,t}
&0_{m-t,n-t} \rightmatrix.
\tag 5.4
$$
Thus, $\Omnt$ is a homogeneous $G$-space. However, the action of
$G$ on $\Mmn$ is not a Poisson action for the standard Poisson
structure on $G$. To remedy this, we take
$$
G= \GLm\times\GLn^\bullet= (\GLm,\pi^m)\times (\GLn, -\pi^n),
$$
where $\pi^m$ and $\pi^n$ denote the standard Poisson structures
on $\GLm$ and $\GLn$ (recall \S1.4). With this change, the
action $G\times\Mmn \rightarrow \Mmn$ is a Poisson action, and
therefore $\Omnt$ is a Poisson homogeneous $G$-space.

Since the Poisson bivectorfield $\pi^{m,n}$ vanishes at
$I^{m,n}_t$, Theorem 1.8(a) shows that the Poisson homogeneous
$G$-space $\Omnt$ is isomorphic to $(G/Q^{m,n}_t,\,
\pi^{G/Q^{m,n}_t})$, where
$$\aligned
Q^{m,n}_t &= \stab_G(I^{m,n}_t) \\
 &= \biggl\{ \biggl( \leftmatrix a &b\\ 0 &d_1\rightmatrix,
\leftmatrix a &0\\ c &d_2\rightmatrix \biggr) \biggm| \matrix a\in
\GLt,\, b\in M_{t,m-t},\, c\in M_{n-t,t}, \\ d_1\in GL_{m-t},\,
d_2\in GL_{n-t} \endmatrix \biggr\}.
\endaligned \tag 5.5
$$
(Note that $\qfrak^{m,n}_t= \Lie(Q^{m,n}_t)$ can be
described in the same manner as (5.5).)
Thus, we can apply Theorem 1.10 to compute the $T$-orbits of
symplectic leaves within $\Omnt$. We sketch the steps in this
subsection, leaving details to the reader. When we compare the
results with those of Section 3 (see Theorem 5.11), we will obtain
an independent derivation, as a corollary of Theorem 3.9.

Write $\gfrak= \glm\times \gln^\bullet$ for the Lie bialgebra of
$G$. Because of the appearance
of $\gln^\bullet$ in the second factor of $\gfrak$, we use the
negative of the Killing form $\langle -,-\rangle$ on that factor.
Thus, the bilinear form to be used in $\gfrak$ is given by 
$$
\langle (x_1,x_2), (y_1,y_2)\rangle= \langle x_1,y_1\rangle
-\langle x_2,y_2\rangle,
$$
and the corresponding form on the double $D(\gfrak) \cong
\gfrak\oplus \gfrak$ (recall (1.6)) is given by
$$
\langle (x_1,x_2,x_3,x_4),(y_1,y_2,y_3,y_4)\rangle= 
\langle x_1,y_1\rangle -\langle x_2,y_2\rangle -\langle
x_3,y_3\rangle +\langle x_4,y_4\rangle.
$$
The duals appearing in the Manin triples $(D(G),\Delta(G),F)$ and
$(D(\gfrak), \Delta(\gfrak), \gfrak^*)$ (recall (1.4) and (1.5))
take the forms
$$
F= \{(a,b,a^{-1},b^{-1}) \mid a\in T_m,\, b\in T_n\}(N_m^+\times
N_n^+\times N_m^-\times N_n^-)
\tag 5.6
$$
and
$$
\gfrak^*= \{(x,y,-x,-y) \mid x\in \hfrak_m,\, y\in \hfrak_n \}+
(\nfrak_m^+\times\nfrak_n^+ \times \nfrak_m^-\times\nfrak_n^-),
\tag 5.7
$$
where we have written $N_l^\pm$ for the unipotent radical of
$B_l^\pm$ to avoid conflict with the notation (3.3).

In view of Theorem 1.8(b), the Drinfeld Langrangian subalgebra
corresponding to the base point $I^{m,n}_t$ in the present
situation has the form $\lfrak^{m,n}_t= \diag(\qfrak^{m,n}_t)
\oplus (\qfrak^{m,n}_t)^\perp$. As is easily computed,
$\lfrak^{m,n}_t$ consists of those 4-tuples
$$\multline
\biggl( \leftmatrix a_1&b_1\\ 0&d_1
\rightmatrix, \leftmatrix a_2&0\\ c_2&d_2 \rightmatrix, \leftmatrix
a_3&b_3\\ 0&d_3 \rightmatrix, \leftmatrix a_4&0\\ c_4&d_4
\rightmatrix \biggr) \\
 \in \leftmatrix \glt &M_{t,m-t}\\ 0 &\gl_{m-t}
\rightmatrix \times \leftmatrix \glt &0\\ M_{n-t,t} &\gl_{n-t}
\rightmatrix \times \leftmatrix \glt &M_{t,m-t}\\ 0 &\gl_{m-t}
\rightmatrix \times \leftmatrix \glt &0\\ M_{n-t,t} &\gl_{n-t}
\rightmatrix
\endmultline
$$
such that $a_1=a_2$, $a_3=a_4$, $d_1=d_3$, and $d_2=d_4$. Now
$\lfrak^{m,n}_t= \Lie(L^{m,n}_t)$ where the algebraic subgroup
$L^{m,n}_t \subseteq D(G)$ can be described in the same manner; we
write it as follows:
$$\multline
L^{m,n}_t= \biggl\{ \biggl( \leftmatrix a_1&b_1\\ 0&d_1
\rightmatrix, \leftmatrix a_1&0\\ c_1&d_2 \rightmatrix, \leftmatrix
a_2&b_2\\ 0&d_1 \rightmatrix, \leftmatrix a_2&0\\ c_2&d_2
\rightmatrix \biggr) \biggm| a_1,a_2\in \GLt, \\
b_1,b_2\in M_{t,m-t},\,
c_1,c_2\in M_{n-t,t},\, d_1\in GL_{m-t},\, d_2\in GL_{n-t}
\biggr\}.
\endmultline\tag 5.8
$$

We now apply Theorem 1.10, and conclude that the $T$-orbits of
symplectic leaves in $\Omnt$ are the irreducible components of the
sets
$$
\P^t_\sigma= \{ r_1I^{m,n}_tr_2^{-1} \mid (r_1,r_2,r_1,r_2)\in
(B_m^+\times B_n^+\times B_m^-\times B_n^-) \sigma L^{m,n}_t \},
\tag 5.9
$$
for $\sigma\in G\times G$. In fact, as we shall see later
(Corollary 5.12), each $\P^t_\sigma$ is a single $T$-orbit of
symplectic leaves. Thus, each $\P^t_\sigma$ is irreducible; we leave
it to the reader to seek a direct proof for this fact.

Next, an application of Theorem A.1 shows that a complete,
irredundant set of representatives for the $(B_m^+\times
B_n^+\times B_m^-\times B_n^-)$, $L^{m,n}_t$ double cosets in
$G\times G$ is given by
$$
S_m^{S^2_{m-t}} \times S_n^{S^1_t} \times
S_m^{S^1_t} \times S_n^{S^2_{n-t}}.
\tag 5.10
$$
Thus, we analyze the $\P^t_\sigma$ for $\sigma$ in the set (5.10).
In particular, we shall find a criterion for $\P^t_\sigma$ to be
nonempty (see Proposition 5.5).
\enddefinition

\proclaim{5.3\. Lemma} Let $\sigma= (y,v,z,u)\in S_m^{S^2_{m-t}}
\times S_n^{S^1_t} \times S_m^{S^1_t} \times S_n^{S^2_{n-t}}$.
Then $\P^t_\sigma$ consists of all matrices $r_1I^{m,n}_tr_2^{-1}$
for $r_1\in\GLm$ and $r_2\in\GLn$ such that
$$
\alignedat2 
r_1 &= b_1^+y= b_3^-z \leftmatrix a&b\\ 0&I_{m-t}
\rightmatrix  &&\qquad\quad (b_1^+\in B_m^+,\,
b_2^+\in B_n^+,\, b_3^-\in B_m^-,\, b_4^-\in B_n^-, \\
r_2 &= b_2^+v= b_4^-u \leftmatrix a&0\\
c&I_{n-t} \rightmatrix  &&\quad\qquad  a\in \GLt,\, b \in
M_{t,m-t},\ c\in M_{n-t,t}).
\endalignedat \tag 5.11
$$
\endproclaim

\demo{Proof} First, consider a matrix $x= r_1I^{m,n}_tr_2^{-1}$,
where $r_1\in\GLm$ and $r_2\in\GLn$ satisfy (5.11). Then
$$\align
(r_1,r_2,r_1,r_2) &= (b_1^+,b_2^+,b_3^-,b_4^-) (y,v,z,u) \biggl(
I_m,\, I_n,\, \leftmatrix a&b\\ 0&I_{m-t}
\rightmatrix, \leftmatrix a&0\\
c&I_{n-t} \rightmatrix \biggr)  \\
 &\in (B_m^+\times B_n^+\times B_m^-\times B_n^-) \sigma L^{m,n}_t,
\endalign
$$
whence $x\in \P^t_\sigma$.

Conversely, if $x\in \P^t_\sigma$, then $x= r_1I^{m,n}_tr_2^{-1}$
for some $r_1\in\GLm$ and $r_2\in\GLn$ such that 
$$
(r_1,r_2,r_1,r_2)= \biggl( b_1^+y \leftmatrix a_1&b_1\\
0&d_1 \rightmatrix,\, b_2^+v \leftmatrix a_1&0\\ c_1&d_2
\rightmatrix,\, b_3^-z \leftmatrix a_2&b_2\\ 0&d_1
\rightmatrix,\, b_4^-u \leftmatrix a_2&0\\ c_2&d_2
\rightmatrix \biggr),
$$
where $b_1^+\in B_m^+$, $b_2^+\in B_n^+$, $b_3^-\in B_m^-$,
$b_4^-\in B_n^-$, and the $a_i$, $b_i$, $c_i$, $d_i$ satisfy the
conditions of (5.8). Set $s_1= \leftmat a_1&b_1\\
0&d_1 \rightmat ^{-1}$ and $s_2= \leftmat a_1&0\\ c_1&d_2
\rightmat ^{-1}$, and observe that $(s_1,s_2,s_1,s_2) \in
L^{m,n}_t$. Hence, the 4-tuple $(r_1s_1,r_2s_2,r_1s_1,r_2s_2)$ lies
in $(B_m^+\times B_n^+\times B_m^-\times B_n^-) \sigma L^{m,n}_t$.
Since $s_1 I^{m,n}_t s_2^{-1}= I^{m,n}_t$, we have
$x= (r_1s_1) I^{m,n}_t (r_2s_2)^{-1}$, and so we may replace
$(r_1,r_2,r_1,r_2)$ by $(r_1s_1,r_2s_2,r_1s_1,r_2s_2)$. Thus,
there is no loss of generality in assuming that 
$$
(r_1,r_2,r_1,r_2)= \biggl( b_1^+y,\, b_2^+v,\, b_3^-z \leftmatrix
a &b\\ 0 &I_{m-t} \rightmatrix,\, b_4^-u \leftmatrix a &0\\
c &I_{n-t} \rightmatrix \biggr)
$$
for some $a\in\GLt$, $b \in
M_{t,m-t}$, $c\in M_{n-t,t}$. Now $r_1$ and $r_2$ satisfy (5.11),
and the proof is complete.
\qed\enddemo

\definition{5.4} Recall that the sets $S_n^{S^1_t}$ and
$S_n^{S^2_{n-t}}$ of minimal length coset representatives for the
subgroups $S^1_t$ and $S^2_{n-t}$ of $S_n$ can be described as
follows:
$$\align 
S_n^{S^1_t} &= \{u\in S_n\mid u(1)<\cdots<u(t)\}\\
S_n^{S^2_{n-t}} &= \{v\in S_n\mid v(t+1)<\cdots<v(n)\}.
\endalign
$$
\enddefinition

\proclaim{Lemma} {\rm (a)} If $v\in S_n^{S^1_t}$, then $v
\leftmat B_t^\pm &0\\ 0&I_{n-t} \rightmat \subseteq B_n^\pm v$.

{\rm (b)} If $u\in S_n^{S^2_{n-t}}$, then $u
\leftmat I_t &0\\ 0&B^\pm_{n-t} \rightmat \subseteq B_n^\pm u$.
\endproclaim

\demo{Proof} 
The lemma follows at once from the fact that 
for given a Weyl group $W$ and a subgroup $W_I$ generated 
by simple reflections for a subset $I$ of simple roots,
an element $w \in W$ belongs to the set $W^{W_I}$ of minimal
length  representatives of the cosets in $W/W_I$ 
if and only if $w(\alpha)$ is a positive root for 
any $\alpha \in I$, cf. \cite{\Car, Proposition 2.3.3}. 
\qed\enddemo

%
%

\proclaim{5.5\. Proposition} Let $\sigma= (y,v,z,u) \in
S_m^{S^2_{m-t}}\times S_n^{S^1_t}\times S_m^{S^1_t}\times
S_n^{S^2_{n-t}}$. Then
$$
\P^t_\sigma= \bigcup \Sb \tau\in S^1_t\\ z\tau\le y,\
v\tau^{-1}\le u \endSb
\bigl( B_m^+yB_m^+\cap B_m^-z\tau
\bigr).I^{m,n}_t.\bigl( \tau^{-1}B_n^-u^{-1}B_n^-\cap v^{-1}B_n^+
\bigr).
\tag 5.12
$$
Further, $\P^t_\sigma \ne\varnothing$ if and only if $z\le y$ and
$v\le u$. \endproclaim

\demo{Proof} By Theorem 2.4 and Proposition 2.2, $B^+yB^+\cap
B^-z\tau$ is nonempty if and only if $z\tau \le y$, and similarly
$B^-u^{-1}B^-\cap \tau v^{-1}B^+$ is nonempty if and only if
$v\tau^{-1}\le u$. Hence, the union in
(5.12) can just as well be taken over all $\tau\in S^1_t$.

Now assume for the moment that (5.12) has been proved. If 
$z\le y$ and $v\le u$, then the intersections $B^+yB^+\cap B^-z$ and
$B^-u^{-1}B^-\cap v^{-1}B^+$ are both nonempty, and (5.12)
yields $\P^t_\sigma \ne\varnothing$.
Conversely, if $\P^t_\sigma \ne\varnothing$, then because of
(5.12), there is some $\tau\in S^1_t$ such that both $B^+yB^+\cap
B^-z\tau$ and $\tau^{-1}B^-u^{-1}B^-\cap v^{-1}B^+$ are nonempty,
whence $z\tau\le y$ and $v\tau^{-1}\le u$. But since $z\in
S_m^{S^1_t}$ and $v\in S_n^{S^1_t}$, we see that
$z\le z\tau$ and $v\le v\tau^{-1}$. Therefore $z\le y$ and $v\le u$,
and the final statement of the theorem is proved.

It remains to prove (5.12).

If $x\in \P^t_\sigma$, then $x= r_1I^{m,n}_tr_2^{-1}$
for some $r_1\in\GLm$ and $r_2\in\GLn$ satisfying (5.11). 
By the $B_t^-$, $B_t^+$ Bruhat decomposition in
$\GLt$, we have $a= a^-\tau (a^+)^{-1}$ for some $a^\pm\in
B_t^\pm$ and $\tau\in S_t$. Set $s_1= r_1\leftmat
a^+ &-a^{-1}b\\ 0 &I_{m-t} \rightmat$ and $s_2= r_2\leftmat
a^+ &0\\ 0 &I_{n-t} \rightmat$, so that $x=
s_1I^{m,n}_t s_2^{-1}$ and
$$\align
s_1 &= b_1^+y\leftmatrix
a^+ &-a^{-1}b\\ 0 &I_{m-t} \rightmatrix = b_3^-z \leftmatrix a^-\tau&0\\
0&I_{m-t} \rightmatrix  \\
s_2 &= b_2^+v\leftmatrix
a^+ &0\\ 0 &I_{n-t} \rightmatrix = b_4^-u \leftmatrix
a^-\tau&0\\ ca^+ &I_{n-t} \rightmatrix.
\endalign
$$
It follows that $s_1\in B_m^+yB_m^+$ and $s_2\in
B_n^-uB_n^-\tau$, where we now view $\tau\in S^1_t\subseteq S_n$.
Since $z\in S_m^{S^1_t}$, Lemma 5.4 implies that $z\leftmat a^-
&0\\ 0 &I_{m-t}\rightmat \in B_m^-z$, whence $s_1\in
B_m^-z\tau$. Similarly, $v\in S_n^{S^1_t}$ implies that $v\leftmat
a^+ &0\\ 0 &I_{n-t}\rightmat \in B_n^+v$, whence $s_2\in
B_n^+v$. Thus,
$$
x= s_1I^{m,n}_ts_2^{-1} \in \bigl( B_m^+yB_m^+\cap B_m^-z\tau
\bigr).I^{m,n}_t.\bigl( \tau^{-1}B_n^-u^{-1}B_n^-\cap v^{-1}B_n^+
\bigr).
$$

Conversely, let $x\in\Mmn$ be a matrix such that
$$
x\in \bigl( B_m^+yB_m^+\cap B_m^-z\tau
\bigr).I^{m,n}_t.\bigl( \tau^{-1}B_n^-u^{-1}B_n^-\cap v^{-1}B_n^+
\bigr)
$$
for some $\tau\in S^1_t$. Then $x= r_1I^{m,n}_tr_2^{-1}$ where
$$
\xalignat2 r_1 &= b_1^+y \leftmatrix a_1&b_1\\ 0&d_1
\rightmatrix= b_3^-z\tau &r_2 &= b_4^-u \leftmatrix
a_2\tau&0\\ c_2\tau&d_2
\rightmatrix= b_2^+v \endxalignat
$$
where $b_1^+\in B_m^+$, $b_2^+\in B_n^+$, $b_3^-\in B_m^-$,
$b_4^-\in B_n^-$, while $a_1\in B_t^+$, $a_2\in B_t^-$, $b_1\in
M_{t,m-t}$,
$c_2\in M_{n-t,t}$, $d_1\in B^+_{m-t}$, $d_2\in B^-_{n-t}$. Since
$y\in S_m^{S^2_{m-t}}$, Lemma 5.4(b) implies that $y
\leftmat I_t&0\\ 0&d_1\rightmat \in B_m^+y$, and so $r_1=
\beta_1^+y \leftmat a_1&b_1\\ 0&I_{m-t}
\rightmat$ for some $\beta_1^+\in B_m^+$. Since $v\in
S_n^{S^1_t}$, Lemma 5.4(a) implies that $v \leftmat a_1&0\\
0&I_{n-t} \rightmat \in B_n^+v$, and so $r_2=
\beta_2^+v  \leftmat a_1&0\\ 0&I_{n-t}
\rightmat$ for some $\beta_2^+\in B_n^+$. Similarly, $r_1=
\beta_3^- z \leftmat a_2\tau &0\\ 0&I_{m-t} \rightmat$ and $r_2=
\beta_4^- u \leftmat a_2\tau&0\\ c'_2&I_{n-t}
\rightmat$ for some $\beta_3^-\in B_m^-$, $\beta_4^- \in B_n^-$,
and $c'_2\in M_{n-t,t}$. Consequently,
$$\align 
(r_1,r_2,r_1,r_2) &= (\beta_1^+y,\beta_2^+v,\beta_3^-z,\beta_4^-u)
\left( \leftmat a_1&b_1\\ 0&I_{m-t}
\rightmat,\, \leftmat a_1&0\\ 0&I_{n-t}
\rightmat,\, \leftmat a_2\tau &0\\ 0&I_{m-t} \rightmat,\, 
\leftmat a_2\tau&0\\ c'_2&I_{n-t}
\rightmat \right) \\
 &\in (B_m^+\times B_n^+\times B_m^-\times B_n^-)\sigma L_t,
\endalign
$$ 
and so $x\in \P^t_\sigma$. Therefore (5.12)
holds. \qed\enddemo

\definition{5.6} In view of Proposition 5.5, the following set
indexes the nonempty $\P^t_\sigma$:
$$
\Sigmnt = \bigl\{ (y,v,z,u) \in S_m^{S^2_{m-t}} \times S_n^{S^1_t}
\times S_m^{S^1_t} \times S_n^{S^2_{n-t}} \bigm| z\le y,\ v\le u
\bigr\}.
\tag 5.13
$$
In order to match the $\P^t_\sigma$ with appropriate $T$-orbits
$\P_w$ of symplectic leaves, we need a bijection between $\Sigmnt$
and the index set $S_{N}^{\ge\wnwmpair}[t]$ defined in (5.2).
Recall from (3.24) and Lemma 3.12 that
$$
S_{N}^{\ge\wnwmpair}= \wmn S_{N}^{[-n,m]}.
$$
Hence, we define
$$\aligned
\Smnt &= \wmn S_{N}^{\ge\wnwmpair}[t] \\
 &= \bigl\{ w\in S_{N}^{[-n,m]} \bigm| |w^{-1}\bigl(
\{1,\dots,m\} \bigr) \cap \{1,\dots,n\}| = t \bigr\},
\endaligned \tag 5.14
$$
so that 
$S_{N}^{\ge\wnwmpair}[t]= \wmn\Smnt$. It is convenient to
first construct a bijection $\Sigmnt \rightarrow \Smnt$. To
describe that, we will need the matrix $I^{m,n}_t \in \Mmn$ and the
analogous matrix $I^{n,m}_t \in M_{n,m}$, as well as
$$\xalignat2 
J^m_t &= \leftmatrix 0_t &0_{t,m-t}\\ 0_{m-t,t} &I_{m-t}
\rightmatrix \in M_m  &J^n_t &= \leftmatrix 0_t &0_{t,n-t}\\
0_{n-t,t} &I_{n-t} \rightmatrix \in M_n \, .
\tag 5.15 \endxalignat 
$$
\enddefinition

\proclaim{5.7\. Lemma} Let $u,v\in S_n$.

{\rm (a)} If $u\in S_n^{S^2_{n-t}}$ and $v\le u$, then $v(j)\ge
u(j)$ for $j=t+1,\dots,n$.

{\rm (b)} If $v\in S_n^{S^1_t}$ and $v(j)\ge u(j)$ for
$j=t+1,\dots,n$, then
$v\le u$.
\endproclaim

\demo{Proof} First consider subsets $U,V\subseteq \onen$ with $|U|=|V|$,
and let $\Util$ and $\Vtil$ denote their complements in $\onen$. We claim
that $V\le U$ if and only if $\Vtil \ge \Util$.

Assume first that $\Vtil\ge\Util$. Label the elements of the four sets in
ascending order:
$$\xalignat2 U &= \{u_1<\cdots<u_r\} &V &= \{v_1<\cdots<v_r\}\\
\Util &= \{\util_1<\cdots<\util_{n-r}\} &\Vtil &=
\{\vtil_1<\cdots<\vtil_{n-r}\}. \endxalignat$$
We have $\vtil_i\ge\util_i$ for all $i$, and must show that $v_j\le u_j$
for all $j$.

Consider the interval $L= \{1,2,\dots,v_j-1\}$ for some $j\le r$.
Since $L$ contains exactly $j-1$ elements of $V$, it contains the
first $v_j-j$ elements of $\Vtil$. So for $i=1,\dots,v_j-j$, we
have $\vtil_i\in L$ and $\util_i\le \vtil_i$, whence $\util_i\in
L$. Thus, $L$ contains at least $v_j-j$ elements of $\Util$, and
hence at most $j-1$ elements of $U$. It follows that $u_j\notin
L$, whence $u_j\ge v_j$. Therefore
$V\le U$, as desired.

The fact that $V\le U$ implies $\Vtil\ge\Util$ follows by reversing the
roles of these sets and their complements.

(a) By assumption, $v(\{1,\dots,t\}) \le u(\{1,\dots,t\})$,
and so the claim above implies that $v(\{t+1,\dots,n\}) \ge
u(\{t+1,\dots,n\})$. Since $u\in S_n^{S^2_{n-t}}$, the least
element of $u(\{t+1,\dots,n\})$ is $u(t+1)$, and consequently
$u(t+1)\le v(t+1)$. Moreover, $u\in S_n^{S^2_{n-r}}$ for $t\le
r<n$, and so the same argument yields $u(r+1)\le v(r+1)$ for $t\le
r<n$.

(b) Our assumption implies that $v(\{t+1,\dots,n\}) \ge
u(\{t+1,\dots,n\})$, and so the claim above yields
$v(\{1,\dots,t\}) \le u(\{1,\dots,t\})$. Since $v(1)<\cdots,v(t)$,
it follows that $v(\{1,\dots,r\}) \le u(\{1,\dots,r\})$ for
$r=1,\dots,t$. Moreover, for $r=t,\dots,n-1$, we have
$v(\{r+1,\dots,n\}) \ge u(\{r+1,\dots,n\})$ and the claim yields
$v(\{1,\dots,r\}) \le u(\{1,\dots,r\})$. Therefore $v\le u$.
\qed\enddemo

\definition{5.8\. Partial permutations} Just as with permutations (cf\.
\S3.1), we view any partial permutation matrix $w$ as both a
matrix and a function (a bijection from its domain to its range). Write
$\dom(w)$ and
$\rng(w)$ for the domain and range of $w$; then the matrix form of
$w$ has a $1$ in position $w(j),j$ for each $j\in\dom(w)$, and a
$0$ in all other positions. Observe that $w\tr$ is the inverse
bijection, from $\rng(w)$ to $\dom(w)$.
\enddefinition

\proclaim{5.9\. Proposition} There is a bijection $\phi: \Sigmnt
\rightarrow \Smnt$ given by
$$\spreadmatrixlines{1ex}
\phi(y,v,z,u) = \left[\matrix \wm yI^{m,n}_tv^{-1} &\wm
yJ^m_tz^{-1}\wm \\ uJ^n_tv^{-1} &uI^{n,m}_tz^{-1}\wm
\endmatrix\right] .
\tag 5.16
$$
\endproclaim

\demo{Proof} Let $(y,v,z,u)\in \Sigmnt$, and let
$$
w = \phi(y,v,z,u)= \leftmatrix w_{11} &w_{12} \\ w_{21} &w_{22}
\rightmatrix ,
$$
where the $w_{ij}$ stand for the blocks shown in (5.16). Since $w$
can be expressed in the form
$$
w= \leftmatrix \wm y &0\\ 0 &u \rightmatrix \leftmatrix I^{m,n}_t
&J^m_t\\ J^n_t &I^{n,m}_t \rightmatrix \leftmatrix v^{-1} &0\\ 0
&z^{-1}\wm \rightmatrix ,
$$
it is clear that $w$ is a permutation matrix, which we identify
with a permutation in $S_{N}$ in the usual way. Observe that
$$
| w^{-1}\bigl( \{1,\dots,m\} \bigr) \cap \{1,\dots,n\} | =
\rank(w_{11}) = t.
$$

By Lemma 5.7(a), $z(j)\ge y(j)$ and $v(j)\ge u(j)$ for $j>t$. Thus,
$w_{21}v(j)= u(j)\le v(j)$ for $j>t$, and so $w_{21}(i)\le i$ for
all $i\in \dom(w_{21})$. It follows that $w(i)\le i+m$ for all
$i$. Similarly, $w_{12}\wm z(j)= \wm y(j) \ge \wm z(j)$ for all
$j>t$ and so $w_{12}(i)\ge i$ for all $i\in \dom(w_{12})$, whence
$w(i)\ge i-n$ for all $i$. Therefore $w\in \Smnt$, which shows
that the rule (5.16) does define a map $\phi$ from $\Sigmnt$ to
$\Smnt$.

Observe that $y(j)= \wm w_{11}v(j)$ for $j\le t$. Since $v(1)<
\cdots< v(t)$ (because $v\in S_n^{S^1_t}$), it follows that the
restriction of $y$ to $\{1,\dots,t\}$ is determined by $w_{11}$.
But $y\in S_m^{S^2_{m-t}}$, and thus $y$ is completely determined
by $w_{11}$. Similarly, $u(j)= w_{22}\wm z(j)$ for $j\le t$ and
$z(1)< \cdots< z(t)$, whence the restriction of $u$ to
$\{1,\dots,t\}$ is determined by $w_{22}$. Since $u\in
S_n^{S^2_{n-t}}$, it follows that $u$ is completely determined by
$w_{22}$.

For $j=t+1,\dots,n$, we have $u(j)= w_{21}v(j)$ and so $v(j)=
w_{21}\tr u(j)$. Since $v\in S_n^{S^1_t}$, it follows that $v$ is
completely determined by $u$ and $w_{21}$. Similarly, for
$j=t+1,\dots,m$, we have $\wm y(j)= w_{12}\wm z(j)$ and so $z(j)=
\wm w_{12}\tr \wm y(j)$. Since $z\in S_m^{S^1_t}$, it follows that
$z$ is completely determined by $y$ and $w_{12}$. Therefore,
$(y,v,z,u)$ is completely determined by $w$, which shows that the
map $\phi$ is injective.

Now consider an arbitrary element $w\in \Smnt$, and write
$$\xalignat2 
w &= \left[\matrix w_{11} &w_{12}\\ w_{21} &w_{22}
\endmatrix\right],
 &&\pmatrix w_{11} \in \Mmn &w_{12}\in M_m \\
w_{21}\in M_n &w_{22}\in M_{n,m} \endpmatrix.
\endxalignat
$$
Each $w_{ij}$ is a partial permutation matrix, and
$$\aligned
\dom(w_{11})\sqcup \dom(w_{21}) &= \rng(w_{21})\sqcup
\rng(w_{22})= \{1,\dots,n\} \\
\dom(w_{12})\sqcup \dom(w_{22}) &= \rng(w_{11})\sqcup
\rng(w_{12})= \{1,\dots,m\}.
\endaligned  \tag 5.17
$$
Further, $\rank(w_{11})=t$ (recall (5.14)), from which we see that
$\rank(w_{12})= m-t$ and $\rank(w_{21})= n-t$, and hence
$\rank(w_{22})=t$. Since
$i-n\le w(i)\le i+m$ for all $i=1,\dots,N$, we have
$w_{12}(j)\ge j$ for all
$j\in
\dom(w_{12})$ and $w_{21}(j)\le j$ for all $j\in \dom(w_{21})$.

Write the elements of $\dom(w_{11})$ in ascending order:
$\dom(w_{11})= \{v_1< \cdots< v_t\}$. Set $y(j)= \wm w_{11}(v_j)$
for $j=1,\dots,t$, and extend (uniquely) to a permutation
$y\in S_m^{S^2_{m-t}}$. Write the elements of $\dom(w_{22})$ in
descending order: $\dom(w_{22})= \{z_1> \cdots> z_t\}$. Set $u(j)=
w_{22}(z_j)$ for $j=1,\dots,t$, and extend (uniquely) to a
permutation $u\in S_n^{S^2_{n-t}}$. 
Next, observe using (5.17) that
$$\align
u(\{t+1,\dots,n\}) &= \{1,\dots,n\}\setminus \rng(w_{22})=
\rng(w_{21})= \dom(w_{21}\tr)  \\
\rng(w_{21}\tr) &= \dom(w_{21})=
\{1,\dots,n\}\setminus \dom(w_{11}).
\endalign
$$
Hence, we can define a permutation $v\in S_n^{S^1_t}$ such that
$v(j)= v_j$ for $j=1,\dots,t$ and $v(j)= w_{21}\tr u(j)$ for
$j=t+1,\dots,n$. Similarly,
$$\align
\wm y(\{t+1,\dots,m\}) &= \{1,\dots,m\}\setminus \rng(w_{11})=
\rng(w_{12})= \dom(w_{12}\tr)  \\
\wm(\rng(w_{12}\tr)) &= \wm(\dom(w_{12}))= \{1,\dots,m\}\setminus
\wm(\dom(w_{22})),
\endalign
$$
and so we can define a permutation $z\in S_m^{S^1_t}$ such that
$z(j)= \wm(z_j)$ for $j=1,\dots,t$ and $z(j)= \wm w_{12}\tr \wm
y(j)$ for $j=t+1,\dots,m$.

We have now defined $(y,v,z,u)\in S_m^{S^2_{m-t}} \times
S_n^{S^1_t} \times S_m^{S^1_t} \times S_n^{S^2_{n-t}}$. For
$j=t+1,\dots,m$, we have
$$
y(j)= \wm w_{12}\wm z(j) \le \wm \wm z(j)= z(j),
$$
and so $z\le y$ by Lemma 5.7(b). Similarly, $u(j)= w_{21}v(j)\le
v(j)$ for $j=t+1,\dots,n$, and so $v\le u$. Thus, $(y,v,z,u) \in
\Sigmnt$. Finally, we analyze the domains and actions of the four
components of $\phi(y,v,z,u)$, as follows.
$$\xalignat2
\wm yI^{m,n}_tv^{-1}\,: &&\text{domain} &= v(\{1,\dots,t\})=
\{v_1,\dots,v_t\}= \dom(w_{11}) \\
 &&v_j &= v(j) \mapsto \wm y(j)= w_{11}(v_j) \\
\wm yJ^m_tz^{-1}\wm\,: &&\text{domain} &= \wm z(\{t+1,\dots,m\})=
\rng(w_{12}\tr)= \dom(w_{12}) \\
 &&\wm z(j) &\mapsto \wm y(j)= w_{12}\wm z(j) \\
uJ^n_tv^{-1}\,: &&\text{domain} &= v(\{t+1,\dots,n\})=
\rng(w_{21}\tr)= \dom(w_{21}) \\
 &&v(j) &\mapsto u(j)= w_{21}v(j) \\
uI^{n,m}_tz^{-1}\wm\,: &&\text{domain} &= \wm z(\{1,\dots,t\})=
\{z_1,\dots,z_t\}= \dom(w_{22}) \\
 &&z_j &= \wm z(j) \mapsto u(j)= w_{22}(z_j).
\endxalignat
$$
This shows that $\phi(y,v,z,u)= \leftmat w_{11} &w_{12}\\
w_{21} &w_{22} \rightmat = w$, and therefore that $\phi$ is
surjective.
\qed\enddemo

\proclaim{5.10\. Corollary} There is a bijection $\Sigmnt
\rightarrow S_{N}^{\ge\wnwmpair}[t]$ given by
$$\spreadmatrixlines{1ex}
(y,v,z,u) \longmapsto w^{N}_\ci \phi(y,v,z,u)= \left[\matrix
\wn uJ^n_tv^{-1} &\wn uI^{n,m}_tz^{-1}\wm\\ yI^{m,n}_tv^{-1}
&yJ^m_tz^{-1}\wm \endmatrix\right] . \qquad\square
\tag 5.18
$$
\endproclaim

We are now ready to state and prove the main theorem of the
section. The description it provides of orbits $\P_w$ of
symplectic leaves requires a union involving more than one set in
general (see Example 5.14). For a class of cases in which only a
single term is required, see Theorem 6.1.

\proclaim{5.11\. Theorem} Let $w\in S_{N}^{\ge\wnwmpair}[t]$
{\rm(}recall {\rm(5.2))}. Then $w= \wmn\phi(\sigma)$ for a unique
$4$-tuple $\sigma= (y,v,z,u)\in \Sigmnt$ {\rm(}recall {\rm(5.13))},
and
$$
\P_w= \P^t_\sigma = \bigcup \Sb \tau\in S^1_t\\
z\tau\le y,\ v\tau^{-1}\le u \endSb
\bigl( B_m^+yB_m^+\cap B_m^-z\tau
\bigr).I^{m,n}_t.\bigl( \tau^{-1}B_n^-u^{-1}B_n^-\cap v^{-1}B_n^+
\bigr).
\tag 5.19
$$
\endproclaim

\demo{Proof} The existence and uniqueness of $\sigma$ are given by
Corollary 5.10, and the second equality in (5.19) by Proposition
5.5. It remains to prove that $\P_w= \P^t_\sigma$, for which we
shall use the description of $\P^t_\sigma$ given in Lemma 5.3.

Observe that, in block form, $w= \leftmat \wn &0\\ 0 &I_m
\rightmat s \leftmat I_n &0\\ 0 &\wm \rightmat$,
where 
$$\spreadmatrixlines{1ex}
s= \left[\matrix uJ^n_tv^{-1} &uI^{n,m}_tz^{-1} \\
yI^{m,n}_tv^{-1} &yJ^m_tz^{-1} \endmatrix\right] .
$$
Hence (recall (3.18)), $\P_w$ consists of those matrices $x\in
\Mmn$ such that
$$\aligned
\leftmatrix I_n &0\\ x &I_m \rightmatrix &\in \leftmatrix \wn &0\\
0&I_m \rightmatrix B^+ \leftmatrix \wn &0\\ 0 &I_m
\rightmatrix s \leftmatrix I_n &0\\ 0 &\wm
\rightmatrix B^+ \leftmatrix I_n &0\\ 0 &\wm
\rightmatrix  \\
 &= \leftmatrix B^-_n &M_{n,m}\\ 0 &B^+_m \rightmatrix s
\leftmatrix B^+_n &M_{n,m}\\ 0 &B^-_m \rightmatrix. 
\endaligned \tag 5.20
$$

If $x\in\Mmn$ satisfies (5.20), then
$$\gather
 \leftmatrix I_n &0\\ x &I_m\rightmatrix = \leftmatrix
\alpha_1&\beta_1\\ 0&\gamma_1 \rightmatrix \leftmatrix uJ^n_t
v^{-1} &uI^{n,m}_t z^{-1}\\ yI^{m,n}_t v^{-1} &yJ^m_t z^{-1}
\rightmatrix \leftmatrix \alpha_2&\beta_2\\ 0&\gamma_2
\rightmatrix \tag  5.21 \\
(\alpha_1\in B_n^-,\, \alpha_2\in B_n^+,\, \gamma_1\in B_m^+,\,
\gamma_2\in B_m^-,\, \beta_1,\beta_2\in M_{n,m}).
\endgather
$$
Set $b= u^{-1}\alpha_1^{-1}\beta_1y\in M_{n,m}$, and
rewrite (5.21) in the form
$$\aligned 
I_n &= \alpha_1u(J^n_t +bI^{m,n}_t )v^{-1}\alpha_2 \\
0 &= \alpha_1u(J^n_t v^{-1}\beta_2+I^{n,m}_t z^{-1}\gamma_2)+
\alpha_1ub(I^{m,n}_t v^{-1}\beta_2+J^m_t z^{-1}\gamma_2) \\
x &= \gamma_1yI^{m,n}_t v^{-1}\alpha_2 \\
I_m &= \gamma_1y(I^{m,n}_t v^{-1}\beta_2+J^m_t z^{-1}\gamma_2). 
\endaligned \tag 5.22
$$
Multiply the first equation of (5.22) on the right by
$\alpha_2^{-1}vJ^n_t $ and by $\alpha_2^{-1}vI^{n,n}_t $, and the
fourth on the left by $I^{m,m}_t y^{-1}\gamma_1^{-1}$ and by $J^m_t
y^{-1}\gamma_1^{-1}$, to obtain
$$\alignedat2 
\alpha_2^{-1}vJ^n_t  &= \alpha_1uJ^n_t 
&\qquad\qquad\qquad \alpha_2^{-1}vI^{n,n}_t  &=
\alpha_1ubI^{m,n}_t  \\ 
I^{m,m}_t y^{-1}\gamma_1^{-1} &= I^{m,n}_t v^{-1}\beta_2
&\qquad\qquad\qquad J^m_t y^{-1}\gamma_1^{-1} &= J^m_t
z^{-1}\gamma_2. 
\endalignedat \tag 5.23
$$
Adding the two equations in each row of (5.23) yields
$$\xalignat2
\alpha_2^{-1}v &= \alpha_1u(J^n_t +bI^{m,n}_t )
&y^{-1}\gamma_1^{-1} &= I^{m,n}_t v^{-1}\beta_2+
J^m_t z^{-1}\gamma_2.
\tag 5.24 \endxalignat 
$$
Now substitute the second equation of (5.24) into the second
equation of (5.22), and multiply on the left by $I^{m,n}_t
u^{-1}\alpha_1^{-1}$, to obtain
$$
I^{m,m}_t z^{-1}\gamma_2+ I^{m,n}_t by^{-1}\gamma_1^{-1} =0.
\tag 5.25
$$ 
The last equation of (5.23) combines with (5.25) to yield
$z^{-1}\gamma_2= (J^m_t -I^{m,n}_t b)y^{-1}\gamma_1^{-1}$, and
consequently
$$
\gamma_1y= \gamma_2^{-1}z(J^m_t -I^{m,n}_t b). \tag 5.26
$$

Write $b= \leftmat b_{11} &b_{12}\\ b_{21} &b_{22} \rightmat \in
\leftmat M_t &M_{t,m-t}\\ M_{n-t,t} &M_{n-t,m-t} \rightmat$.
Since, as we see from (5.26), the matrix $J^m_t -I^{m,n}_t b=
\leftmat -b_{11} &-b_{12}\\ 0&I_{m-t} \rightmat$ is invertible,
$b_{11}\in\GLt$. Now set
$$\aligned 
r_1 &= \gamma_1y= \gamma_2^{-1}z(J^m_t -I^{m,n}_t b)=
\gamma_2^{-1}z \leftmatrix -b_{11} &-b_{12}\\ 0&I_{m-t}
\rightmatrix \\ 
r_2 &= \alpha_2^{-1}v= \alpha_1u(J^n_t +bI^{m,n}_t
)= \alpha_1u \leftmatrix b_{11} &0\\ b_{21} &I_{n-t} \rightmatrix.
\endaligned \tag 5.27
$$
Since $z \leftmat -I_t &0\\ 0&I_{m-t}
\rightmat \in zT_m= T_mz$, we have $r_1\in B_m^-z \leftmat
b_{11} &b_{12}\\ 0&I_{m-t} \rightmat$. Thus, $r_1$ and $r_2$
satisfy (5.11), and so 
$x = \gamma_1yI^{m,n}_t v^{-1}\alpha_2= r_1I^{m,n}_t r_2^{-1}\in \P^t_\sigma$.

Conversely, if $x\in \P^t_\sigma$, then, making use of the relation $z
\leftmat -I_t &0\\ 0&I_{m-t}
\rightmat \in T_mz$ as above, $x= r_1I^{m,n}_t r_2^{-1}$ where
$$\gather
r_1 = \gamma_1y= \gamma_2^{-1}z \leftmatrix -b_{11}
&-b_{12}\\ 0 &I_{m-t} \rightmatrix \qquad\qquad r_2 =
\alpha_2^{-1}v= \alpha_1u \leftmatrix b_{11} &0\\ b_{21} &I_{n-t}
\rightmatrix \tag 5.28  \\
(\gamma_1\in B_m^+,\, \gamma_2\in B_m^-,\, \alpha_2\in B_n^+,\,
\alpha_1\in B_n^-,\, b_{11}\in\GLt,\, b_{12}\in M_{t,m-t},\,
b_{21}\in M_{n-t,t}).
\endgather
$$
In particular,
$$
x= \gamma_1 y I^{m,n}_t v^{-1}\alpha_2 \, . \tag 5.29
$$
Set $b= \leftmat b_{11} &b_{12}\\ b_{21}
&0 \rightmat \in M_{n,m}$; then (5.28) can be rewritten as
$$\xalignat2 
r_1 &= \gamma_1y= \gamma_2^{-1}z (J^m_t -I^{m,n}_t b)
&r_2 &= \alpha_2^{-1}v= \alpha_1u (J^n_t +bI^{m,n}_t ).
\tag 5.30 \endxalignat 
$$ 
The first equation of (5.30) implies that
$$
z^{-1}\gamma_2= (J^m_t -I^{m,n}_t b)y^{-1}\gamma_1^{-1}. \tag
5.31
$$ 
The second equation of (5.30), together with (5.31), yields
$$\alignedat2 
\alpha_2^{-1}vJ^n_t  &= \alpha_1uJ^n_t 
&\qquad\qquad\qquad\quad \alpha_2^{-1}vI^{n,n}_t  &=
\alpha_1ubI^{m,n}_t  \\ 
J^m_t y^{-1}\gamma_1^{-1} &=
J^m_t z^{-1}\gamma_2. \endalignedat \tag 5.32
$$

Now set $\beta_1= \alpha_1uby^{-1}$ and $\beta_2= -\alpha_2\alpha_1
u (I^{n,m}_t +bJ^m_t )z^{-1}\gamma_2$ in $M_{n,m}$. From (5.32) and
the definitions of $\beta_1$ and $\beta_2$, we get
$$
(\alpha_1uJ^n_t +\beta_1yI^{m,n}_t )v^{-1}\alpha_2 =
(\alpha_2^{-1}vJ^n_t + \alpha_2^{-1}vI^{n,n}_t ) v^{-1}\alpha_2
=I_n, \tag 5.33
$$ 
which implies 
$$
\alpha_2^{-1}v = \alpha_1uJ^n_t +\beta_1yI^{m,n}_t =
\alpha_1u(J^n_t +bI^{m,n}_t ), \tag5.34
$$
 as well as
$$\multline 
(\alpha_1uJ^n_t +\beta_1yI^{m,n}_t )v^{-1}\beta_2
+ (\alpha_1uI^{n,m}_t +\beta_1yJ^m_t )z^{-1}\gamma_2 = \\
 \alpha_2^{-1}\beta_2+
\alpha_1u(I^{n,m}_t +bJ^m_t )z^{-1}\gamma_2 =0. \endmultline
\tag 5.35
$$  
Note that (5.34) implies that $v^{-1}\alpha_2\alpha_1u= (J^n_t
+bI^{m,n}_t )^{-1}$, and so, from (5.31) and the definition of
$\beta_2$, we have
$$\aligned 
I^{m,n}_t v^{-1}\beta_2 &= -I^{m,n}_t  (J^n_t +bI^{m,n}_t )^{-1}
(I^{n,m}_t +bJ^m_t )  (J^m_t -I^{m,n}_t b)y^{-1}\gamma_1^{-1} \\
 &= -I^{m,n}_t  \leftmatrix b_{11} &0\\ b_{21} &I_{n-t}
\rightmatrix^{-1} \leftmatrix I_t  &b_{12}\\ 0&0 \rightmatrix 
\leftmatrix -b_{11}
&-b_{12}\\ 0 &I_{m-t} \rightmatrix y^{-1}\gamma_1^{-1} \\
 &= - \leftmatrix b_{11}^{-1} &0\\ 0 &0_{m-t,n-t} \rightmatrix
\leftmatrix -b_{11} &0\\ 0&0_{n-t,m-t} \rightmatrix
y^{-1}\gamma_1^{-1} = I^{m,m}_t y^{-1}\gamma_1^{-1}. \endaligned
\tag 5.36
$$
Consequently, with the help of (5.32), we get
$$
\gamma_1y(I^{m,n}_t v^{-1}\beta_2+ J^m_t z^{-1}\gamma_2)= \gamma_1y
(I^{m,m}_t y^{-1}\gamma_1^{-1}+ J^m_t y^{-1}\gamma_1^{-1}) =I_m.
\tag 5.37
$$

Combine (5.33), (5.29), (5.35) and (5.37) to see that (5.21) holds,
whence (5.20), and therefore $x\in \P_w$.
\qed\enddemo

Theorem 5.11 verifies the main conclusions of \S5.2, as
follows.

\proclaim{5.12\. Corollary} The $T$-orbits of symplectic leaves
within $\Omnt$ are precisely the sets $\P^t_\sigma$ {\rm(}recall
{\rm (5.9))} for $\sigma\in \Sigmnt$ {\rm(}recall {\rm (5.13))}.
\endproclaim

\demo{Proof} Equation (5.3), Corollary 5.10, and Theorem 5.11.
\qed\enddemo

\definition{5.13\. Example} We recalculate Example 4.5 from the
viewpoint of Theorem 5.11. Here $m=n=3$ and $t=1$. Via the proof
of Proposition 5.9, one finds that the unique 4-tuple
$\sigma= (y,v,z,u) \in \Sigma^{3,3}_1$ such that
$w^6_\ci\phi(\sigma) = w$ is given by
$$
(y,v,z,u)= \left( \leftmat 0&1&0\\ 0&0&1\\ 1&0&0 \rightmat,
\leftmat 1&0&0\\ 0&0&1\\ 0&1&0 \rightmat, \leftmat 1&0&0\\ 0&1&0\\
0&0&1
\rightmat, \leftmat 0&1&0\\ 0&0&1\\ 1&0&0 \rightmat \right).
$$
Since $S^1_1$ consists only of the identity, Theorem 5.11 yields
$$
\P_w= \P^1_\sigma= \bigl( B_3^+yB_3^+\cap B_3^-z \bigr) .
I^{3,3}_1 . \bigl( B_3^-u^{-1}B_3^-\cap v^{-1}B_3^+ \bigr). \tag
5.38
$$

It follows from Proposition 4.1 that
$$
B_3^+yB_3^+= \{ x\in GL_3 \mid x_{31}\ne 0 \text{\ and\ }
\rank(x_{[2,3;1,2]}) =1 \},
$$
and consequently (since $z$ is the identity)
$$
B_3^+yB_3^+\cap B_3^-z= \left\{ x\in
\leftmatrix \CCx &0 &0\\ \CCx &\CCx &0\\ \CCx &\CCx &\CCx
\rightmatrix \biggm| 
\rank(x_{[2,3;1,2]}) =1 \right\}. \tag 5.39
$$
On the other hand,
$$
B_3^-u^{-1}B_3^-= \{ x\in GL_3 \mid x_{13}\ne 0 \text{\ and\ }
\rank(x_{[1,2;2,3]}) =1 \},
$$
and so
$$
B_3^-u^{-1}B_3^-\cap v^{-1}B_3^+= B_3^-u^{-1}B_3^-\cap \leftmatrix
\CCx &\CC &\CC\\ 0 &0 &\CCx\\ 0 &\CCx &\CC \rightmatrix=
\leftmatrix \CCx &0 &\CCx\\ 0 &0 &\CCx\\ 0 &\CCx &\CC \rightmatrix.
\tag 5.40
$$
We conclude from (5.38), (5.39), and (5.40) that
$$
\P_w= \leftmatrix \CCx &0 &0\\ \CCx &0 &0\\ \CCx &0 &0
\rightmatrix \leftmatrix \CCx &0 &\CCx\\ 0 &0 &0\\ 0 &0 &0
\rightmatrix= \left\{ x\in \leftmatrix \CCx &0 &\CCx\\ \CCx &0
&\CCx\\ \CCx &0 &\CCx \rightmatrix \biggm| \rank(x) = 1 \right\},
$$
as calculated in Example 4.5. \enddefinition

Next, we offer an example in which the union in (5.19) runs over
two disjoint nonempty sets.

\definition{5.14\. Example} Define $\sigma= (y,v,z,u)\in
\Sigma^{3,3}_2$ as follows:
$$
(y,v,z,u)= \left( \leftmat 0&0&1\\ 1&0&0\\ 0&1&0 \rightmat,
\leftmat 1&0&0\\ 0&1&0\\ 0&0&1 \rightmat, \leftmat 1&0&0\\ 0&1&0\\
0&0&1 \rightmat, \leftmat 0&0&1\\ 1&0&0\\ 0&1&0 \rightmat \right).
$$
The nontrivial element of $S^1_2$ can be given as $\tau= \leftmat
0&1&0\\ 1&0&0\\ 0&0&1 \rightmat$, and we observe that $z\tau\le y$
and $v\tau^{-1}\le u$. Next, we calculate that
$$\align
B_3^+yB_3^+ &= \{x\in GL_3\mid x_{31}=0;\ x_{21},x_{32}\ne 0\} \\
B_3^-u^{-1}B_3^- &= \{x\in GL_3\mid x_{13}=0;\ x_{12},x_{23}\ne
0\},
\endalign
$$
and consequently
$$\xalignat2
B_3^+yB_3^+ \cap B_3^-z &= \leftmatrix \CCx &0 &0\\ \CCx &\CCx
&0\\ 0 &\CCx &\CCx \rightmatrix
 &B_3^-u^{-1}B_3^- \cap v^{-1}B_3^+ &= \leftmatrix \CCx &\CCx &0\\
0 &\CCx &\CCx\\ 0 &0 &\CCx \rightmatrix  \\
B_3^+yB_3^+ \cap B_3^-z\tau &= \leftmatrix 0 &\CCx &0\\ \CCx &\CC
&0\\ 0 &\CCx &\CCx \rightmatrix
 &\tau^{-1}B_3^-u^{-1}B_3^- \cap v^{-1}B_3^+ &= \leftmatrix \CCx
&\CC &\CCx\\ 0 &\CCx &0\\ 0 &0 &\CCx \rightmatrix.
\endxalignat
$$
Thus, we find that
$$\multline
\bigl( B_3^+yB_3^+ \cap B_3^-z \bigr) .I^{3,3}_2. \bigl(
B_3^-u^{-1}B_3^- \cap v^{-1}B_3^+ \bigr) \\
 \shoveright{ = \leftmatrix \CCx &0
&0\\ \CCx &\CCx &0\\ 0 &\CCx &0 \rightmatrix  \leftmatrix \CCx
&\CCx &0\\ 0 &\CCx &\CCx\\ 0 &0 &0 \rightmatrix = \left\{ x\in
\leftmatrix \CCx &\CCx &0\\ \CCx &\CC &\CCx\\ 0 &\CCx &\CCx
\rightmatrix \biggm| \rank(x)=2 \right\}  } \\
\shoveleft{ \bigl( B_3^+yB_3^+ \cap B_3^-z\tau \bigr) .I^{3,3}_2.
\bigl( \tau^{-1}B_3^-u^{-1}B_3^- \cap v^{-1}B_3^+ \bigr)  } \\
 = \leftmatrix 0 &\CCx &0\\ \CCx &\CC &0\\ 0 &\CCx &0
\rightmatrix \leftmatrix \CCx &\CC &\CCx\\ 0 &\CCx &0\\ 0 &0 &0
\rightmatrix = \leftmatrix 0 &\CCx &0\\ \CCx &\CC &\CCx\\ 0 &\CCx
&0 \rightmatrix .  
\endmultline
$$
The union of these two disjoint sets equals $\P^t_\sigma$.
\enddefinition

\head 6. Row- and column-echelon forms \endhead

We show that, up to Zariski closure, the $T$-orbits of symplectic
leaves in $\Mmn$ are matrix products of orbits with specific
row- and column-echelon forms. Further, the quasi-affine varieties
of matrices with fixed row-echelon (or column-echelon) forms are
unions of orbits of symplectic leaves of a particularly nice form. 
Throughout the section, overbars will denote Zariski closures within
matrix varieties. As in Section 5, we fix the positive integers $m$ and
$n$ as well as a nonnegative integer $t\le \min\{m,n\}$, and we
concentrate on
$T$-orbits of symplectic leaves within $\Omnt$ (recall (5.3)).

Recall (\S3.9) that the action of $T$ on $\Mmn$ is given by
viewing $T= T_m\times T_n$ and letting $(a,b).x= axb^{-1}$ for
$a\in T_m$, $b\in T_n$, and $x\in \Mmn$. We shall use the
analogous actions of $T_m\times T_t$ and $T_t\times T_n$ on
$M_{m,t}$ and $M_{t,n}$, respectively.

\proclaim{6.1\. Theorem} Let $w\in S_{N}^{\ge\wnwmpair}[t]$
{\rm(}recall {\rm(5.2))}, write $w= \wmn\phi(\sigma)$ for a
unique $\sigma= (y,v,z,u)\in \Sigmnt$ {\rm(}recall {\rm(5.13))},
and set
$$\aligned
\C_{y,z} &= \bigl( B_m^+yB_m^+\cap B_m^-z \bigr).I^{m,t}_t
\subseteq M_{m,t}
\\
\R_{u,v} &= I^{t,n}_t.\bigl(
B_n^-u^{-1}B_n^-\cap v^{-1}B_n^+ \bigr) \subseteq M_{t,n}.
\endaligned  \tag 6.1
$$
Then $\C_{y,z}$ {\rm(}respectively, $\R_{u,v}${\rm)} is a
$(T_m\times T_t)$-orbit {\rm(}respectively, $(T_t\times
T_n)$-orbit{\rm)} of symplectic leaves within
$M_{m,t}$ {\rm(}respectively, $M_{t,n}${\rm)}, and
$$
\C_{y,z}.\R_{u,v} \subseteq \P_w \subseteq
\overline{\C_{y,z}.\R_{u,v}}.
\tag 6.2
$$
In particular, $\overline{\P}_w= \overline{\C_{y,z}.\R_{u,v}}$.
\endproclaim

\demo{Proof} We have $\C_{y,z}.\R_{u,v} \subseteq \P_w$ by Theorem
5.11 (take $\tau=1$ in (5.19)).

Next, viewing $(y,1,z,1)$ as an element of $\Sigma_t^{m,t}$, we
see by Theorem 5.11 that
$$
\P^t_{(y,1,z,1)} = \bigl( B_m^+yB_m^+\cap B_m^-z
\bigr).I_t^{m,t}.(B_t^-\cap B_t^+)= \bigl( B_m^+yB_m^+\cap B_m^-z
\bigr).I_t^{m,t} = \C_{y,z}.
\tag 6.3
$$
Thus, $\C_{y,z}$ is a $(T_m\times
T_t)$-orbit of symplectic
leaves in $M_{m,t}$. Similarly, $\R_{u,v}$ is a $(T_t\times
T_n)$-orbit of symplectic leaves in $M_{t,n}$. In particular, it
follows that their closures
$\overline{\C}_{y,z}$ and
$\overline{\R}_{u,v}$ are Poisson subvarieties of
$M_{m,t}$ and $M_{t,n}$, stable under the respective tori
$T_m\times T_t$ and $T_t\times T_n$.

Let $\mu :M_{m,t}\times M_{t,n}\rightarrow \Mmn$ denote the
morphism given by matrix multiplication, and observe that $\mu$ is
a Poisson map. Since $\overline{\C_{y,z}\times \R_{u,v}}=
\overline{\C}_{y,z}\times \overline{\R}_{u,v}$ (e.g., \cite{\Nor,
Corollary to Theorem 28, p\. 45}), we have
$\mu\bigl( \overline{\C}_{y,z}\times \overline{\R}_{u,v} \bigr)
\subseteq \overline{\C_{y,z}.\R_{u,v}}$. Moreover, as
$\overline{\C}_{y,z}\times \overline{\R}_{u,v}$ is a closed Poisson
subvariety of $M_{m,t}\times M_{t,n}$, the closure $Z$ of
$\mu\bigl( \overline{\C}_{y,z}\times \overline{\R}_{u,v} \bigr)$
is a Poisson subvariety of $\Mmn$, and $Z\subseteq
\overline{\C_{y,z}.\R_{u,v}}$. Note also that if the action of
$T_m\times T_t\times T_t\times T_n$ on $M_{m,t}\times M_{t,n}$ is
restricted to $T_m\times \langle1\rangle\times
\langle1\rangle\times T_n \cong T$, then $\mu$ is $T$-equivariant.
Since $\C_{y,z}\times \R_{u,v}$ is $T$-stable, it follows that
$\mu\bigl( \overline{\C}_{y,z}\times \overline{\R}_{u,v} \bigr)$ is
$T$-stable, and thus $Z$ is a $T$-stable subvariety of $\Mmn$.

Now $\C_{y,z}.\R_{u,v} \subseteq \P_w\cap Z$, so that
$\P_w\cap Z$ is nonempty. Choose $a\in \P_w\cap Z$
and let $\L$ denote the symplectic leaf containing $a$; then
$\P_w= T_m.\L.T_n$. On the other hand, as $Z$ is a $T$-stable
closed Poisson subvariety of $\Mmn$, it is a union of
$T$-orbits of symplectic leaves. Consequently, $T_m.\L.T_n
\subseteq Z$, and therefore
$\P_w \subseteq Z\subseteq \overline{\C_{y,z}.\R_{u,v}}$.
\qed\enddemo

\definition{6.2\. Remark}
Theorem 6.1 can be interpreted as a tensor product result
concerning prime Poisson ideals in coordinate rings, as follows.
First, note that the ideal $P_w$ defining the $T$-stable closed
Poisson subvariety $\overline{\P}_w \subseteq \Mmn$ is a $T$-stable
Poisson ideal in $\O(\Mmn)$, where the action of $T$ on $\O(\Mmn)$
by automorphisms is induced from the $T$-action on $\Mmn$ in the
usual way. It can be shown that $P_w$ is a prime ideal, and that
all $T$-stable prime Poisson ideals of $\O(\Mmn)$ have this form.
Similarly, the defining ideal of $\overline{\C}_{y,z}$
(respectively,
$\overline{\R}_{u,v}$) is a $(T_m\times T_t)$-stable
(respectively, $(T_t\times T_n)$-stable) prime Poisson ideal
$P_{y,z}\subseteq \O(M_{m,t})$ (respectively,
$P_{u,v}\subseteq \O(M_{t,n})$). The statement that
$\overline{\P}_w= \overline{\C_{y,z}.\R_{u,v}}$ is equivalent to
the statement that $P_w$ equals the kernel of the homomorphism
$$
\O(\Mmn) @>{\,\mu^*\,}>> \O(M_{m,t})\otimes\O(M_{t,n})
@>{\,\text{quo}\otimes\text{quo}\,}>> \bigl(
\O(M_{m,t})/P_{y,z} \bigr) \otimes \bigl( \O(M_{t,n})/P_{u,v}
\bigr),
$$
where $\mu^*$ is the comorphism of the 
matrix multiplication map from 
$M_{m,t}\times M_{t,n}$ to $\Mmn$. Consequently,
$$
P_w= (\mu^*)^{-1} \bigl( (P_{y,z}\otimes \O(M_{t,n})) +
(\O(M_{m,t})\otimes P_{u,v}) \bigr).
$$
Such tensor product decompositions were proved to hold for
$T$-stable prime ideals in the generic quantized coordinate ring
of $n\times n$ matrices,
$\O_q(M_n)$, by Goodearl and Lenagan \cite{\GLijm, Theorem 3.5}.
Their development can be used, {\it mutatis mutandis\/} (e.g., by
replacing additive commutators with Poisson brackets), to prove
results of the type above. (While that route only gives
information about closures of $T$-orbits of symplectic leaves in
$\Mmn$, it does have the advantage of working over an arbitrary
base field of characteristic zero.)
\enddefinition

\definition{6.3\. Column-echelon and row-echelon forms} We next
wish to observe that the sets
$\C_{y,z}$ and $\R_{u,v}$ in (6.1) consist of matrices with a
single column-echelon (respectively, row-echelon) form. Note that
to specify a particular column-echelon form for rank $t$ matrices
in $M_{m,t}$, we just need to specify the rows in which the
highest nonzero entries of columns $1,\dots,t$ occur;
column-echelon form requires that the list of these row indices is
strictly increasing.

Let $\Inc^m_t$ denote the set of all strictly increasing sequences
in $\{1,\dots,m\}$ of length $t$, that is,
$$
\Inc^m_t= \{ \bfe= (e_1,\dots,e_t) \in \{1,\dots,m\}^t \mid e_1<
\cdots< e_t\},
$$
and define $\Inc^n_t$ analogously. For $\bfr\in \Inc^m_t$ and
$\bfc\in \Inc^n_t$, define
$$\aligned 
\C^m_\bfr &= \{a\in M_{m,t} \mid a_{r_jj}\ne 0 \text{\ for\ }
j=1,\dots,t \text{\ and\ } a_{ij}=0 \text{\ when\ } i<r_j\} \\
\R^n_\bfc &= \{a\in M_{t,n} \mid a_{ic_i}\ne 0 \text{\ for\ }
i=1,\dots,t \text{\ and\ } a_{ij}=0 \text{\ when\ } j<c_i\}. 
\endaligned \tag 6.4
$$
For example,
$$
\R^6_{(2,4,5)} = \leftmatrix 0 &\CCx &\CC &\CC &\CC &\CC \\
0 &0 &0 &\CCx &\CC &\CC \\ 0 &0 &0 &0 &\CCx &\CC \rightmatrix,
$$
the variety of $3\times6$ matrices in row-echelon form with pivot
columns 2, 4, and 5.

Consider a permutation $z\in S_m^{S^1_t}$. Then $z(1)<
\cdots< z(t)$, whence $\bfr= (z(1),\dots,z(t))$ lies
in $\Inc^m_t$. Given an accompanying $y\in S_m^{S^2_{m-t}}$ with
$z\le y$, we thus see that
$$
\C_{y,z} \subseteq B_m^-z.I^{m,t}_t = \C^m_\bfr.
\tag 6.5
$$
Similarly, if $v\in S_n^{S^1_t}$ and $u\in S_n^{S^2_{n-t}}$ with
$v\le u$, then $\bfc= (v(1),\dots,v(t))\in \Inc^n_t$ and
$$
\R_{u,v} \subseteq I^{t,n}_t.v^{-1}B_n^+ = \R^n_\bfc.
\tag 6.6
$$
The inclusions (6.5) and (6.6) exhibit orbits of symplectic leaves
contained within $\C^m_\bfr$ and $\R^n_\bfc$, indexed by the
following sets.

For $\bfr\in \Inc^m_t$ and $\bfc\in \Inc^n_t$, define
$$\aligned
\Sigma^{m,t}_\bfr &= \{(y,z) \in S_m^{S^2_{m-t}}\times S_m^{S^1_t}
\mid z\le y \text{\ and\  } z(j)= r_j \text{\ for\ } j=1,\dots,t\}
\\
\Sigma^{t,n}_\bfc &= \{(u,v)\in S_n^{S^2_{n-t}}\times
S_n^{S^1_t} \mid v\le u \text{\ and\  } v(i)= c_i \text{\ for\
} i=1,\dots,t\}.
\endaligned \tag 6.7
$$
(There is no ambiguity in this notation in the one overlapping
case, namely when $m=t=n$ and $\bfr=\bfc$, since then $\bfr=\bfc=
(1,2,\dots,t)$ and so $z=v=1$.) We now show that the orbits of
symplectic leaves indexed by
$\Sigma^{m,t}_\bfr$ and $\Sigma^{t,n}_\bfc$ cover $\C^m_\bfr$ and
$\R^n_\bfc$, as follows.
\enddefinition

\proclaim{6.4\. Theorem} If $\bfr \in \Inc^m_t$, then $\C^m_\bfr$
is a disjoint union of $(T_m\times T_t)$-orbits of symplectic
leaves of $M_{m,t}$, indexed by $\Sigma^{m,t}_\bfr$, as follows:
$$
\C^m_\bfr = \bigsqcup_{(y,z) \in \Sigma^{m,t}_\bfr} \bigl(
B_m^+yB_m^+\cap B_m^-z \bigr).I^{m,t}_t.
 \tag 6.8
$$
\endproclaim

\demo{Proof} Recall from (6.3) that $\C_{y,z}= \P^t_{(y,1,z,1)}$
for $(y,z)\in \Sigma^{m,t}_\bfr$, where each $(y,1,z,1)$ is viewed
as an element of $\Sigma_t^{m,t}$. Hence, the sets $\C_{y,z}$ are
$(T_m\times T_t)$-orbits of symplectic leaves of $M_{m,t}$, and
they are pairwise disjoint. Further, (6.5) shows that each such
$\C_{y,z}$ is contained in $\C^m_\bfr$. Thus, $\C^m_\bfr$ contains
the disjoint union displayed in (6.8), and it only remains to
prove equality.

Given $a\in \C^m_\bfr$, note that $\rank(a)=t$. By Theorem 3.9 and
equation (5.3),
$a\in \P_w$ for some $w\in S_{m+t}^{\ge(\wt,\wm)}[t]$. Now apply
Corollary 5.10 and Theorem 5.11 (with $n=t$), to get $w= w^{m+t}_\ci
\phi(\sigma)$ for some $\sigma= (y,v,z,u)\in \Sigma_t^{m,t}$ and
$\P_w= \P^t_\sigma$. Note that since $v\in S_t^{S^1_t}$, it must
be the identity. Write $w=
\leftmatrix w_{11} &w_{12}\\ w_{21} &w_{22} \rightmatrix$ as in
(4.1) (with $n=t$), and observe from (5.18) that $\wm
w_{12}\tr\wt= zI_t^{m,t}u^{-1}$. Hence, Corollary 4.3 implies that
$$
\P_w \subseteq B_m^- zI_t^{m,t}u^{-1} B_t^-.
\tag 6.9
$$

Let $s\in M_{m,t}$ be the (unique) partial permutation matrix such
that $s(j)= r_j$ for $j=1,\dots,t$. Then
$$
\C^m_\bfr= B_m^-s \subseteq B_m^-sB_t^-.
\tag 6.10
$$
      From (6.9) and (6.10), we obtain $B_m^- zI_t^{m,t}u^{-1} B_t^-
\cap B_m^-sB_t^- \ne \varnothing$. Since $zI_t^{m,t}u^{-1}$ and $s$ are
partial permutation matrices, it follows that $zI_t^{m,t}u^{-1} = s$.
(See \S7.1 below for more detail.) In particular,
$su(j)= z(j)$ for $j=1,\dots,t$. Since $s(1)< \cdots< s(t)$ and
$z(1)< \cdots< z(t)$, it follows that $u(1)< \cdots< u(t)$. But
$u$ is a permutation in $S_t$, and so $u=1$. Thus, $\sigma=
(y,1,z,1)$, whence $\P_w= \P^t_\sigma= \C_{y,z}$ by (6.3).
Moreover, $z(j)= s(j)= r_j$ for $j=1,\dots,t$, whence $(y,z) \in
\Sigma^{m,t}_\bfr$.

Therefore $a\in \C_{y,z} \subseteq \C^m_\bfr$, and the proof is
complete.
\qed\enddemo

\proclaim{6.5\. Corollary} If $\bfc \in \Inc^n_t$, then $\R^n_\bfc$
is a disjoint union of $(T_t\times T_n)$-orbits of symplectic
leaves of $M_{t,n}$, indexed by $\Sigma^{t,n}_\bfc$, as follows:
$$
\R^n_\bfc = \bigsqcup_{(u,v)\in \Sigma^{t,n}_\bfc} 
I^{t,n}_t.\bigl( B_n^-u^{-1}B_n^-\cap v^{-1}B_n^+ \bigr).
 \tag 6.11
$$
\endproclaim

\demo{Proof} Note that matrix transposition provides a Poisson
isomorphism from $\C^n_\bfc$ onto $\R^n_\bfc$. Moreover, this map
sends $(T_n\times T_t)$-orbits to $(T_t\times T_n)$-orbits. Note
also that the transpose of a permutation matrix is its inverse.
Therefore, (6.11) follows from (6.8).
\qed\enddemo

\head 7. Generalized double Bruhat cells \endhead

\definition{7.1. Bruhat decompositions in $\Mmn$} In the theory of
reductive algebraic monoids (cf\. \cite{\Ren}), the role of the
Weyl group is taken over by what is now called the {\it Renner
monoid\/}. In the case of the algebraic monoid $M_n$, the Renner
monoid is naturally identified with the monoid of all $n\times n$
{\it partial permutation matrices\/}, that is, $0,1$-matrices with
at most one nonzero entry in each row or column \cite{\Ren, pp\.
326-7}. The Bruhat decomposition of a reductive algebraic monoid
$M$ corresponding to any Borel subgroup $B$ of the group of
invertible elements of $M$ partitions
$M$ into Bruhat cells $BwB$ where $w$ runs through the Renner
monoid \cite{\Ren, Corollary 5.8}. Thus, for any Borel subgroup
$B$ of $\GLn$, the monoid $M_n$ is a disjoint union of Bruhat
cells $BwB$, where $w$ runs through the partial permutation
matrices in $M_n$.

As is well known and easily checked, the above Bruhat
decomposition of $M_n$ holds for the rectangular matrix variety
$\Mmn$ as well. Namely, if $\Stilmn$ denotes the set of partial
permutations in $\Mmn$, then
$$
\Mmn= \bigsqcup_{w\in\Stilmn} B_m^+wB_n^+ =
\bigsqcup_{w\in\Stilmn} B_m^-wB_n^-. \tag 7.1
$$
Consequently, $\Mmn$ is also the disjoint union of the {\it
generalized double Bruhat cells\/}
$$
\Bww= B_m^+w_1B_n^+ \cap B_m^-w_2B_n^-  \tag7.2
$$
for $w_1,w_2\in \Stilmn$. The latter generalize
the standard double Bruhat cells
for $\GLm$, which are obtained when $n=m$
and $w_1, w_2 \in S_m \subset \widetilde{S}^{m,m}$. 

Each double Bruhat cell $\Bww$ is a locally closed
subset of $\Mmn$ because it is an intersection 
of two orbits of algebraic groups.
As is surely well known, $\Bww$ is also smooth
and irreducible, but we could not locate a reference 
in the literature. 
We indicate in Proposition 7.2 and Theorem 7.4 how
these properties follow from our results.
\enddefinition

\proclaim{7.2\. Proposition} Let $w_1,w_2\in \Stilmn$.

{\rm (a)} The generalized double Bruhat cell $\Bww= B_m^+w_1B_n^+
\cap B_m^-w_2B_n^-$ is nonempty if and only if there exists some
$w\in S_N^{\ge\wnwmpair}$ of the form $w= \leftmat * &\wn
w_2\tr\wm \\ w_1 &* \rightmat$. 

{\rm (b)} When $\Bww$ is nonempty, it 
is a smooth locally closed subvariety
of $\Mmn$ which is in addition a $T$-stable complete
Poisson subvariety. In fact,
$$
\Bww= \bigsqcup \big\{ \P_w \bigm| w\in \leftmatrix M_n &\wn
w_2\tr\wm\\ w_1 &M_m \rightmatrix \cap S_N^{\ge\wnwmpair}
\bigr\}.  \tag7.3
$$
\endproclaim

\demo{Proof} 
The smoothness of $\Bww$ in the case when it is nonempty
can be obtained as follows. First, note that the Bruhat cells
$B_m^+w_1B_n^+$ and $B_m^-w_2B_n^-$ are smooth, because they are orbits of
the algebraic groups $B_m^+\times B_n^+$ and $B_m^-\times B_n^-$. Secondly,
$\Bww$ lies  within a single $\GLm \times \GLn$ orbit $\Omnt$ in $\Mmn$
for the action $(g_1, g_2). m = g_1 m g_2^{-1}$, cf\.
\S5.2. Now the intersection of $B^+_m w_1 B^+_n$ and
$B^-_m w_2 B^-_n$ in $\Omnt$ is transversal because
the Lie algebras of $B^+_m$ and $B^-_m$ span $\glm$, 
hence $\Bww$ is smooth.

The rest of the proposition follows from
Corollary 4.3 and Theorem 3.9. 
\qed\enddemo

We will describe the partition (7.3) in terms of the $T$-orbits of
symplectic leaves $\P^t_\sigma$ (recall (5.12)) 
more explicitly in Theorem 7.4 below.
Additional criteria for
$\Bww$ to be nonempty are given in Theorem 7.4 and Corollary 7.7.

For the remainder of this section,

\centerline{\it Fix a nonnegative integer $t\le \min\{m,n\}$,}

\noindent and let $\Stilmn^t$ denote the subset of $\Stilmn$
consisting of partial permutations of rank $t$.

\proclaim{7.3\. Lemma} Every partial permutation in $\Stilmn^t$
can be uniquely represented in the form
$$
yI^{m,n}_tv^{-1}  \tag7.4
$$
for some $y\in S_m^{S^2_{m-t}}$ and $v\in S_n^{S^1_t S^2_{n-t}}$,
and also uniquely in the form
$$
zI^{m,n}_tu^{-1}  \tag7.5
$$
for some $z\in S_m^{S^1_t S^2_{m-t}}$ and $u\in S_n^{S^2_{n-t}}$.
\endproclaim

\demo{Proof} The second statement follows from the first by
noting that $S_m^{S^2_{m-t}}= S_m^{S^1_t S^2_{m-t}} S^1_t$ and
$S_n^{S^2_{n-t}}= S_n^{S^1_t S^2_{n-t}} S^1_t$, and that $\tau
I^{m,n}_t= I^{m,n}_t \tau$ for all $\tau\in S^1_t \subseteq
S_m,\, S_n$.

To prove the first statement, we first show that each element of
$\Stilmn^t$ can be represented in the form (7.4). 
This follows from the facts that
$$\alignat2
\Stilmn^t &= S_m I^{m,n}_t S_n \\
\tau_1 I^{m,n}_t &= I^{m,n}_t \tau_2= I^{m,n}_t
&&\qquad\text{for all\ } \tau_1\in S^2_{m-t} \text{\ and\ }
\tau_2\in S^2_{n-t} \\
\tau I^{m,n}_t &= I^{m,n}_t \tau  &&\qquad\text{for all\ } \tau\in
S^1_t \subseteq S_m,\, S_n. \endalignat
$$

The lemma will now follow if we prove that the sets $\Stilmn^t$
and $S_m^{S^2_{m-t}}\times S_n^{S^1_t S^2_{n-t}}$ have the same
number of elements. The cardinality of the second set is
$\frac{m!}{(m-t)!} \frac{n!}{t!(n-t)!}= t!{m\choose t}{n\choose
t}$ because each coset in $S_m/S^2_{m-t}$ or $S_n/S^1_tS^2_{n-t}$
has a unique minimal length representative. Observe that a
partial permutation $w\in \Stilmn^t$ is uniquely defined by
prescribing its domain $\dom w$, range $\rng w$ (both of
cardinality $t$), and a bijective mapping from $\dom w$ to $\rng
w$. Therefore the cardinality of $\Stilmn^t$ is ${m\choose
t}{n\choose t}t!$. \qed\enddemo

\proclaim{7.4\. Theorem} Fix two partial permutations $w_1,w_2\in
\Stilmn^t$ with {\rm(}unique{\rm)} decompositions
$$\xalignat2
 w_1 &= yI^{m,n}_tv^{-1}  &w_2 &= zI^{m,n}_tu^{-1}  
\tag7.6  \endxalignat
$$
for some $y\in S_m^{S^2_{m-t}}$, $v\in S_n^{S^1_t S^2_{n-t}}$,
$z\in S_m^{S^1_t S^2_{m-t}}$, and $u\in S_n^{S^2_{n-t}}$
{\rm(}cf\. Lemma {\rm 7.3}{\rm)}. Then the following hold.

{\rm (a)} The generalized double Bruhat cell $\Bww= B_m^+w_1B_n^+
\cap B_m^-w_2B_n^-$ is nonempty if and only if $z\le y$ and
$v\le u$. 

If $z\le y$ and $v\le u$, then:

{\rm (b)} The partition of $\Bww$ into $T$-orbits of symplectic
leaves is given by
$$
\Bww= \bigsqcup\, \biggl\{ \P^t_{(y,v\tau_2,z\tau_1,u)} \biggm|
\matrix \tau_1\in S^2_{m-t}\subseteq S_m,\ z\tau_1\le y \\
\tau_2\in S^2_{n-t}\subseteq S_n,\ v\tau_2\le u \endmatrix \biggr\}. 
\tag7.7
$$

{\rm (c)} The $T$-orbit of symplectic leaves
$\P^t_{(y,v,z,u)}$ is an open and dense subset of
$\Bww$.

{\rm (d)} $\Bww$ is a smooth irreducible locally closed
subvariety of $\Mmn$.  \endproclaim

For the proof of Theorem 7.4 we will need two lemmas. Recall the set
$\Sigmnt$ from (5.13).

\proclaim{7.5\. Lemma}
For any $\sigma = (y,v,z,u) \in \Sigma^{m,n}_t$, we have
$$
\P_\sigma^t \subseteq B^+ y I^{m,n}_t v^{-1} B^+ 
\cap B^- z I^{m,n}_t u^{-1} B^-.
$$
\endproclaim

\demo{Proof}
The lemma follows from Lemma 5.3 because, for $r_1$, $r_2$ as in (5.11),
$$
r_1 I^{m,n}_t r_2^{-1} = b^+_1 y I^{m,n}_t v^{-1} (b^+_2)^{-1} \in B^+
y I^{m,n}_t v^{-1} B^+
$$
and
$$\align
r_1 I^{m,n}_t r_2^{-1} &= b^-_3 z \leftmatrix a & b \\ 0 & I_{m-t}
\rightmatrix I^{m,n}_t  \leftmatrix a & 0 \\ c & I_{n-t}
\rightmatrix^{-1} u^{-1} (b^-_4)^{-1} \\
  &=  b^-_3 z I^{m,n}_t u^{-1} (b^-_4)^{-1} \in B^- z I^{m,n}_t u^{-1}
B^-. \qquad\square
\endalign
$$
\enddemo

\proclaim{7.6\. Lemma}
Set 
$$
\til{\Sigma}^{m,n}_t =
\{ (y, v_0, z_0, u) \in 
S_m^{S_{m-t}^2} \times S_n^{S_t^1 S_{n-t}^2}
\times S_m^{S_t^1 S_{m-t}^2} \times S_n^{S_{n-t}^2}
\mid z_0 \leq y,\, v_0 \leq u \}.
$$
Then
$$
\Sigma_t^{m,n} = \biggl\{ (y, v_0 \tau_2, z_0 \tau_1, u) \biggm| 
\matrix (y, v_0, z_0, u) \in \til{\Sigma}^{m,n}_t,\ \tau_1\in
S^2_{m-t}\subseteq S_m, \\  \tau_2\in S^2_{n-t}\subseteq S_n,\
z_0\tau_1\le y,\  v_0\tau_2\le u \endmatrix \biggr\}.
\tag 7.8
$$
\endproclaim

\demo{Proof}
It is clear that every element of $\Sigma^{m,n}_t$
has the form $(y, v_0 \tau_2, z_0 \tau_1, u)$ for some
$$
(y, v_0, z_0, u) \in S_m^{S_{m-t}^2} \times S_n^{S_t^1 S_{n-t}^2}
\times S_m^{S_t^1 S_{m-t}^2} \times S_n^{S_{n-t}^2}
$$
and some
$$\xalignat2
\tau_1 &\in S_{m-t}^2 \subseteq S_m 
 &\tau_2 &\in S_{n-t}^2 \subseteq S_n
\endxalignat
$$
such that 
$$\xalignat2
z_0 \tau_1 &\leq y  &v_0 \tau_2 &\leq u. 
\endxalignat
$$
But $z_0 \in S_m^{S_t^1 S_{m-t}^2}$ and
$\tau_1 \in S_{m-t}^2$ imply that $z_0 \leq z_0 \tau_1$
and therefore $z_0 \leq y$. Analogously, one obtains
that $v_0 \leq v_0 \tau_2$ and as a consequence of it
$v_0 \leq u$. Therefore 
$$
(y, v_0, z_0, u) \in \til{\Sigma}^{m,n}_t.
$$

This proves that $\Sigma^{m,n}_t$ is contained in the 
set on the right hand side of (7.8). The opposite inclusion is 
straightforward. 
\qed
\enddemo

\demo{Proof of Theorem {\rm7.4}}
Combining Lemma 7.6 and Corollary 5.12, one obtains
$$
\O^{m,n}_t = 
\bigsqcup_{(y, v_0, z_0, u) \in \til{\Sigma}^{m,n}_t}
\bigsqcup\, \biggl\{ \P^t_{(y, v_0 \tau_2, z_0 \tau_1, u)} \biggm|
\matrix \tau_1\in S^2_{m-t}\subseteq S_m,\ z_0\tau_1\le y \\
\tau_2\in S^2_{n-t}\subseteq S_n,\ v_0\tau_2\le u \endmatrix \biggr\}.
\tag 7.9
$$
At the same time,
$$
\O^{m,n}_t = \bigsqcup_{w_1, w_2 \in \Stilmn^t}
B^+ w_1 B^+ \cap B^- w_2 B^-.
\tag 7.10
$$
      From Lemma 7.5, for each $T$-orbit of leaves on the right hand side 
of (7.9) one derives:
$$
\P^t_{(y, v_0 \tau_2, z_0 \tau_1, u)} \subseteq
B^+ y I^{m,n}_t v_0^{-1} B^+ \cap 
B^- z_0 I^{m,n}_t u^{-1} B^-.
$$  
Comparing (7.9) and (7.10) now proves at once
parts (a) and (b).

(c) Because of (7.7), it suffices to show that
$$
\P^t_{(y,v\tau_2,z\tau_1,u)} \subseteq \overline{
\P^t_{(y,v,z,u)} }  \tag7.11
$$
for $\tau_1\in S^2_{m-t}$ and $\tau_2\in S^2_{n-t}$ such
that $z\tau_1\le y$ and
$v\tau_2\le u$. Fix such $\tau_1$, $\tau_2$, recall the bijection
$\Sigmnt \rightarrow S_N^{\ge\wnwmpair}[t]$ given in Corollary 5.10, and
set
$$\align
 \wbar &= \wN\phi(y,v,z,u)=  \left[\matrix
\wn uJ^n_tv^{-1} &\wn uI^{n,m}_tz^{-1}\wm\\ yI^{m,n}_tv^{-1}
&yJ^m_tz^{-1}\wm \endmatrix\right]  \\
w &= \wN\phi(y,v\tau_2,z\tau_1,u)=  \left[\matrix
\wn uJ^n_t\tau_2^{-1}v^{-1} &\wn uI^{n,m}_tz^{-1}\wm\\
yI^{m,n}_tv^{-1} &yJ^m_t\tau_1^{-1}z^{-1}\wm \endmatrix\right].
\endalign
$$
By Theorems 5.11 and 3.13, (7.11) is equivalent to $w\le \wbar$.

First, note that $w(j)= \wbar(j)= n+w_1(j)$ for $j\in v\bigl(
\{1,\dots,t\} \bigr)$. Now $w$ and $\wbar$ both map $v\bigl(
\{t+1,\dots,n\} \bigr)$ bijectively onto $\wn u\bigl(
\{t+1,\dots,n\} \bigr)$, and for $\wbar$ this restriction is
order-reversing because $u,v\in S_n^{S^2_{n-t}}$. It follows that
$w\bigl( \{1,\dots,j\} \bigr)\le \wbar\bigl( \{1,\dots,j\}
\bigr)$ for $j=1,\dots,n$. Similarly, $w$ and $\wbar$ agree on
$n+\wm z\bigl( \{1,\dots,t\} \bigr)$, and the restriction of
$\wbar$ to $n+\wm z\bigl( \{t+1,\dots,m\} \bigr)$ is
order-reversing, from which we conclude that $w\bigl(
\{1,\dots,j\} \bigr)\le \wbar\bigl( \{1,\dots,j\}
\bigr)$ for $j=n+1,\dots,N$. Therefore $w\le \wbar$, as required.

(d) The irreducibility of $\Bww$ follows from part
(c) since $\P^t_{(y,v,z,u)}$ is irreducible by Theorem 3.9. 
\qed\enddemo

\proclaim{7.7\. Corollary} For partial permutations $w_1,w_2\in
\Stilmn^t$, the generalized double Bruhat cell $\Bww= B_m^+w_1B_n^+
\cap B_m^-w_2B_n^-$ is nonempty if and only if
$$
\dom(w_1) \le \dom(w_2) \qquad\text{and}\qquad \rng(w_1) \ge
\rng(w_2)  \tag7.12
$$
{\rm(}recall \S{\rm 3.11)}.
\endproclaim

\demo{Proof} Let $w_1 = yI^{m,n}_tv^{-1}$ and $w_2 =
zI^{m,n}_tu^{-1}$ for $y$, $v$, $z$, and $u$ as in Theorem 7.4.

If $\Bww$ is nonempty, then by the theorem, $z\le y$ and $v\le
u$. Hence,
$$\aligned
 \dom(w_1) &= v\bigl( \{1,\dots,t\} \bigr) \le u\bigl(
\{1,\dots,t\} \bigr) = \dom(w_2)  \\
\rng(w_1) &= y\bigl( \{1,\dots,t\} \bigr) \ge z\bigl(
\{1,\dots,t\} \bigr) = \rng(w_2).  \endaligned  \tag7.13
$$
Conversely, assume that $\dom(w_1) \le \dom(w_2)$ and $\rng(w_1)
\ge \rng(w_2)$, so that (7.13) holds. It follows, as shown in the
proof of Lemma 5.7, that $v\bigl( \{t+1,\dots,n\} \bigr) \ge
u\bigl( \{t+1,\dots,n\} \bigr)$. Since $u,v\in S_n^{S^2_{n-t}}$,
we obtain $v(j)\ge u(j)$ for $j=t+1,\dots,n$. But then, since
$v\in S_n^{S^1_t}$, Lemma 5.7(b) implies that $v\le u$. Similarly,
$z\le y$. \qed\enddemo

\head Appendix A. Double coset representatives\endhead

\definition{A.1}
Let $G$ be a complex reductive algebraic group with fixed 
positive/negative Borel subgroups $B^\pm$ and maximal torus
$T= B^+ \cap B^-$. Fix a parabolic subgroup $P$ of  
$G$, containing a Borel subgroup $B \supset T$ of $G$ 
with the property that for each simple factor $F$ of $G$,  either 
$B \cap F = B^+ \cap F$ or $B \cap F = B^- \cap F$. 

Denote by $L_0$ the Levi factor of $P$ containing
$T$ and by $N$ the unipotent radical of $P$. So, we have the the
Levi decomposition $P \cong L_0 \ltimes N$. Denote by $\ol{N}$ the
unipotent subgroup of $G$ dual to $N$. 

We will assume that $L_0$ is decomposed as a product of two
reductive  subgroups 
$$
L_0 = L_1 \times L_2  \tag\roman{A.1}
$$  
such that there is an isomorphism
$$ 
\Theta : L_1 @>{\cong}>> L_2  \tag\roman{A.2}
$$ 
with the property that for every simple factor 
$F_1$ of $L_1$,
$$
\Theta(F_1 \cap B^\pm) = F_2 \cap B^+  \tag\roman{A.3}
$$
for some simple factor $F_2$ of $L_2$
and an appropriate choice of the sign.

Denote the Weyl group of $G$ by $W$ and the 
Weyl groups of $L_i$ $(i=0,1,2)$ by $W_i$, 
considered as subgroups of 
$W$. Clearly $W_0 = W_1 \times W_2$. Denote the composition of the
projections $P @>>>P/N \cong L_0$ and
$L_0 @>>> L_i$ $(i=1,2)$ by $ \pi_i : P @>>> L_i$.

Finally, define the following subgroup of $P$:
$$ 
R = \{ p \in P \mid \Theta \pi_1(p) = \pi_2(p) \}.  \tag\roman{A.4}
$$

In this Appendix we give a classification of all
$(B^+, R)$ double cosets of $G$. Recall that $W^{W_i}$ denotes the
set of  (unique) minimal length representatives of cosets from
$W/W_i$, see \cite{\Car, Proposition 2.3.3} for details.
For an  element $w \in W$, we will denote by $\dot{w}$
a representative of it in the  normalizer of $T$ in $G$.
\enddefinition

\proclaim{Theorem} 
In the above setting,
every $(B^+, R)$ double coset of $G$ is of the
form
$$
B^+ \dot{w} R, \quad \text{for some} \quad w \in W^{W_2}.
$$  
For distinct $w \in W^{W_2}$, 
the above double cosets are distinct.
\endproclaim 

Let us note that in the case when $L_1$ and $L_2$ have more than 
one simple factor, it is possible to obtain $R$ as a subgroup 
of $P$ in several different ways by changing $L_1$ and $L_2$. 
In such a case, Theorem A.1 produces different sets of representatives 
for the $(B^+, R)$ double cosets of $G$. As is clear from Lemma 3.8,
sometimes one of these sets has better properties than the 
others.

For the proof of Theorem A.1, we will need the following lemma.

\proclaim{A.2. Lemma} {\rm (a) (Bruhat Lemma)} All $(B^+, P)$
double cosets in $G$ are uniquely pa\-ram\-e\-trized by $W^{W_0}$, by 
$v \in W^{W_0} \mapsto B^+ \dot{v} P$.

{\rm (b)} For any $v \in W^{W_0}$,
$$
B^+ \dot{v} =  \dot{v} \ol{N}_v B_0^+ N_v
$$ 
where $B_0^+ = B^+ \cap L_0$ and  
$$\xalignat2
N_v &= N \cap Ad^{-1}_{\dot{v}}(B^+)   
 &\ol{N}_v  &= \ol{N} \cap Ad^{-1}_{\dot{v}}(B^+).  \endxalignat
$$  

{\rm (c)} There is a bijection of sets
$$\xalignat2
W^{W_0}\times W_1 &\rightarrow W^{W_2},  &(v,u) &\longmapsto vy. 
\endxalignat
$$ 

{\rm (d)} 
Set $Q=R \cap L_0 = \{ l_1 \Theta(l_1) \mid l_1 \in L_1 \}$. 
All $(B_0^+, Q)$ double cosets of $L_0$ are uniquely parametrized  by
$W_1$, by $w_1 \mapsto B_0^+ \dot{w}_1 Q$.
\endproclaim

\demo{Proof} Part (a) is well known.

Part (b) follows from the well known
description of minimal  length representatives:
$$
W^{W_0} = \{ w \in W \mid w(\alpha) \; \text{is a positive root} 
\; \text{for any positive root} \; \alpha \; \text{of} \; L_0 \}.
$$ 
See, e.g., \cite{\Car, Proposition 2.3.3}.
Part (c) is a consequence of $W_0 = W_1 \times W_2$.

To prove part (d), we first show that it suffices to establish (d) in
the case that
$$
\Theta(L_1\cap B^+) = L_2\cap B^+.  \tag\roman{A.5}
$$
For each simple factor $F_1$ of $L_1$, the assumption (A.3) can be
written in the form
$$
\Theta \Ad_{\dot u} (F_1\cap B^+) = F_2\cap B^+,
$$
where $u$ is either the identity or the longest element of the Weyl
group of $F_1$. Hence, there exists an element $u_1\in W_1$ such that
$u_1^2=1$ and
$$
\Theta \Ad_{\dot u_1}(L_1\cap B^+)= L_2\cap B^+.
$$
The map $\Thetatil= \Theta\circ \Ad_{\dot u_1}|_{L_1}$ is an
isomorphism of $L_1$ onto $L_2$, and the subgroup $\Qtil$ of $G$
obtained by changing $\Theta$ to $\Thetatil$ in the definition of $Q$
can be written as
$$
\Qtil= \{ \ltil_1 \Thetatil(\ltil_1) \mid \ltil_1\in \Ad_{\dot
u_1}(L_1) \}= \{ \Ad_{\dot u_1}(l_1) \Theta(l_1) \mid l_1\in L_1 \}=
\Ad_{\dot u_1}(Q).
$$
If (d) holds for $\Thetatil$, then, since $W_1=W_1u_1$, we may express
the result as
$$
L_0= \bigsqcup_{w_1\in W_1} B_0^+ \dot{w}_1 \dot{u}_1 \Qtil,$$
and consequently
$$
L_0= L_0 \dot{u}_1= \bigsqcup_{w_1\in W_1} B_0^+ \dot{w}_1 Q,
$$
as desired. Thus, we may assume (A.5), as claimed.

Recall the fact that if $F_1$, $F_2$
are subgroups of  a group $C$, then the set of $(F_1 \times F_2,\,
\Delta(C))$ double cosets of $C \times C$ (where $\Delta(C) \subseteq
C\times C$ denotes the diagonal copy of $C$) is in one to one
correspondence with the set of 
$(F_1, F_2)$ double cosets of $C$, by 
$(F_1 \times F_2) (y_1, y_2) \Delta(C) \mapsto F_1 y_1 y_2^{-1} F_2$.
If we identify $L_0$ with $L_1\times L_1$ via $\Theta$, then $Q$ is
identified with $\Delta(L_1)$, and because of (A.5), $B_0^+$ is
identified with $B_1^+\times B_1^+$, where $B_1^+= L_1\cap B^+$. Since
the $(B_1^+,B_1^+)$ double cosets of $L_1$ are uniquely parametrized by
$W_1$, the $(B_1^+\times B_1^+,\, \Delta(L_1))$ double cosets of
$L_1\times L_1$ are uniquely parametrized by $W_1\times \{1\}$, and
part (d) follows.
\qed\enddemo

\demo{Proof of Theorem {\rm A.1}} Since $P = L_0 R$, the Bruhat lemma
implies that every 
$(B^+, R)$ double coset of $G$ is of the form 
$B^+ \dot{v} l_0 R$ for
some 
$v \in W^{W_0}$ and $l_0 \in L_0$. In addition, the Bruhat lemma
also  implies that if $B^+ \dot{v} l_0 R= B^+ \dot{v}' l_0' R$ 
for some $v, v' \in W^{W_0}$ and $l_0, l_0' \in L_0$, then $v'=v$.

      From the facts that $R= QN= NQ$, that $L_0$ normalizes 
$N$, and $N_v \subseteq N$, we get
$$
N_v l_0 R= N_v l_0 Q N = l_0 Q N.
$$
Thus, from part (b) of the above lemma, we have
$$ 
B^+ \dot{v} l_0 R = \dot{v} \ol{N}_v ( B^+_0 l_0 Q) N.
$$ 
Since $\ol{N}_v \subseteq \ol{N}$ and $\ol{N} L_0 N$ is the Cartesian
product  of the subsets $\ol{N}$, $L_0$ and $N$ of $G$, we get that 
for $v \in W^{W_0}$ and $l_0, l_0' \in L_0$, 
$$ 
B^+ \dot{v} l_0 R = B^+ \dot{v} l_0' R \; 
\iff \;  B^+_0 l_0 Q = B^+_0 l_0' Q.
$$ 
Part (d) of the lemma now implies that all $(B^+, R)$ double
cosets of $G$ are uniquely parametrized by $W^{W_0} \times W_1$, by 
$(v, u) \mapsto B^+ \dot{v} \dot{u} R$. The theorem, 
finally, follows from part (c) of the lemma.
\qed\enddemo

\head Acknowledgements\endhead

We thank Allen Knutson, Jiang-Hua Lu, Jon McCammond, James McKernan,
Nicolai Reshetikhin, Manfred Schocker and Alan Weinstein for helpful
discussions and correspondence.

\enddocument